%% file: poly.tex
\documentclass[reqno,10pt]{amsart}
\usepackage{tabularx}
\usepackage{cancel}
\usepackage{multirow}
\usepackage{array}
\usepackage{wasysym}
\usepackage{arydshln}
\usepackage{hhline}
\usepackage{dashbox}
\usepackage{todonotes}
\usepackage[all,2cell]{xy}
\usepackage[dvipsnames]{xcolor}
\usepackage{graphicx, amsmath, amssymb, amsthm}
\usepackage{newclude}
\usepackage{adjustbox}
\usepackage{mathrsfs}
\usepackage{sseq}
\usepackage[hidelinks,backref]{hyperref}
\usepackage{enumitem}
\usepackage[margin=1in]{geometry} 

\usepackage{fancyhdr}
\input{macros.tex}

\makeatletter
\providecommand\@dotsep{5}
\renewcommand{\listoftodos}[1][\@todonotes@todolistname]{%
  \@starttoc{tdo}{#1}}
\makeatother

\SelectTips{cm}{}
\UseAllTwocells
\NoCompileMatrices

\title{\tops{$\RO(G)$}{RO(G)}-graded norms for de Rham-Witt forms}

\author{Yuri~J.~F. Sulyma}
\email{yuri.sulyma@protonmail.com}

\begin{document}

\begin{abstract}
  We make two case studies of $\mathrm{RO}(G)$-graded homotopy in connection with $p$-adic cohomology. First, we identify the $\mathrm{RO}(\mathbf T)$-graded $\mathrm{TC}^-$ and $\mathrm{TP}$ of quasiregular semiperfectoid rings as multi-Rees algebras for the Nygaard and $\phi^\bullet(I)$-adic filtrations. We additionally compute $E^{h\mathbf T}_\bigstar$ and $E^{t\mathbf T}_\bigstar$ when $E$ is a general even $\mathbf E_\infty$-ring with $\mathbf T$-action, and compute $(\mathrm{TR}^n)^{h\mathbf T/C_{p^n}}_\bigstar$ in terms of iterated Nygaard filtrations. Next, we compute the $\mathrm{RO}(\mathbf Z_p)^\land$-graded $\mathrm{TR}$ of smooth algebras over a perfectoid base, providing a pushout of the works of Hesselholt and Angeltveit-Gerhardt along the base case $\mathbf F_p$. Using this, we extend Angeltveit's Norm map of Witt vectors to the de Rham-Witt complex, which participates in a ``Witt coordinates'' decomposition of the $p$-primary Adams operations.
\end{abstract}

\maketitle
\tableofcontents


\begin{flushright}
  \begin{minipage}{0.465\textwidth}
    \itshape
    \raggedleft
     ``I remember a psychiatrist once telling me that I gamble in order to escape the reality of life, and I told him that’s why everyone does everything.''

    \medskip
    \normalfont
    Norm Macdonald, Canadian mathematician 
    
    \vspace{0.2cm} 
    \hrule
  \end{minipage}
\end{flushright}

\include*{intro}
\include{notation}
\include{qrsp}
\include{dRW}

\include{norms}

\bibliographystyle{amsalpha}
\bibliography{../bibliography}

\end{document}

%% file: macros.tex
\newcommand{\<}{\left\langle}
\let\>=\undefined
\newcommand{\>}{\right\rangle}

\newcommand{\pos}[1]{(\!(#1)\!)}

\DeclareMathOperator{\gr}{gr}

\DeclareMathOperator{\Spf}{Spf}

\newcommand{\bigoplusp}{\mathop{\widehat\bigoplus}\limits}

\newcommand{\can}{\mathrm{can}}
\newcommand{\rmW}{\mathrm W}
\newcommand{\dRW}[1][]{\mathrm W_{#1}\Omega}
\newcommand{\DRW}[1][]{\mathbf W_{#1}\Omega}
\newcommand{\AOmega}{A\Omega}
\newcommand{\qOmega}{\mathrm{q}\Omega}

\let\ifdraft=\iftrue

\newcommand{\mot}{\mathrm{mot}}
\newcommand{\triv}{\mathrm{triv}}
\renewcommand{\th}{\mathrm{th}}

\usepackage{relsize}
\usepackage[bbgreekl]{mathbbol}
\DeclareSymbolFontAlphabet{\mathbb}{AMSb} 
\DeclareSymbolFontAlphabet{\mathbbl}{bbold}

\newcommand{\prism}{{\mathlarger{\mathbbl{\Delta}}}}

\newcommand{\ds}{\displaystyle}

\newcommand{\id}{\mathrm{id}}

\renewcommand{\top}{\mathrm{top}}
\newcommand{\alg}{\mathrm{alg}}

\renewcommand{\a}{\m a}
\renewcommand{\u}{\m u}

\newcommand{\defeq}{\mathrel{:=}}

\newcommand{\tops}{\texorpdfstring}

\newcommand{\red}[1]{{\color{red}{#1}}}

\DeclareMathOperator{\fib}{fib}

\DeclareMathOperator{\dlog}{dlog}

\DeclareMathOperator{\ind}{ind}
\DeclareMathOperator{\fil}{fil}
\newcommand{\Fil}{\mathrm{Fil}}
\DeclareMathOperator{\res}{res}

\newcommand{\HKR}{\mathrm{HKR}}

\DeclareMathOperator{\Sym}{Sym}

\DeclareMathOperator{\RO}{RO}
\DeclareMathOperator{\RU}{RU}
\renewcommand{\~}{\widetilde}
\newcommand{\hsl}[1]{[\![#1]\!]_\lambda}

\newcommand{\lf}{\left\lfloor}
\newcommand{\rf}{\right\rfloor}
\newcommand{\lc}{\left\lceil}
\newcommand{\rc}{\right\rceil}

\newdir{ >}{{}*!/-10pt/@{>}}
\newcommand{\pullback}{\ar@{}[dr]|<<{\mbox{\Huge$\lrcorner$}}}
\newcommand{\pushout}{\ar@<2pt>@{}[ul]|<{\mbox{\Huge$\ulcorner$}}}

\newcommand{\Mod}{\operatorname{Mod}}

\newcommand{\Sp}{\mathrm{Sp}}

\newcommand{\C}{\mathbf C}

\newcommand{\E}{\mathbf E}
\newcommand{\F}{\mathbf F}

\newcommand{\G}{\mathbf G}
\renewcommand{\H}{\mathrm H}

\renewcommand{\L}{\mathbf L}
\newcommand{\N}{\mathbf N}

\renewcommand{\O}{\mathcal O}

\newcommand{\Q}{\mathbf Q}
\newcommand{\R}{\mathbf R}
\let\sec=\S
\renewcommand{\S}{\mathbf S}
\newcommand{\T}{\mathbf T}

\newcommand{\V}{\mathbf V}
\newcommand{\W}{\mathbf W}
\newcommand{\Z}{\mathbf Z}

\newcommand{\K}{{\mathrm{K}}}

\newcommand{\HH}{\mathrm{HH}}
\newcommand{\HC}{\mathrm{HC}}
\newcommand{\HP}{\mathrm{HP}}

\newcommand{\THH}{\mathrm{THH}}
\newcommand{\TR}{\mathrm{TR}}
\newcommand{\TC}{\mathrm{TC}}
\newcommand{\TCmin}{\mathrm{TC}^-}
\newcommand{\TP}{\mathrm{TP}}
\newcommand{\TF}{\mathrm{TF}}

\newcommand{\bigTR}{\mathbf{TR}}

\newcommand{\Nyg}{\mathcal N}
\newcommand{\m}{\underline}
\newcommand{\mpi}{\underline\pi}
\newcommand{\rog}{\bigstar}

\newcommand{\Ainf}{{ \mathrm A_{\mathrm{inf}} }}
\newcommand{\SW}[1]{\S_{\rmW(#1)}}
\newcommand{\SAinf}[1]{\S_{\Ainf(#1)}}
\newcommand{\rel}{\hspace{-1pt}\mathbin|\hspace{-1pt}}

\DeclareMathOperator{\im}{im}

\newcommand{\tnsr}{\otimes}

\DeclareMathOperator{\Hom}{Hom}

\DeclareMathOperator{\Tor}{Tor}
\DeclareMathOperator{\Ext}{Ext}

\newcommand{\Sol}{\mathrm{Sol}}

\newcommand{\morph}{\mathop{\longrightarrow}\limits}
\newcommand{\xra}{\xrightarrow}

\newcommand{\plainlim}{\operatorname{lim}}
\renewcommand{\lim}{\mathop{\operatorname*{lim}\limits_{\longleftarrow}}\limits}

\newcommand{\epi}{\twoheadrightarrow}
\newcommand{\mono}{\hookrightarrow}

\newcommand{\iso}{\overset\sim\longrightarrow}

\newcommand{\isom}{\cong}

\newtheorem{theorem}{Theorem}[section]
\newtheorem*{theorem*}{Theorem}

\newtheorem{lemma}[theorem]{Lemma}
\newtheorem{corollary}[theorem]{Corollary}
\newtheorem{proposition}[theorem]{Proposition}

\theoremstyle{definition}
\newtheorem{definition}[theorem]{Definition}
\newtheorem{remark}[theorem]{Remark}

\newtheorem{method}[theorem]{Method}
\newtheorem{proofsketch}[theorem]{Proof Sketch}
\newtheorem{conjecture}[theorem]{Conjecture}
\newtheorem{example}[theorem]{Example}
\newtheorem{warning}[theorem]{Warning}

\newtheorem{motivation}[theorem]{Motivation}

\newtheorem{notation}[theorem]{Notation}

\newcommand{\fade}[2][20]{{\color{black!#1}#2}}

\newcommand{\prismcpl}{\hat\prism}

\newcommand{\mup}[1]{\ar@/_1em/[u]_-{#1}}
\newcommand{\mdown}[1]{\ar@/_1em/[d]_-{#1}}
\newcommand{\mcup}[2][]{\ar@[#1]@/_1em/[u]_-{\color{#1}#2}}
\newcommand{\mcdown}[2][]{\ar@[#1]@/_1em/[d]_-{\color{#1}#2}}

\newcommand{\Lewis}[4]{\xymatrix{
  #1 \mdown{#3}\\#2 \mup{#4}
}}

\newcommand{\Lewiss}[7]{\xymatrix{
    #1 \mdown{#4}\\ 
    #2 \mdown{#5} \mup{#6}\\ 
    #3 \mup{#7}
}}

%% file: intro.tex
%
\section{Introduction}
\label{sec:intro}

This paper establishes some foundational computations for studying $p$-adic cohomology through the lens of genuine equivariant homotopy theory. Recall that the interface between these is topological Hochschild homology (here viewed as a proper-genuine $\T$-spectrum), and its various types of fixed point spectra. For nice input rings, the homotopy groups of these are given by important arithmetic invariants, often with additional \emph{filtrations} and \emph{twists}.

\begin{example}
If $R$ is a quasiregular semiperfectoid ring, the ($p$-completed) $\TCmin$ and $\TP$ of $R$ are related to its (Nygaard-filtered, Nygaard-completed, Breuil-Kisin twisted, absolute) prismatic cohomology via
\begin{align*}
  \TCmin_{2i}(R;\Z_p) &= \Nyg^{\ge i}\hat\prism_R\{i\}\\
  \TP_{2i}(R;\Z_p) &= \hat\prism_R\{i\}
\end{align*}
(and the odd groups vanish) by work of Bhatt-Morrow-Scholze \cite{BMS2}. In the $q$-de Rham case, the twist can be understood as tensoring with a free module on the symbol $(p^\infty)_q = \frac{q^{p^\infty}-1}{q-1} = \frac{\log_\prism(q^p)}{q-1}$ (see Notation \ref{not:canonical} and \cite[Proposition 2.6.1]{APC}). Antieau-Krause-Nikolaus have announced a generalization of this result in terms of their prismatic cohomology relative to a $\delta$-ring $A$ \cite{RelativePrismatic} whenever $A$ descends along $\S_p\to\Z_p$.
\end{example}

\begin{example}
  When $S$ is $p$-completely smooth over any perfectoid ring $R$ (i.e.\ the $p$-completion of a smooth $R$-algebra), Hesselholt \cite{LarsPTyp} showed that $\TR^n(S;\Z_p)$ is given by the $p$-complete relative de Rham-Witt complex \cite{LangerZinkDRW} adjoined with a polynomial generator in degree 2:
  \begin{equation}
    \label{eq:hesselholt}
    \TR^n_*(S;\Z_p) \isom \dRW[n]^*_{\Spf S/\Spf R}[\sigma_n]\qquad |\sigma_n|=2
  \end{equation}
  Here the $p$-complete de Rham-Witt complex is
  \[ \dRW[n]^*_{\Spf S/\Spf R} \defeq \lim_s \dRW[n]^*_{(S/p^s)/(R/p^s)} \]
  as in \cite[Definition 10.11]{BMS1}; we will subsequently abuse notation and simply write this as $\dRW[n]^*_{S/R}$.
  The identification \eqref{eq:hesselholt} can be made canonical by interpreting $\sigma_n$ as tensoring with the symbol $\xi_{n+1}$ ($= (p^{n+1})_{q^{1/p^{n+1}}}$ in the $q$-de Rham case), see Notation \ref{not:canonical}.

  In particular, the Restriction maps on the twists get multiplied by powers of $\xi_{n+1}/\xi_n=\phi^{-n}(\xi)$, so there are potentially $\plainlim^1$ terms when we take the limit over this result to understand $\TR(S;\Z_p)$. Fortunately, these vanish in the main cases of interest: if $R=k$ is a perfect $\F_p$-algebra, then
  \[ \TR_*(S) = \dRW^*_S, \]
  and if $R=\O_C$ is the ring of integers of a spherically complete perfectoid field $C$, then we have
  \[ \TR_*(S;\Z_p) = \dRW^*_{S/\O_C}[\beta], \qquad |\beta| = 2 \]
   by \cite[Lemma 7.22]{MathewK1}. This identification can be made canonical by interpreting $\beta$ as tensoring with the ideal $(q-1)$ of $\Ainf(\O_C)$.
\end{example}

By Kan extending the Postnikov filtration (or its double-speed acceleration), one gets filtrations on these spectra for arbitrary animated rings, whose associated gradeds are given by derived versions of the same arithmetic invariants.

The genuine equivariant\footnote{$\RO(G)$-graded homotopy exists already for Borel $G$-spectra, but the notion is emphasized less there.} perspective provides a significant enlargement of the homotopy groups of a spectrum: a $G$-spectrum $E$ naturally has homotopy groups $\pi_V^H E$ for each subgroup $H\le G$ and each real virtual representation $V\in\RO(H)$ of $H$. This has a great deal of functoriality in $H$, leading to \emph{$\RO(G)$-graded homotopy Mackey functors $\mpi_\rog E$}. When $E$ is a $G$-$\E_\infty$-ring, we even get an \emph{$\RO(G)$-graded Tambara functor} in the sense of Angeltveit-Bohmann \cite{GradedTambara}.

\begin{remark}
  For all notation and background, the reader can consult
  \begin{itemize}
    \item \cite[Part 1]{SulSliceTHH} and \cite{HillHandbook} for equivariant homotopy theory in general;
    \item \cite[\sec4.2]{SulSlopes} for poly-representations;
    \item \cite[\sec2]{PrismsTambara1} for Tambara functors;
    \item \sec\ref{sec:notation} for all the notation used in this paper.
  \end{itemize}
  In particular, Notation \ref{not:canonical} elaborates what we mean by ``canonical''. In short, $\RO(G)$-graded homotopy groups naturally have the structure of \emph{spans} relative to the $\Z$-graded homotopy groups, and it is as such that we identify them.
\end{remark}

In this paper, we will generalize the above two computations above to $\RO(G)$-grading. Here are some reasons for doing so.

\begin{motivation}
  In general, the $\RO(G)$-graded homotopy groups of a $G$-spectrum tend to give a simpler or more holistic picture, in which the $\Z$-graded groups sit as a funny diagonal. For example, while the integer-graded $\TCmin$ mixes the Nygaard filtration with Breuil-Kisin twists, the $\RO(\T)$-graded groups will give us the Nygaard filtration on the nose:
  \[
    \TCmin_{2i}(R;\Z_p) = \Nyg^{\ge i}\hat\prism_R\{i\}
    \quad\text{but}\quad
    \TCmin_{i\lambda}(R;\Z_p) = \Nyg^{\ge i}\hat\prism_R.
  \]

  This example showcases another important point: the $\RO(G)$-graded groups considered here do not yield ``new'' arithmetic invariants. Rather, they vastly increase the supply of additional \emph{filtrations} and \emph{twists} that can be passed through the topological dictionary.
\end{motivation}

\begin{motivation}
  The original motivation for studying $\RO(G)$-graded $\TR$ specifically comes from algebraic $K$-theory: for any commutative ring $A$, we may consider the reduced $p$-adic $K$-theory of truncated polynomial algebras over $A$,
  \[ \~\K(A[x]/x^e;\Z_p) \defeq \fib\Big( \K(A[x];\Z_p) \to \K(A;\Z_p) \Big) \]
  By \cite[Proposition 8]{HesselholtHandbook}, this can be understood in terms of ``$\RO(G)$-graded big $\TR$'': there is a fiber sequence
  \begin{equation}
    \label{eq:syntomic-drw}
    \~\K(A[x]/x^e;\Z_p) \to \bigTR(\Sigma^{\hsl 1}\THH(A;\Z_p)) \to \bigTR(\Sigma^{\hsl{1/e}}\THH(A;\Z_p))
  \end{equation}
  and the latter two terms are extremely computable. The notation $\hsl\nu$ is recalled in Definition \ref{def:regular-poly}.

  A basic computation is that for $R$ perfectoid,
  \[ \bigTR_{2i-\hsl\nu}(R; \Z_p) = \W_{\lf i/\nu\rf}(R); \]
  this is essentially \cite[Proposition 9.1]{HMFinite} + \cite[Theorem 4.2.10]{HMCyclicPolytopes}, but see \cite[Corollary 4.6]{SulSlopes} for the present formulation (which is mildly more general). This immediately implies the syntomic cohomology
  \[ \~{\Z_p(i)}(R[x]/x^e) = \frac{\W_{ei}(R)}{\V^e\W_i(R)}[-1]. \]

  In \cite{SulSlopes}, this syntomic cohomology is computed by first computing $\Nyg^{\ge *}\prism_{k[x]/x^e}$ and $\prism_{k[x]/x^e}$, but the big Witt vectors appear somewhat miraculously. This method makes the appearance of big Witt vectors clearer, and also generalizes to perfectoid rings. The calculation of $\~\Z_p(i)(R[x]/x^e)$ for $R$ perfectoid was done independently by Riggenbach  using different methods \cite[Theorem 1.1]{RiggenbachTruncated}, and we will build on some of his tools in \sec\ref{sec:qrsp}.

  We therefore view this method of computing syntomic cohomology of monoid algebras as superior to the one used in \cite[\sec3]{SulSlopes}. The main disadvantage is that it is extremely difficult to establish the fiber sequence \eqref{eq:syntomic-drw} \cite{HMCyclicPolytopes}. However, this should be much easier to do in the modern prismatic approach, and to extend to more general monoid algebras. We will not pursue this here, but our work is intended to be a repository of computations to accompany such an endeavor.

  The $\RO(G)$-graded $\TR$ groups relevant to algebraic $K$-theory are the ``negative range'' (integer minus actual representation), which is the easiest range to compute. In this range, the answer has in some sense been known since the author was two and a half years old:
  \begin{itemize}
    \item the base case was proved in \cite[Proposition 9.1]{HMFinite}, and re-proven in \cite[Theorem 8.3]{ROS1TR} from the perspective of the homotopy orbits to TR spectral sequence (our main tool in this paper);
    \item the base case essentially implies the result for smooth algebras over it by \cite[Theorem 11]{HesselholtHandbook};
    \item the perfectoid base case was done in \cite{SulSlopes}---but in this range it behaves exactly like the characteristic $p$ case.
  \end{itemize}
  However, it took the author over thirty years to understand this.

  Even in the negative range, we get two new things. Our canonical identifications make it possible to compute the functoriality of $\~\K(A[x]/x^e;\Z_p)$ (or more general monoid algebras) in $e$. We also obtain an interesting filtration on these groups by computing with cell structures instead of the homotopy orbits to $\TR$ spectral sequence. These applications will be deferred to a companion paper.
\end{motivation}

While the $\RO(G)$-graded groups themselves are not ``new'', they provide a home for something which is.

\begin{motivation}[Tambara structure]
  Given $K\le H\le G$, there are $\RO(G)$-graded Frobenius (restriction) and Verschiebung (transfer) maps that go
    \[\xymatrix{
        \pi^H_V E \ar@<0.5ex>[r]^-F & \pi^K_{\res^H_K(V)} \ar@<0.5ex>[l]^-V E
    }\]
    Since the restriction of a trivial representation is trivial, these restrict to $F$ and $V$ maps on the $\Z$-graded homotopy groups.

    Genuine equivariant homotopy theory emphasizes another map, the \emph{Norm\iffalse\footnote{Frobenius is capitalized after German mathematician Ferdinand Frobenius, Verschiebung is capitalized because it's a German noun, and people capitalize Restriction out of consistency. The Norm map is named for Canadian mathematician Norm Macdonald.}\fi{} map}, which goes
    \[ N^H_K\colon \pi^K_V E \morph \pi^H_{\ind^H_K(V)} E \]
    when $E$ is a $G$-$\E_\infty$-ring spectrum. Since the induction of a trivial representation is never trivial, except for the zero representation, studying the norm in nonzero $\Z$-graded degrees forces us to contemplate $\RO(G)$-graded degrees. Understanding these norm maps was our primary motivation for these calculations.

  As for the zero representation, the $\pi_0$ Norm was introduced in \cite{AngeltveitNorm} for Witt vectors (and therefore perfect prisms), and for transversal prisms in \cite{PrismsTambara1}. Two applications of the latter are a conceptual construction of de Rham-Witt comparison maps for prismatic cohomology \cite{Molokov}, \cite[\sec4.1]{PrismsTambara1}, and (at odd primes) a refinement of the prismatic logarithm \cite{MaoLogarithm}. It would be very interesting to try to generalize these applications to $\RO(G)$-graded norms.
\end{motivation}

\begin{warning}
  Unfortunately, much of the literature refers to the sum-over-conjugates map
  \[ X_{hG} \to X^{hG} \]
  as the norm, where $X_{hG}$ denotes the homotopy orbits in spectra. By analogy with Galois theory, this should really be called the \emph{trace}. The true norm is the analogue of the above where the orbits are taken in $\E_\infty$-monoids and we take a product over conjugates rather than a sum.
\end{warning}

\begin{motivation}
  $\RO(G)$-graded homotopy groups also serve as calculational input for studying the equivariant slice filtration. Understanding the role of norms is also very important for understanding the slice filtration, since it is tailor-made to interact well with them. The author's original entrypoint into this subject was to attempt to modify the BMS filtration by replacing the Postnikov filtration with the slice filtration, which amounts to studying the slice filtration on $\TCmin(S;\Z_p)$ and $\TP(S;\Z_p)$ when $S$ is quasiregular semiperfectoid. This was our primary motivation for the calculations in \sec\ref{sec:qrsp}, and is now within reach.
\end{motivation}

\subsection{Borel spectra and prismatic cohomology}
We first compute the $\RO(\T)$-graded $\TCmin$ and $\TP$ of quasiregular semiperfectoid rings (\sec\ref{sec:qrsp}). Let us quickly review the representation theory of $\T$ and set some notation. The complex representation ring of $\T$ is $\RU(\T)=\Z[\lambda^\pm]$, with $\lambda$ the natural representation of $\T$ on $\C$; explicitly, $\lambda^i=\C$ with $z\in\T$ acting as multiplication by $z^i$. Complex conjugation induces a real isomorphism $\lambda^i \isom \lambda^{-i}$, so the real representation ring is $\RO(\T)=\Z[\lambda]$. Since we are only concerned with the homotopy type of representation spheres, we will work with $\Z[\lambda]$; but $\lambda^0$ will still mean $\C^{\mathrm{triv}}$ rather than $\R^{\mathrm{triv}}$. In the $p$-typical context, we have $S^{\lambda^i} \simeq S^{\lambda^j}$ iff $v_p(i)=v_p(j)$, so we set $\lambda_i=\lambda^{p^i}$ and $\lambda_\infty=\lambda^0$ (so we should really say $\RO(C_{p^\infty})$-graded). 
In addition to the irreducible decomposition
\[ \alpha = k_0 \lambda_0 + k_1 \lambda_1 + \dotsb + k_n \lambda_n + k_\infty \lambda_\infty, \]
it will be convenient to encode representations via the eventually-constant sequence
\[ d_n(\alpha) = \dim_\C \alpha^{C_{p^n}} \]
which we write as $(d_0,\dotsc,d_n;d_\infty)$. This encoding was introduced by Gerhardt and Angeltveit-Gerhardt \cite{TeenaTR,ROS1TR}, and as we will see, it is suited to working in ``ghost coordinates''.

Recall the ideals
\[ I_n \defeq I\phi(I)\dotsm\phi^{n-1}(I) \]
from \cite[Notation 2.2.2]{APC} (warning: our notation is off-by-one from that used in \cite{PrismsTambara1}).

\begin{theorem}
  \label{thm:qrsp}
  Let $A$ be a $\delta$-ring admitting a spherical lift $\S_A$, and let $S$ be a relatively quasiregular semiperfectoid $A$-algebra with associated prism $(\prism_{S/A},I)$. Let
  \begin{align*}
    \alpha
    &=k_0\lambda_0 + \dotsb + k_n\lambda_n + k_\infty \lambda_\infty\\
    &= (d_0,\dotsc,d_n;d_\infty)
  \end{align*}
  be a virtual $C_{p^\infty}$-representation. Then there are canonical identifications
  \begin{align*}
    \TCmin_\alpha(S\rel\S_A)
      &= \Nyg^{\ge d_0}\hat\prism^{(1)}_{S/A}\{k_\infty\} \tnsr \bigotimes_{i=1}^n I_i^{k_i}
    & \TP_\alpha(S\rel\S_A)
    &= \hat\prism^{(1)}_{S/A}\{k_\infty\} \tnsr \bigotimes_{i=1}^n I_i^{k_i}\\
    &= \hat\prism^{(1)}_{S/A} \tnsr \bigotimes_{i=0}^\infty \phi^{i-1}(I^{d_i})
    &&= \hat\prism^{(1)}_{S/A} \tnsr \bigotimes_{i=1}^\infty \phi^{i-1}(I)^{d_i}
  \end{align*}
  To interpret the second line: we have absorbed the Breuil-Kisin twist into the product via $\prism\{1\} = \bigotimes_{i=1}^\infty \phi^{i-1}(I)$, and note that $\phi^{i-1}(I^{d_i}) = \phi^{i-1}(I)^{d_i}$ for $i>0$.
\end{theorem}

\begin{warning}
  Given a $G$-spectrum $E$, there are two notations for the integer-graded homotopy groups of $E^G$: either $\pi_n E^G$ or $E^G_n$. We can imagine these as
  \begin{align*}
    \pi_n E^G &\defeq [S^n, E^G]_\Sp\\
    E^G_n &\defeq [S^n, E]_{\Sp^G}
  \end{align*}
  In other words, for the first version we take $G$-fixed points of $E$ and then take $\pi_n$ of the resulting ordinary spectrum, while in the second we compute the mapping space in $\Sp^G$ directly.

  In the $\RO(G)$-graded context, the first version \emph{does not make sense}: for $\alpha\in\RO(G)$, we define
  \[ E^G_\alpha \defeq [S^\alpha, E]_{\Sp^G} \]
  (where $S^\alpha$ is the one-point compactification of $\alpha$ when $\alpha$ is an actual representation), and this cannot be recovered from the plain spectrum $E^G$. We will often abuse notation by writing $\pi_\alpha E^G$, but it should be remembered that this is technically incorrect.

  Thus, when we write $\TCmin_\alpha$ and $\TP_\alpha$, these are defined as
  \begin{align*}
    \TCmin_\alpha &\defeq [S^\alpha, \THH^{E\T_+}]_{\Sp^\xi} &
    \TP_\alpha &\defeq [S^\alpha, \widetilde{E\T}\tnsr \THH^{E\T_+}]_{\Sp^\xi}\\
    &= \pi_0 (\Sigma^{-\alpha} \THH)^{h\T} &
    &= \pi_0 (\Sigma^{-\alpha} \THH)^{t\T}
  \end{align*}
  with $\Sp^\xi$ the category of cyclonic\footnote{We use the notation $\Sp^\xi$ because $\xi$ looks like a tornado.} spectra (i.e.\ proper-genuine $\T$-spectra).
\end{warning}

%

The proof of Theorem \ref{thm:qrsp} is very easy, and the argument works much more generally; all we need to do is formulate it. Recall the following generalization of cyclotomic polynomials.
\begin{definition}[{\cite[Remark 2.4.3]{GNS}}]
  \label{def:fgl-cyclotomic}
  Let $\G$ be a formal group law with $n$-series $[n]_\G(t)$. Define
  \[ \Phi_{n,\G}(t) = \prod_{d\mid n}[d]_\G(t)^{\mu(n/d)} \]
  where $\mu$ is the M\"obius function. Note that by M\"obius inversion we have
  \begin{align*}
    [n]_\G(t) &= \prod_{d\mid n}\Phi_{d,\G}(t)
  \end{align*}
\end{definition}

\begin{theorem}[Remark \ref{rmk:even}]
  \label{thm:even}
  Let $E$ be a $\T$-equivariant even $\E_\infty$-ring with associated formal group $\G$. Let $A\{i\} = \pi_{2i} E^{t\T}$, and let $\Fil^\bullet A\{i\}$ be the abutment filtration of the Tate spectral sequence. Fix a choice of complex orientation $t\in A\{-1\}$, define $\<n\>_\G=\frac{[n]_\G(t)}t\in A$ and $\Phi_{n,\G}=\Phi_{n,\G}(t)\in A\ (n>1)$ as above, and assume that $\<n\>_\G$ is a non-zerodivisor in $A$ for all $n$.

  Let $\alpha$ be a virtual complex $\T$-representation with irreducible decomposition
  \[
    \alpha = \sum_{i=0}^\infty k_i \lambda^i
  \]
  and dimension sequence $d_n = \dim_\C \alpha^{C_n},\, d_\infty = \dim_\C \alpha^\T$.

  Then the $\RO(\T)$-graded homotopy groups of $E^{h\T}$ and $E^{t\T}$ are given canonically by
  \begin{align*}
    E^{h\T}_\alpha
      &= \Fil^{d_1}A\{k_0\}\tnsr_A \bigotimes_{n>1} (\<n\>_\G A)^{\tnsr k_n}\\
      &= \Fil^{d_1} A \tnsr_A \bigotimes_{n>1} (\Phi_{n,\G} A)^{\tnsr d_n}\\
      &\fade{\defeq \Fil^{d_1} A\{d_\infty\} \tnsr_A \bigotimes_{n>1} (\Phi_{n,\G} A)^{\tnsr (d_n-d_\infty)}}\\
    E^{t\T}_\alpha
      &= A\{k_0\}\tnsr_A \bigotimes_{n>1} (\<n\>_\G A)^{\tnsr k_n}\\
      &= \bigotimes_{n>1} (\Phi_{n,\G} A)^{\tnsr d_n}\\
      &\fade{\defeq A\{d_\infty\} \tnsr_A \bigotimes_{n>1} (\Phi_{n,\G} A)^{\tnsr (d_n-d_\infty)}}
  \end{align*}
\end{theorem}

We can get more interesting invariants by going ``more genuine'' and considering the ``higher topological negative cyclic homology'' spectra $\TC^{(-n)} \defeq (\TR^n)^{h\T/C_{p^n}}$. We study this in \sec\ref{sub:iterated-nygaard}.

\subsection{Polygonic spectra and de Rham-Witt complexes}
We next shift our attention to $\TR$. Since the Restriction maps involve geometric fixed points, $\TR$ is not graded by representations of $C_{p^\infty}$ or $\T$, but rather \emph{poly-representations} of $\Z_p$ or $\widehat\Z$ (introduced in \cite[\sec4.2]{SulSlopes}). A poly-representation $\m\alpha$ consists of representations $\m\alpha(\Z/n)\in\RO(\Z/n)$ for each $n$, compatible under the maps
\begin{align*}
  \RO(\Z/mn) &\to \RO\left(\frac{\Z/mn\Z}{m\Z/mn\Z}\right)\\
              &= \RO(\Z/m)\\
  \alpha &\mapsto \alpha^{m\Z/mn\Z}
\end{align*}
The name \emph{poly-}representation is chosen partly because it is a family of many representations, but more importantly to emphasize the connection with polygonic spectra \cite{KMNPolygonic}: the representation spheres $S^{\m\alpha(\Z/n)}$ assemble into a polygonic spectrum $S^{\m\alpha}$ (which in particular is $\tnsr$-invertible). We will mostly work with $p$-typical poly-representations, in which case $\widehat\Z$ is replaced by $\Z_p$. We will study the polyrepresentation-graded $\TR$ of smooth algebras $S$ over a perfectoid base $R$.

\begin{remark}
Say that a continuous representation $V$ of $\widehat\Z$ is \emph{quasifinite-dimensional}\footnote{Of course it doesn't make sense for a dimension to be quasifinite; this term is chosen for consistency with the \emph{quasifinitely genuine} $\Z$-spectra from \cite{KMNPolygonic}.} if
\[ V(\Z/n) \defeq V^{n\widehat\Z} \fade{{}\curvearrowleft \Z/n} \]
is finite-dimensional for all $n>0$ (but $V^e$ could be countably-infinite-dimensional). Then a poly-representation is equivalently an element of the Grothendieck group $\RO(\widehat\Z)^\land$ of quasifinite-dimensional $\widehat\Z$-representations, a completion of the usual representation ring $\RO(\widehat\Z)$.
\end{remark}

Flipping and reversing the dimensional encoding, we write a $p$-typical poly-representation as
\[ \m \alpha = (\m d_0, \m d_1, \dotsc), \]
where $\m d_i=\dim_\C\m\alpha(C_{p^n})$ should be thought of $\dim_\C(\text{``$\alpha^{1/C_{p^i}}$''})$, and indeed we sometimes write $d_{-i}(\m\alpha)$ instead of $\m d_i$.


The pioneering work on $\RO(G)$-graded $\TR$ is by Gerhardt and Angeltveit-Gerhardt \cite{TeenaTR,ROS1TR}, who gave an algorithm to compute $\TR^n_\rog(k)$ as $W_n(k)$-modules when $k$ is a perfect $\F_p$-algebra. We revisited their work in \cite{SulROGTHH} to compute the full rings $\TF_\rog(k)$ and $\TF_\rog(R^+;\Z_p)$ when $k$ is a perfect $\F_p$-algebra and $R^+$ is a $p$-torsionfree perfectoid ring. This in principle implies the computation of $\TR_\rog$ via the long exact sequences
\[ \TF_{\alpha+\lambda_i}(R;\Z_p) \xra{a_{\lambda_i}} \TF_\alpha(R;\Z_p) \to \dbox{$\TR^i_\alpha(R;\Z_p)$} \xra\partial \TF_{\alpha+\lambda_i-1}(R;\Z_p), \]
but it is nontrivial to unpack this; we will instead proceed directly by adapting the homotopy orbits to TR spectral sequence. The results of \cite{SulROGTHH} may be informally summarized as follows:
\begin{itemize}
  \item in even degrees, $\TF_\rog(R;\Z_p)$ is torsionfree and satisfies an $\RU(\T)$-graded form of B\"okstedt periodicity;
  \item in the characteristic $p$ case, there are additional classes witnessing the ``collision'' of the distinguished elements $d=\phi(d)=\phi^2(d)=\dotsb$;
  \item the odd degrees contain deep caverns with dangerous monsters\footnote{Infamously, \cite{SulROGTHH} contains an error related to extension problems in these odd degrees. The author now knows how to fix this and will update the paper soon. This error does not affect any of the computations in this paper on account of restricting to the contiguous range.}.
\end{itemize}

In light of the extraordinary complexity of even the base case, we will not attempt to compute the full $\TR_\rog(S)$, but rather focus our attention on certain safe regions:
\begin{definition}
  Consider the following swathes of poly-representations.

  \begin{itemize}
    \item the \emph{negative range} consists of $\rog=*-V$ with $*\in\Z$ and $V$ an actual poly-representation; equivalently, $\m\alpha=(\m d_0,\m d_1,\dotsc)$ is in the negative range if $\m d_i \ge \m d_{(i+1)}$ for all $i$. This is the range that arises in algebraic $K$-theory calculations. 
    \item the \emph{positive range} consists of $\rog=V-*$ with $*\in\Z$ and $V$ an actual poly-representation; equivalently, $\m\alpha=(\m d_0,\m d_1,\dotsc)$ is in the positive range if $\m d_i \le \m d_{(i+1)}$ for all $i$. This range is of interest because the Norm from the $\Z$-graded de Rham-Witt complex lands in here, allowing us to extend the Witt vector Norm \cite{AngeltveitNorm,PrismsTambara1} to the de Rham-Witt complex. (Also, a tiny portion of it is used in computing the slice filtration \cite[Lemma 4.11]{SulSliceTHH}.)
  \item the \emph{contiguous range} consisting of those $\m\alpha=(\m d_0,\m d_1,\dotsc)$ for which $\{i \mid d_i(\m\alpha)\ge 0\}$ is a contiguous interval $[-\ell, -k]$ or $(-\infty, k]$. We say that $\m\alpha$ is \emph{non-negatively concentrated in $[-\ell, -k]$} or $(-\infty, k]$.
  \end{itemize}
\end{definition}

The negative and positive ranges are exactly the ranges that can be handled by the cellular method; however, this is a very difficult way to proceed. We instead skip straight to the contiguous range and employ the \emph{homotopy orbits to TR spectral sequence} (HOTRSS) introduced in \cite{ROS1TR} and upgraded in \cite{SulROGTHH}. The homotopy orbits to TR spectral sequence can handle all of $\TR_\rog$, but the answer becomes unwieldy outside of the contiguous range.

We need two more ingredients to state the result. Recall that the \emph{canonical filtration} of the de Rham-Witt complex is given by
\begin{align*}
  \Fil^{k+1} \dRW[n]^j_{S/R}
    &= \ker(\dRW[n]^j_{S/R} \xra{R^{n-k}} \dRW[k]^j_{S/R})\\
    &= \im(V^{k+1} + dV^{k+1})
\end{align*}
for $n\in\N\cup\{\infty\}$. (Remember that in our notation Witt vectors start at $\rmW_0$.) Next, in the $\Z$-graded case, Hesselholt's theorem \cite{LarsPTyp} says that $\TR_*$ splits into contributions from individual $\dRW[\bullet]^j$ terms. This is also true in the fancy-graded case (either from the proof of Hesselholt's theorem, or by the homotopy orbits to TR spectral sequence), so we let $\TR_{\rog,j}(S)$ denote the $\dRW^j$ contribution to $\TR_{\rog+j}(S)$. Pulling out $j$ this way, the computation becomes essentially the same for the base case $S=R$, although we have to be slightly more careful with $\Fil$.

\begin{theorem}[Contiguous range of $\TR_\rog$, Theorem {\ref{thm:tr-ctg}}]
  \label{thm:intro-tr-ctg}
  Let $\m\alpha = (\m d_0,\dotsc, \m d_n)$ be a $\Z/p^n$-polyrepresentation which is non-negatively concentrated in $[-\ell, -k]$.

  If $S$ is $p$-completely smooth over a $p$-torsionfree perfectoid ring $R^+$, then there are canonical identifications
  \begin{align*}
    \TR^n_{\m\alpha,j}(S\rel\SAinf{R^+})
      &=
        \Fil^k \dRW[\ell]^j_{S/R^+}
        \tnsr
        \bigotimes_{i=0}^\ell \phi^{-i}(\xi)^{\m d_i}A\\
    \TR^n_{\m\alpha-1,j}(S\rel\SAinf{R^+})
      &=
        \frac{\dRW[k-1]^j_{S/R^+}}{\phi^{-k}(\xi)^{\m d_k}\dotsm\phi^{-\ell}(\xi)^{\m d_\ell}}
        \tnsr
        \bigotimes_{i=0}^{k-1} \phi^{-i}(\xi)^{\m d_i}A
  \end{align*}
  with $A=\Ainf(R^+$). This also holds for $n=\infty$ if $\ell<\infty$ or if $R^+=\O_C$.

  If $S$ is smooth over a perfect $\F_p$-algebra $k$, then there are canonical identifications
  \begin{align*}
    \TR^n_{\m\alpha,j}(S\rel\SW k)
      &= \Fil^{k-(\m d_k+\dotsb+\m d_\ell)} \dRW[\ell]^j_S\tnsr p^{\m d_0 + \dotsb + \m d_\ell}W\\
    \TR^n_{\m\alpha-1,j}(S\rel\SW k)
      &= \frac{\dRW[k-1]^j_S}{p^{\m d_k+\dotsb + \m d_\ell}}\otimes p^{\m d_0 + \dotsb + \m d_{(k-1)}}W
  \end{align*}
  with $W=\rmW(k)$. This also holds for $n=\infty$.
\end{theorem}

While the restriction to the contiguous range is disappointing, note that the contiguous range is already massively larger than the range of degrees that have been previously studied; besides, the point of this paper is to open the doors to more $\RO(G)$-graded ideas in the field, not to serve as the final word on it. Restricting to the contiguous range is most distressing for the following reason:
\[ \{\m\alpha \mid \m\alpha + i \text{ is contiguous for all } i \in \Z\} \]
is easily seen to be the union of the positive and negative ranges, so only in these ranges can we fully mix our result with the Postnikov filtration.

\begin{remark}
  The notations $\Z/n$ and $C_n$ are often used interchangeably for a cyclic group of order $n$. In this paper, we will think of $\Z/n$ as a quotient of $\widehat\Z$, while we think of $C_n$ as embedded in $\T$ (so perhaps a better notation would be $\Z/n\{-1\}$ as in Example \ref{ex:tor-ext}, the Pontryagin dual $(\Z/n)^\land$, or $\mu_n$, although the latter can be confused for a group scheme). As we will see in \sec\ref{sec:dRW}, the passage between these two perspectives involves the Cartier isomorphism (Example \ref{ex:cartier}); it would be interesting to know if there is a deeper connection between Pontryagin duality and the Cartier isomorphism.
\end{remark}

\begin{remark}[Poly-de Rham-Witt complex]
  Our computations imply a new expansion of the de Rham-Witt complex. In general, there is a map
  \[ \dRW[\bullet]^*_{S/R} \to \TR^\bullet_*(S \rel \SAinf R) \]
  because $\dRW[\bullet]^*$ is defined exactly as the universal example of ``the kind of structure $\TR^\bullet_*$ has''. This is quite vague, but in fact different authors get away with axiomatizing different fragments of the $d,F,V,R$ structure; and when $\TR$ becomes more complicated than the classical algebraic case, we can always generalize the definition of this structure, as with the $p$-typical Cartier complexes of \cite[Definition 3.34]{TCart} or the $\eta$-deformed ($p$-typical) Cartier complexes of \cite[Definition 4.17]{CycSyn}.

  The polyrepresentation-graded groups $\TR^\bullet_\rog$ possess $d$, $F$, $V$, and $R$ maps just as the $\Z$-graded groups do, but now also have $N$ maps. We can define the four variants ($p$-typical vs big, truncated vs full) of the \emph{poly-de Rham Witt complex}
  \begin{gather*}
    \dRW^\rog_{S/R}, \qquad \dRW[\bullet]^\rog_{S/R},\\
    \DRW^\rog_{S/R}, \qquad \DRW[\bullet]^\rog_{S/R}
  \end{gather*}
  as the universal examples of this structure. We will not study the poly-de Rham-Witt complex in this paper, but our computations provide an upper bound for it. Again, the poly-de Rham-Witt complex is a filtered and twisted version of the usual de Rham-Witt complex, rather than a totally new object, with the norms being the most interesting new piece of structure.
\end{remark}

\subsection{Norms}
While the $\RO(G)$-graded groups are interesting in their own right, our main reason for computing them is to be able to understand the Norm (even for nonzero integer degrees). Since the Norm is not very well-known, we begin by recalling the Norm map of Witt vectors from \cite{AngeltveitNorm,PrismsTambara1}.

\subsubsection{Norm map of Witt vectors}
This is a multiplicative function $N\colon \rmW(R) \to \rmW(R)$ characterized in ghost coordinates by
\begin{align*}
  FN(f) &= f^p\\
  N &\equiv \id \bmod V
\end{align*}
It is conceptually clarifying to also consider iterates of the Norm, characterized in ghost coordinates by
\begin{align*}
  F^k N^k(f) &= f^{p^k}\\
  F^{k-r} N^k(f) &\equiv f^{p^{k-r}} \bmod V, \quad 1\le r\le k\hspace{-5.1em}
\end{align*}
This descends to a map $N\colon \rmW_n(R) \to \rmW_{n+1}(R)$, in which case $N^n\colon R = \rmW_0(R) \to \rmW_n(R)$ is the familiar multiplicative lift, and in Witt coordinates
\begin{equation}
  \label{eq:witt-coords}
  (a_0, \dotsc, a_n) = N^n a_0 + VN^{n-1} a_1 + \dotsb + V^n a_n.
\end{equation}
This essentially says that $\rmW_\bullet(R)$ is the free $\Z_p$-Tambara functor on $R$, where $R$ is viewed as a ring with (trivial) $\Z_p$-action \cite{BrunWitt}.

Recall from \cite[Proposition 1]{JoyalDeltaWitt} and \cite[Lemma B.1]{DrinfeldFormalGroup} that a $\delta$-ring admits operations $\delta_\bullet$ and $\theta_\bullet$ satisfying
\begin{align*}
  \phi^n(f) &= \sum_{k=0}^n p^k \delta_k(f)^{p^{n-k}} &
  \phi(f^{p^{n-1}}) &= f^{p^n} + p^n \theta_n(f)
\end{align*}
(In particular, $\delta_0 = \id$ and $\delta_1 = \theta_1 = \delta$.) These give implicit and explicit formulas for iterates of the Norm:
\begin{align*}
  \id &= N + V\delta_1 &
  N &= \id - V\theta_1\\
  \id &= \sum_{k=0}^n V^k N^{n-k} \delta_k &
  N^k &= \id - \sum_{k=1}^n V^k \theta_k
\end{align*}
In particular, we can view $\id = \sum V^k N^{n-k} \delta_k$ as a ``Witt coordinates'' decomposition of the identity, generalizing \eqref{eq:witt-coords} to infinite level.

\subsubsection{de Rham-Witt norms}
Now let $S$ be $p$-completely smooth over a perfectoid ring $R$, let $A=\Ainf(R)$, and assume for the moment that $p>2$. By Theorem \ref{thm:intro-tr-ctg}, there are Norm maps that go
\begin{align}
  \label{eq:norm-dRW-map}
  N\colon \dRW[n]^j_{S/R} &\to \dRW[n+1]^{j}_{S/R} \fade{{}\tnsr_A ({}\dotsb)} + \Fil^1 \dRW[n+1]^{pj}_{S/R} \fade{{}\tnsr_A ({}\dotsb)}\\
  N^k\colon \dRW[n]^j_{S/R} &\to \sum_{r=0}^k \Fil^r \dRW[n+k]^{p^rj}_{S/R} \fade{{}\tnsr_A ({}\dotsb)}
\end{align}
where the sums are direct sums for $j>0$, and we are suppressing the ideal twists for simplicity; the precise twists are stated in Theorem \ref{thm:norms-dRW}. These are characterized in ghost coordinates by
\begin{align*}
  FN(\omega) &= \omega^p & F^k N^k(\omega) &= \omega^{p^k}\\
  N &\equiv \psi^p \bmod \Fil^1 & F^{k-r} N^k(\omega) &\equiv \psi^{p^r}(\omega)^{p^{k-r}} \bmod \Fil^1
\end{align*}

Here $\psi^\ell$ is the Adams operation $\psi^\ell\colon \Omega^*_{S/R} \to \Omega^*_{S/R}$ given by $\psi^\ell(\omega\in\Omega^j_{S/R}) = \ell^j \omega$. These are well-known for $p\nmid \ell$ (see \cite[Proposition 9.14]{BMS2} and \cite[Proposition 6.4.12]{RaksitDDR}), but we need the $p$-primary versions, which we study in \sec\ref{sub:adams}.

By reduction to the $p$-torsionfree case, it suffices to find a natural formula which satisfies the above congruences. We would like to try
\[ N(\omega) = \psi^p(\omega) - V\delta(\omega), \]
but it is not clear how to define $\delta$ on differential forms. Also, this does not have the right shape for the codomain of \eqref{eq:norm-dRW-map}.

The missing formula is
\begin{equation}
  \label{eq:norm-formula-2}
  N(\omega) =  (p-V(1))\frac{\psi^p(\omega)}p + V\frac{\omega^p}p
\end{equation}
It is easy to check that, as a formal expression, this satisfies the necessary congruences. This formula does not make sense in the degree 0 case, but it does in positive degree. Indeed, if $\omega$ is a differential form of positive degree, then $\omega^p$ and $\psi^p(\omega)$ are divisible by $p$, in fact canonically so in the $p$-torsionfree case, which we can reduce to by functoriality. Beware that the $V$ in $V\frac{\omega^p}p$ is an $\RO(G)$-graded Verschiebung.

With hindsight, we can manually extend $\psi^k$ to the de Rham-Witt complex via $\psi^k(\omega\in\dRW^j_{S/R}) = k^j\omega$, and then set
\[ \delta(\omega) \defeq \frac{F\psi^p(\omega) - \omega^p}p \]
which makes sense for the same reasons as above---though it is not a homogeneous operation! We can then rearrange \eqref{eq:norm-formula-2} into the unified formula
\[ \psi^p(\omega) = N(\omega) + V\delta(\omega) \]
which works also for $\omega\in\dRW^0$.

We can similarly extend $\delta_\bullet$ and $\theta_\bullet$ to the de Rham-Witt complex, and view the Norms along with Verschiebungs as providing a ``Witt coordinates'' decomposition of the Adams operations $\psi^{p^k}$.

\begin{warning}
Although this extension is conceptually appealing, the topological perspective only justifies the definition of $\psi^k$ up to $\Fil^1$ and only justifies $\delta$ up to $\ker V$.
\end{warning}

\begin{warning}
  It is not obvious that \eqref{eq:norm-formula-2} simplifies to $N(\omega)=V\frac{\omega^p}p$ in characteristic $p$, since some of the factor $(p-V(1))$ gets absorbed into the ideal twist. However, a small calculation shows that this factor still vanishes in characteristic $p$, at least for $p>2$ (Corollary \ref{cor:norms-Fp}).
\end{warning}

\begin{remark}
  \label{rmk:norms-2}
  When $p=2$, everything above goes through for de Rham-Witt forms of even degree. Producing Norms for odd-degree differential forms at $p=2$ requires contending with sign representations, which we defer to future joint work with Jason Yeung.
\end{remark}

\begin{remark}
  Let us rename the Norm discussed above (inherited from the Norm on $\TR$) to the \emph{topological Norm}, and denote it by $N_\top$. We define the \emph{algebraic Norm} $N_\alg$ by undoing the twists (see the explicit description of the twists below):
  \[\xymatrix{
      \dRW[n]^j_{S/R} \ar[r]^-{N_\top} \ar@{..>}[dr]_-{N_\alg} & \dRW[n+1]^j_{S/R} \tnsr ({}\cdots) + \Fil^1 \dRW[n+1]^{pj}_{S/R} \tnsr ({}\cdots) \ar[d]^-{\text{\tiny$\begin{bmatrix}1&0\\0&\xi^{j(p-1)/2}\end{bmatrix}$}}\\
                                       & \dRW[n+1]^j_{S/R} + \Fil^1 \dRW[n+1]^{pj}_{S/R} 
  }\]

  This is somewhat easier to work with, and we do not understand the significance of being able to refine $N_\alg$ to $N_\top$. For example, in characteristic $p$ it is obvious that $N_\alg(\omega) = V\frac{\omega^p}p$ for $j>0$; while this is not obvious for $N_\top(\omega)$, a nontrivial calculation shows that it is still true, so $N_\top$ is not seeing anything more. It is conceivable that $N_\top$ and $N_\alg$ will differ more significantly for odd-degree forms at $p=2$.
\end{remark}

\begin{example}
  \label{ex:norm-dx}
  Let $p$ be odd, and let $S$ be $p$-completely smooth over a perfectoid ring $R$.
  If $|\omega|=1$, then
  \begin{align*}
    N^k(\omega) &= (p^k - Vp^{k-1})\omega\\
    \delta_k(\omega) &= F^k\omega \fade{{} = \delta^k(\omega)}\\
    \theta_k(\omega) &= 0 \quad \text{ for } k>1
  \end{align*}
  Note that this implies the values of $N,\delta_k,\theta_k$ on all of $\dRW^*$ by the addition/multiplication formulas.

  Let us be more precise about the twist. Setting $h=\frac{p-1}2$, the Norm goes
  \[ N_{\fade\top} \colon \dRW[n]^1_{S/R} \to \dRW[n+1]^1_{S/R} \tnsr \phi^{-1}\left(\xi_{n+1}^h(h)_\prism!\right) \]
  The symbol $(i)_\prism!$ is defined in Definition \ref{dfn:prismatic-factorial}, and in the $q$-de Rham context the twist becomes
  \[ \phi^{-1}\left((p^{n+1})^h_{q^{1/p^{n+1}}} (h)_q!\right) = \phi^{-1}\prod_{i=1}^h (p^{n+1}j)_{q^{1/p^{n+1}}} \]
  In particular, the factor $(h)_{q^{1/p}}!$ means that the $F$ map out of the codomain is $h!$ times the usual map $F\colon \dRW[n+1]^1_{S/R} \to \dRW[n]^1_{S/R}$. But $FN\omega = \omega^p = 0$ anyway, so we do not understand the importance of this twist here.

  Ignoring the factorial unit, the refinement of $N_\alg$ to $N_\top$ amounts to saying that $p-V(1)$ is divisible by $\phi^{-1}\left(\xi_{n+1}^{(p-1)/2}\right)$ (Remark \ref{rmk:snell} and \cite[Lemma 2.1]{MaoLogarithm}). Taking $n=0$ for simplicity, this means that
  \[ N^1_{0,\top}(\omega) \approx \phi^{-1}(\xi)^{h+1}\omega \tnsr \phi^{-1}(\xi)^h \in \dRW[1]^1_{S/R} \tnsr \phi^{-1}\left(\xi^h (h)_\prism!\right) \]
  where $\approx$ means up to units.
  In characteristic $p>2$, we have $\frac{p+1}2 \ge \frac{3+1}2 = 2$, so the term $\phi^{-1}(\xi)^{h+1} = p^{\frac{p+1}2}$ vanishes in $\rmW_1(k)$, showing that indeed $N^1_{0,\top}(\omega)=0$.
\end{example}

\begin{example}
  \label{ex:norm-even}
  Let $S=R[x_1, y_1, \dotsc, x_p, y_p]^\land_p$ over an arbitrary perfectoid ring $R$. Let
  \[
    \tau = \sum_{i=1}^p d x_i\, d y_i \in \dRW^2_{S/R} \qquad
    \omega = \prod_{i=1}^p d x_i\, d y_i \in \dRW^{2p}_{S/R}
  \]
  We have $\tau^p = p!\omega$, so the algebraic Norm is
  \begin{align*}
    N_\alg(\tau)
      &= (p^2-V(p))\tau + V(p-1)!\omega\\
      &= V(p-1)!\omega \qquad\text{in characteristic $p$}
  \end{align*}
  Similarly,
  \begin{align*}
    \delta(\tau)
    &=
    \frac{F\psi^p(\tau) - \tau^p}p\\
    &= p\sum_{i=1}^p x_i^p y_i^p \dlog x_i \dlog y_i - (p-1)!\omega
  \end{align*}

  The topological Norm goes
  \[
    N_\top\colon
    \dRW[n]^2_{S/R}
    \to
    \bigoplus\left(
      \begin{array}{r@{}l}
      \dRW[n+1]^2_{S/R} &{}\tnsr \phi^{-1}\left(\xi_{n+1}^{p-1}(p-1)_\prism!\right)\\
      \Fil^1 \dRW[n+1]^{2p}_{S/R} &{}\tnsr \dfrac{\phi^{-1}(p-1)_\prism!}{\xi^{p-1}}
    \end{array}\right)
  \]
  so is given by
  \[
    N_\top(\tau) =
    \frac
      {p^2-V(p)}
      {\phi^{-1}\left(\xi_{n+1}^{p-1}(p-1)_\prism!\right)}
    \tau \fade{{}\tnsr \phi^{-1}\left(\xi_{n+1}^{p-1}(p-1)_\prism!\right)}
    +
    V\omega \fade{{}\tnsr\dfrac{\phi^{-1}(p-1)_\prism!}{\xi^{p-1}}} 
  \]
  (remembering that $F(1\tnsr\xi) = 1$, since $F$ is a \emph{divided} Frobenius.)
\end{example}

\subsubsection{Breuil-Kisin and prismatic norms}

Back to the Borel case, one would also like to compute the Tambara structure on $\THH(S\rel\S_A)^{\{h,t\}C_{p^\bullet}}_\rog$, extending the author's $\pi_0$ prismatic norms from \cite{PrismsTambara1}. Although we have managed to compute a few of these, they are far more complicated than expected, and we do not yet see an easy conceptual story like we have for the de Rham-Witt norms or $\pi_0$ prismatic norms. In the interest of this paper seeing the light of day, we postpone these to a future revision of this article. (If one restricts to perfectoid rings $R$, it is relatively easy to see that the norms on $\THH(R;\Z_p)^{\{h,t\}C_{p^\bullet}}_\rog$ take generators to generators, which suffices for many applications; the difficulty is in pinning down units.)

\subsection{Overview}

We review all our notation in \sec\ref{sec:notation}. We study the $\RO(G)$-graded homotopy of even Borel spectra in \sec\ref{sec:qrsp}: primarily the $\RO(C_{p^\infty})$-graded $\TCmin$ and $\TP$ of quasiregular semiperfectoid rings, and some parallel results for even $\E_\infty$-rings with $\T$-action. We study the de Rham-Witt setting in \sec\ref{sec:dRW}. We study $\RO(G)$-graded norms in \sec\ref{sec:norms}.

A web interactive for exploring the homotopy orbits to TR spectral sequence is available at

\begin{center}
  \url{https://ysulyma.github.io/papers/poly/}
\end{center}

\subsection{Acknowledgments}

We thank Benjamin Antieau, Andrew Blumberg, Elden Elmanto, Lars Hesselholt, Jacob Lurie, Zhouhang Mao, Shubhodip Mondal, and Noah Riggenbach for helpful conversations related to this work. We are especially grateful to Antieau for discussions about the Adams operation $\psi^p$; Mao for explaining Example \ref{ex:tor-ext} to us; and to Riggenbach for clarifying the difference between absolute and relative Ngyaard filtrations for us.

%% file: notation.tex

\section{Notation}
\label{sec:notation}

Here we collect all the notation used in the paper (as well as some notation which does not appear in this paper but is used in related works by the author). In fact, the notation system we introduce should be regarded as the ``fundamental technical lemma'' from which all else in the paper follows. (This was also largely the case in \cite{SulSlopes}.) We will also chronicle the evolution of the notation.

We use $=$ for canonical isomorphisms, and $\cong$ for non-canonical isomorphisms; the meaning of ``canonical'' is elaborated in Notation \ref{not:canonical}. We also write $x\approx y$ to mean $x=uy$ for a unit $u$.

The ring we use are
\begin{align*}
  R &= \text{general perfectoid ring}\\
  R^+ &= \text{$p$-torsionfree perfectoid ring}\\
  k &= \text{perfect $\F_p$-algebra}\\
  C &= \text{spherically complete perfectoid field}\\
  S &=
    \begin{cases}
      \text{quasiregular semiperfectoid ring} &\text{in \sec\ref{sec:qrsp}}\\
      \text{$p$-completely smooth $R$-algebra} & \text{in \sec\ref{sec:dRW}}
    \end{cases}\\
  A &=
  \begin{cases}
    \text{base $\delta$-ring} & \text{in \sec\ref{sec:qrsp}}\\
    \Ainf(R) & \text{in \sec\ref{sec:dRW}}
  \end{cases}\\
  W &= \rmW(k)
\end{align*}
We index $p$-typical Witt vectors and $p$-typical $\TR$ starting from $0$, so that $\TR^n = \THH^{C_{p^n}}$ and $\rmW_n(\F_p)=\Z/p^{n+1}$. We write $\bigTR$ for big $\TR$ and $\W$ for big Witt vectors. We write $\phi$ for the Frobenius on $A$, $\varphi$ for the Frobenius on prismatic cohomology, and $\~\varphi$ for a divided Frobenius.

If $(A, I)$ is a prism, then $I_n = I\phi(I)\dotsm\phi^{n-1}(I)$ \cite[Notation 2.2.2]{APC}. We recall the following elements from \cite[Lemma 3.12]{BMS1}: if $R$ is a perfectoid ring with corresponding perfect prism $(A, I)$, then
\begin{align*}
  \~\xi & \text{ is a generator of }I\\
  \xi &= \phi^{-1}(\~\xi) \text{ is a generator of }\Nyg^{\ge 1}\prism_R = \phi^{-1}(I)\\
  \~\xi_n &= \~\xi\phi(\~\xi)\dotsm\phi^{n-1}(\~\xi) \text{ is a generator of } I_n\\
  \xi_n &= \xi\phi^{-1}(\xi)\dotsm\phi^{-(n-1)}(\xi) \text{ is a generator of } \phi^{-n}(I_n)
\end{align*}
Beware that this notation is off-by-one from that used in \cite{PrismsTambara1}.

In the $q$-de Rham case, we have the element $q=[\epsilon]\in\Ainf(R)$ from from \cite[Proposition 3.17]{BMS1}, allowing us to make these more concrete as $q$-analogues:
\begin{align*}
  \mu &= 1-q\\
  \~\xi_n &= (p^n)_q = \frac{\phi^n(\mu)}\mu\\
  \xi_n &= (p^n)_{q^{1/p^n}} = \frac{\mu}{\phi^{-n}(\mu)}
\end{align*}

We will sometimes use the \emph{relative} theory
\[ \TR(S \mid \SW k) \quad\text{or}\quad \TR(S \mid \SAinf{R}) \]
where $\SAinf R := \SW{R^\flat}$, and $\SW k$ is the functor of spherical Witt vectors \cite[Example 5.2.7]{Ell2}, \cite{AntieauSWitt}, \cite[\sec2]{ChromaticNullstellensatz}, \cite{AlternativeSWitt}. This is equivalent to $\TR(S;\Z_p)$ by \cite[Lemma 5.20]{Zhouhang}, an avatar of the identification between absolute and relative prismatic cohomology under a perfect prism \cite[Lemma 4.8]{Prismatic}, \cite[Proposition 4.4.12]{APC}, \cite[Remark 11.9]{BMS2}. However, the relative notation is helpful for emphasizing the base in results that depend on it, and for reducing the context needed when reading our formulas on the go.

The circle group is denoted by $\T$. We write $\lambda^n$ for the 1-dimensional complex $\T$-representation in which $z\cdot w=z^n w$ for $z\in\T$ and $w\in\C$; this is the $n^{\text{th}}$ power of $\lambda=\lambda^1$ in the complex representation ring $\RU(\T)$. We set $\lambda_i = \lambda^{p^i}$, with $\lambda_\infty = \lambda^0 = \C^\triv$. Every real representation of $\T$, and every real orientable representation of $C_n$, is of the form $\alpha$ or $\alpha-1$ for a complex representation $\alpha$. 

\begin{notation}[Dimensional encoding]
  The first major maneuver is to index a representation $\alpha$ by its (eventually constant) dimension sequence $d_\bullet(\alpha)=\dim_\C \alpha^{C_{p^\bullet}}$, rather than the coefficients of its irreducible decomposition. This insight is due to Angeltveit-Gerhardt \cite{ROS1TR}, and would be a good lesson for the rest of equivariant homotopy theory to learn---for example, it changes the ``double cone'' pattern seen in computations into a simple quadrant pattern.

    To appreciate the difference, note that the dimensional encoding is suited for computing with the isotropy separation square, while the irreducible encoding is suited for computing with cell structures. The isotropy separation method is multiplicative, and this is a major advantage over the cellular method, which is not. See \cite{Zengv1} for a detailed comparison of these methods for computing $\mpi_\rog(\Z)$ over $C_{p^2}$.

    Thus, we write a virtual representation $\alpha$ as $\alpha=(d_0,\dotsc,d_n;d_\infty)$, where $d_i = \dim_\C(\alpha^{C_{p^i}})$, and $d_\infty$ is the tail of the eventually constant sequence. Note that $(d_0,\dotsc,d_n;d_\infty) = (d_0,\dotsc,d_n,d_\infty;d_\infty)$; these two encodings give slightly different spectral sequences, but the end result is of course the same. If $d_i(\alpha)=d_\infty$ for $i\ge n$, we sometimes refer to $\alpha$ as a $C_{p^n}$-representation, although this is slightly abusive.

  To convert between this encoding and the irreducible decomposition $\alpha=k_0\lambda_0 + \dotsb + k_n \lambda_n + k_\infty \lambda_\infty$, the formulas are
  \begin{align*}
    d_i &= \sum_{i\le j\le \infty} k_j\\
    k_i &= d_i - d_{i+1} \qquad i <\infty\\
    k_\infty &= d_\infty
  \end{align*}

  Beware that the above encoding does not distinguish between $\lambda^n$ and $\lambda^{p^{v_p(n)}}$. This is harmless for most purposes, but it does twist restriction ($F$) maps by the unit $n/p^{v_p(n)}$, and one needs to keep track of this when trying to pin down the Norm exactly (rather than just up to units).
\end{notation}

\begin{notation}[HOTRSS]
  The next powerful move is to unravel the induction computing $\TR^n_\rog$ into a spectral sequence, the homotopy orbits to TR spectral sequence. This again is due to Angeltveit-Gerhardt \cite[\sec3]{ROS1TR}. We gave a more conceptual construction of this spectral sequence (that makes it easy to see the multiplicative and Mackey structure) in \cite[\sec3]{SulROGTHH}, and it is this formulation that we will use. We will always write this as a two-column spectral sequence in Adams grading, where the left column computes $\pi_{\alpha-1}$ and the right column computes $\pi_\alpha$.
\end{notation}

\begin{notation}[Fool's gold]
  To make the multiplicative and Mackey structures transparent, we shift to indexing using the gold elements. We learned this from Dylan Wilson and \cite{Zeng}. Actually, one of the major advances of this paper is to \emph{remove} the gold elements---but first let us just recall what they are.

  Given an actual representation $V$, we write $a_V \in \pi^G_{-V} \S_G$ for the \emph{Euler class} of $V$, obtained by suspending the inclusion $\{0\}\mono V$ (see \cite[Definition 3.11]{HHR}). In particular we can map $a_V$ to any $G$-spectrum $X$. Note that this manifestly does not depend on any choices, i.e.\ it is ``canonical''. We will not actually get rid of the $a_V$, but they will be hidden in what follows.

  In favorable cases there are also \emph{Thom classes} $u_V \in \pi^G_{V-\dim_\R V}(X)$ giving an isomorphism
\[ u_V\colon \Sigma^V X \iso \Sigma^{\dim_\R V} X. \]
Such an isomorphism typically only exists when $X$ is Borel-complete; however, if $X$ is a genuine $G$-spectrum and $X^h=X^{EG_+}$ is its Borel completion, the classes $u_V$ (but not $u_V^{-1}$) often lift along $X\to X^h$. For the case we need, the $u_V$ are constructed in \cite[\sec4]{RiggenbachDoublePoints}. Other good sources to learn about these classes are \cite[\sec5.2]{Zeng} and \cite{HZhZ}. Crucially, the $u_V$ are \emph{non-canonical}, so these are the ones we want to get rid of.

  The $n^\th$ power maps $S^{\lambda^k} \to S^{\lambda^{nk}}$ give divisibility relations $a_{\lambda^k} \mid a_{\lambda^{nk}}$, and we set
  \begin{align*}
    a_i &\defeq a_{\lambda_{i-1}}^{-1} a_{\lambda_i}\\
    u_i &\defeq u_{\lambda_{i-1}} u_{\lambda_i}^{-1}
  \end{align*}
  When $i=0$, we interpret these formulas by setting $a_{\lambda_{-1}} \defeq 1$ and $u_{\lambda_{-1}} = \sigma$, the B\"okstedt generator. (There is some interesting philosophy here). Setting $T=\THH(R;\Z_p)$, these classes are summarized in the following diagram:
  \[\xymatrix@=5em{
    & T \ar@/_.5em/[dl]_-{a_{\lambda_{-1}} \defeq \id} \ar[d]^{a_{\lambda_0}} \ar[dr]^{a_{\lambda_1}} \ar@/^1em/[drr]^-{a_{\lambda_2}}\\
      T \ar@<.5ex>[r]^-{a_0} & \Sigma^{\lambda_0} T \ar@<.5ex>[l]^-{u_0} \ar@<.5ex>[r]^-{a_1} & \Sigma^{\lambda_1} T \ar@<.5ex>[l]^-{u_1} \ar@<.5ex>[r]^-{a_2} & \cdots \ar@<.5ex>[l]^-{u_2}\\
                            & \Sigma^2 T \ar@/^.5em/[ul]^-{u_{\lambda_{-1}} \defeq \sigma} \ar[u]_-{u_{\lambda_0}} \ar[ur]_-{u_{\lambda_1}} \ar@/_1em/[urr]_-{u_{\lambda_2}}
  }\]
  Note that the generator $t$ of $\TP_{-2}(R;\Z_p)$ is given by $t=a_{\lambda_0} u_{\lambda_0}^{-1}$.

  The gold classes are related by the ``$q$-gold relations''
  \begin{alignat*}{2}
    \phi^{-1}(\~\xi_{i+1})a_{\lambda_i} &= 0\\
    a_{\lambda_j} u_{\lambda_i} &\approx \phi^i(\~\xi_{j-i}) a_{\lambda_i} u_{\lambda_j} &\quad& (-1\le i<j)\\
    a_i u_i &\approx \phi^i(\xi)\\
    \intertext{from \cite[Lemma 4.9]{SulSliceTHH}, as well as the ``infinite products''}
    u_{\lambda_{i-1}} &\approx u_i u_{i+1} u_{i+2} \dotsm
  \end{alignat*}
  These are up to units since the $u_V$ themselves are only defined up to units. The infinite products can be made precise by viewing $\TF(R;\Z_p) = \lim_F \TR^n(R;\Z_p)$ as a condensed spectrum with the limit topology, or by using Breuil-Kisin twists.

  In \cite{SulROGTHH}, we used the classes
  \[ \vartheta^\alpha_r \defeq \left(\prod_{i=0}^r a_i^{-d_i}\right)\left(\prod_{i=r+1}^\infty u_i^{d_i}\right) \]
  using the above convention for infinite products. More canonically (Notation \ref{not:canonical}), this corresponds to tensoring with the ideal $\bigotimes_{i=r}^\infty \phi^{i}(I)$.
\end{notation}

\begin{remark}[Polar rings]
  It is often said that the multiplication on $\RO(G)$ is a red herring, and that the notion of $\RO(G)$-graded homotopy only sees the abelian group structure of $\RO(G)$. Above, we noticed that this is not entirely true: $n^\th$ power maps give divisibility relations $a_V \mid a_{V^{\tnsr n}}$.

  Bauer \cite{BauerPolar} has defined the notion of a \emph{$p$-polar ring}: an abelian group equipped with an operation ``$x\mapsto x^p$'', but not a full multiplication; in particular, he notes that $p$-typical Witt vectors are naturally defined for such. One could remove the $p$-typical restriction and define \emph{polar rings} as abelian groups equipped with operations ``$x\mapsto x^n$'', and then say that $\RO(G)$-graded homotopy sees $\RO(G)$ as a polar ring.
\end{remark}

\begin{conjecture}[$G$-polar derivations and the $q$-de Rham complex]
  More generally, one could define \emph{naive polar $G$-rings} as modules with $G$-action equipped with a notion of ``$x \mapsto x^{G/H}$'', and \emph{genuine polar $G$-rings} as Mackey functors equipped with norm maps but not a full multiplication. Here is a reason for doing so: Hill \cite{HillTambaraAQ} has defined a cotangent complex for Tambara functors, based on a notion of $G$-derivation. A $G$-derivation is a levelwise derivation which additionally satisfies compatibility with Norm maps, which we state heuristically as
  \[ d(x^{G/H}) = G/H \cdot \text{``$x^{G/H} \dlog x$''}. \]

  We can turn the $q$-de Rham complex into a complex of Mackey functors via $C_n \mapsto \frac{\qOmega^*}{(q^n-1)}$ or $C_n \mapsto \frac{\qOmega^*}{(n)_q}$, and we can turn $\qOmega^0$ into a Tambara functor using $q$-twisted powers \cite[\sec5.3]{PrismsTambara1}. The $q$-derivative then satisfies the above compatibility with norms, but it is not a levelwise derivation. We conjecture that, suitably defined, the cotangent complex of $G$-polar rings (which should remove the derivation condition, but keep the norm condition) recovers $\qOmega^1$; here the group is $G=\Q/\Z$, and perhaps the cotangent complex should be relative to the Tambara functor $\m\RU = \left[C_n \mapsto \frac{\Z[q]}{q^n-1}\right]$.
\end{conjecture}

\begin{notation}[Poly-representations]
  To study $\TR_\rog$ rather than just $\TR^n_\rog$, one switches to using poly-representations, introduced in \cite[\sec4.2]{SulSlopes} in order to conceptualize \cite[Proposition 4.2.3]{HMCyclicPolytopes}.

  In the $\TR$ setting, we denote a poly-representation $\m\alpha$ by $(\m d_0, \m d_1, \m d_2,\dotsc)$, with $\m d_n = \dim_C \m\alpha(C_{p^n})$. We have poly-graded gold classes
  \begin{align*}
    \m a_{\lambda_i} \colon C_{p^n} &\mapsto a_{\lambda_{n-i}}\\
    \m u_{\lambda_i} \colon C_{p^n} &\mapsto u_{\lambda_{n-i}}\\
    \m a_i \colon C_{p^n} &\mapsto a_{n-i}\\
    \m u_i \colon C_{p^n} &\mapsto u_{n-i}
  \end{align*}
  The polygonic analogue of $\vartheta^\alpha_r$ is
  \[ \vartheta^{\m\alpha}_r = \m u_0^{\m d_0}\u_1^{\m d_1}\dotsm\u_r^{\m d_r} \a_{r+1}^{-\m d_{r+1}}\dotsm \]
  which more canonically corresponds to tensoring with $\bigotimes_{i=0}^r \phi^{-i}(\xi)^{\m d_i}A$ (Definition \ref{def:poly-ideal}).
\end{notation}

\begin{warning}
  We should really write $d_{-n}(\m\alpha)$, but this notation gets too bulky; the underline notation emphasizes the connection to Mackey functors, and $\m d_n$ still looks pretty close to $d_{-n}$. However, notational macros like this can obscure the bigger picture: for integral poly-representations, one should write $d_{1/n}(\m\alpha)$.
\end{warning}

In particular, interpreting $\m a_i = a_{-i}$ and $\m u_i = u_{-i}$, the polygonic gold relation is still just $a_i u_i = \phi^i(\xi)$, but for negative $i$. In $\TR_\rog(\O_C;\Z_p)$, we have the additional relation
\[
  \beta a_{\m\lambda_i} = \phi^{-i}(\mu) u_{\m\lambda_i}
\]
where $\beta$ is the Bott element.

We highlight some significant families of poly-representations; these are not used in the current paper, but are important for understanding and used in related works.

\begin{definition}
  The \emph{irreducible} poly-representations are
  \[ \lambda^{1/k}\colon C_n \mapsto \begin{cases} \lambda^{n/k} & k\mid n\\ 0 & k\nmid n\end{cases} \]
  In the $p$-typical case, we set $\m\lambda_i = \lambda^{1/p^i}$. We have $\m d_j(\m\lambda_i) = [j \ge i]$, where $[P]$ is the Iverson bracket.
\end{definition}

\begin{definition}
  \label{def:regular-poly}
  For a positive rational number $\nu\in\Q_{>0}$, the \emph{regular poly-representation of slope $\nu$} is
  \[ \hsl\nu\colon C_n \mapsto [\lc n\nu\rc]_\lambda \defeq \lambda^0 + \dotsb + \lambda^{\lc n\nu\rc-1}. \]
  These are so-named because $\hsl1$ specializes to the complex regular representation of each $C_n$. Beware that the irreducible decomposition of $\hsl\nu$ is infinite. This has dimension sequence $\m d_i(\hsl\nu) = \lc p^i\nu\rc$, or integrally $d_{1/n}(\hsl\nu)=\lc n\nu\rc$.
\end{definition}

The notational evolutions reviewed thus far have already appeared in previous papers. Now we explain the new maneuvers introduced in this paper.

\begin{notation}
  To extend the homotopy orbits to TR spectral sequence from perfectoid rings to smooth algebras over such, the main simplification is to pull out
  \[ \TR^\bullet_{\rog,j}(S;\Z_p) \defeq \text{contribution of $\dRW[\bullet]^j_S$ to $\TR^\bullet_{\rog+j}(S;\Z_p)$.} \]
  That we can do this follows by [the proof of] Hesselholt's theorem \cite{LarsPTyp} (see also \cite[Theorem 11]{HesselholtHandbook})
  \[ \TR^n_*(S;\Z_p) =  \dRW[n]^*_{S/R} \otimes_{\rmW(R)} \TR^n_*(R;\Z_p), \]
  which in this notation becomes
  \[ \TR^n_{2i,j}(S;\Z_p) = \dRW[n]^j_{S/R} \tnsr \sigma_n^i \]
  where $\sigma_n$ is the B\"okstedt generator of $\TR^n_2(R;\Z_p)$. With this notation, the computation for a smooth algebra $S$ works mostly the same as for the base perfectoid ring $R$.
\end{notation}

\begin{notation}[Canonical identifications]
  \label{not:canonical}
  Since the Thom classes $u_{\lambda^i}$ are invertible in $\TCmin$, we have \emph{non-canonical} isomorphisms
  \begin{align*}
    \TCmin_\alpha(S;\Z_p) &\cong \TCmin_{2 d_0(\alpha)}(S;\Z_p)\\
    \TP_\alpha(S;\Z_p) &\cong \TP_{2 d_0(\alpha)}(S;\Z_p)\\
                      &\cong \TP_0(S;\Z_p)
  \end{align*}
  for every $\alpha\in\RU(\T)$. The content of our computation is to provide \emph{canonical} identifications. Since it is one of the main focuses of this paper, and the author was confused about it for several years, let us clarify what ``canonical'' really means.

  Usually we are told it means ``without making choices'', but this gets a bit philosophical. In particular, we want to say that ``$p\Z_p$ is isomorphic to $\Z_p$, but not canonically''. One could try to justify this by arguing that there are $\Z_p^\times$ possible generators, but this falls apart on noting that $p\Z_p = p\N \tnsr_\N \Z_p$ and $|\N^\times|=1$.

  On the other hand, $p\Z_p$ has a natural map to $\Z_p$, and $(p\Z_p \to \Z_p)$ is certainly not isomorphic to $(\Z_p \morph^\id \Z_p)$ as objects of the slice category $\Mod_{/\Z_p}$. In general, let $T$ be a $G$-spectrum, and write a virtual $G$-representation $\alpha$ as $\alpha = i + U - V$, where $U$ and $V$ are actual $G$-representations with $U^G=V^G=0$, and $i\in\Z$. Then the equviariant Euler classes give $\pi_\alpha^G T$ the structure\footnote{Strictly speaking, this depends on the choice of $U$ and $V$; there is an instructive sermon about this in \cite[\sec6]{Adams}, see also \cite{DuggerCoherence,GradedTambara} for a careful discussion. We avoid this issue by fixing the basis $\{\lambda^i\}$ of $\RU(\T)$ and never sharing irreducibles between $U$ and $V$.} of an object in $\mathrm{Span}(\mathrm{Mod})_{\pi^G_i(T)/}$:
  \[ \xymatrix{
    & \ar[dl]_-{a_U} \pi^G_{i+U} T \ar[dr]^-{a_V}\\
    \pi_i^G T && \pi^G_{i+U-V}T
  } \]
  Now by ``canonical'' identifications of the $\RO(G)$-graded groups, at minimum we mean we are going to identify the maps in this span. There are lots more maps coming from the equivariant structure ($F,\,V,\,N,\,R,\,\dots$), and we want to identify these as well.

  We do this not by writing down a bunch of matrices, but by expressing the answer as tensored with various line bundles---a sort of dimensional analysis---so that all the maps are simply ``the obvious ones''. For example, in the $\TR$ case there are $F$ and $R$ maps, which act on the B\"okstedt generator by
  \begin{align*}
    F(\sigma_n) &\approx \sigma_{n-1}\\
    R(\sigma_n) &\approx \phi^{-n}(\xi)\sigma_{n-1}
  \end{align*}
  These can be explained by the natural maps
  \[\xymatrix{
      \xi_{n+1} A \ar@<1ex>[r]^-F \ar@<-1ex>[r]_-R & \xi_n A
  }\]
  where $R$ is the inclusion, and $F(f)=\phi(f)/\~\xi$ is a divided Frobenius. Thus, we can further rewrite Hesselholt's theorem as
  \[ \TR^n_{2i,j}(S;\Z_p) = \dRW[n]^j_{S/R} \tnsr_A \xi_{n+1}^i A \]
  Actually, the formulas for $F$ and $R$ only justify this once we know the base case $\THH_{2i,j}(S;\Z_p) = \Omega^j_{S/R} \tnsr_A \xi^j A$. This follows from the comparison with the cotangent complex, as explained in Example \ref{ex:tor-ext}.

  The module $\xi_nA/\xi_n^2A$ is often denoted by $(A/\xi_n)\{1\}$, but we will need a much larger range of twists than this, so we will instead write this as $(A\tnsr \xi_n)/\xi_n$ or $(A/\xi_n) \tnsr \xi_n$. In general, if $M$ is an $A$-module, a notation like $M \tnsr \xi^2 \phi^{-1}(\xi)^{-3}$ means
  \[ M \tnsr_A \xi^2A \tnsr_A \Hom_A\left(\phi^{-1}(\xi)^3A,A\right) \]

We will sometimes need to make sense of this notation for certain \emph{infinite} products. In the prismatic context, the infinite product $I\phi(I)\dotsm$ can be made sense of as the Breuil-Kisin twist, as in \cite[\sec2.2]{APC}. In the de Rham-Witt context, the symbol $\xi\phi^{-1}(\xi)\dotsm$ means the ideal
  \[ \bigcap_{n=0}^\infty \xi_n \Ainf. \]
  We only consider this intersection when the perfectoid base is spherically complete, in which case this intersection is equal to $(1-q)\Ainf$ by \cite[Lemma 3.23]{BMS1}. We do not know what the intersection looks like in any other cases (besides $\phi^{-r}(\xi)\phi^{-(r+1)}(\xi)\dotsm = \phi^{-r}(1-q)$).

  When $R$ is a perfectoid ring, the canonical identification $\TF_{\lambda}(R;\Z_p) = \xi\Ainf(R)$ was explained to the author by Mike Hill at MSRI in 2019. Unreadable canonical identifications of $\TF_\rog(R;\Z_p)$ and $\TR^\bullet_\rog(R;\Z_p)$ in the positive and negative ranges appeared in the author's PhD thesis. We discovered the $q$-gold relation \cite[Lemma 4.9]{SulSliceTHH} a few months later; this essentially implies the canonical identifications, but it took years for the author to understand this. The first canonical identifications in the generality of quasiregular semiperfectoid rings are due to Riggenbach \cite{RiggenbachTruncated}.
\end{notation}

\begin{example}
  \label{ex:tor-ext}
  Let $A$ be a ring, $f\in A$ a non-zerodivisor, $I=(f)$ the principal ideal generated by $A$, and $\bar A=A/I$. It is well-known that for any $A$-module $M$,
  \begin{align*}
    \Tor^A_1(\bar A, M) &\cong {}_I M \fade{{}\defeq \{m \in M \mid Im = 0\}}\\
    \Ext^1_A(\bar A, M) &\cong M/I
  \end{align*}
  as one sees by resolving $A/I$ as $A\morph^f A$.

  It is less well-known that these isomorphisms are non-canonical: the chosen resolution depends on $f$, while $\bar A$ only depends on $I$. Let $\bar A\{1\} = I/I^2$, and for any $\bar A$-module $N$ let $N\{n\} = N \tnsr_{\bar A}(\bar A\{1\})^{\tnsr n}$. In general, we have
  \begin{align*}
    \gr^i_I A &= \begin{cases}\bar A\{i\} & i\ge 0\\0 & i <0\end{cases} &
    \gr^i_I (A[1/I]) &= \bar A\{i\}
  \end{align*}

  Now a \emph{canonical} free resolution of $\bar A$ is $I\to A$, which gives \emph{canonical} identifications
  \begin{align*}
    \Tor^A_1(\bar A, M) &= {}_I M \{1\}\\
    \Ext^1_A(\bar A, M) &= M/I\{-1\}.
    \intertext{Similarly, there are canonical identifications}
    \mathrm L_{\bar A\mid A} &= \bar A\{1\}[1]\\
    \HH_{2i}(\bar A \rel A) &= \Gamma^i \bar A\{i\}.
  \end{align*}
  In particular, if there is additional structure like a Galois action or a Frobenius floating around, these twists will keep track of the correct actions.

  Please mention this point the next time you teach homological algebra.
\end{example}

\begin{notation}[The power of negative thinking] In both \cite{SulROGTHH} and \cite{SulSlopes}, there is a lot of fussing over whether a given expression is negative or not. Usually this is because we have some filtration like $\xi^i A$, and when $i$ is negative we want to interpret this as $\xi^{\max(i, 0)}A$ (i.e.\ make the filtration constant in negative degrees), rather than the $\xi$-adic filtration on the localization $A[\xi^\pm]$. By interpreting filtrations this way, and similarly defining $\rmW_i(R)=0$ for $i<0$, we can avoid handling $2^\infty-1$ cases in exchange for complicated expressions like $\dRW[-(\m d_1+\m d_2)]^j$ that will simplify in all practical situations.
\end{notation}

We close with a final maneuver which we have not adopted, but could be useful in future.

\begin{notation}[Accepting what we cannot know] In \sec\ref{sec:dRW}, one reason for restricting to the contiguous range is that we do not know how to compute the kernels of maps such as
  \[ \frac A\xi \xra{\phi^{-1}(\xi)^a} \frac A{\xi\phi^{-1}(\xi),\, \phi^{-2}(\xi)^b} \]
  where $a,b\ge0$ are arbitrary. But we could simply give this kernel a name, and compute other things in terms of it, and leave it to the commutative algebraists to give a more explicit description of its structure.
\end{notation}

%% file: qrsp.tex
%
\section{\tops{$\RO(\T)$}{RO(T)}-graded \tops{$\TP$}{TP} and \tops{$\TC^-$}{TC-} of quasiregular semiperfectoid rings}
\label{sec:qrsp}

\newcommand{\vprism}{\prismcpl^{(1)}}

In this section we study the $\RO(\T)$ (really $\RO(C_{p^\infty})$)-graded $\TCmin$ and $\TP$ of quasiregular semiperfectoid rings $S$. We compute the relevant groups in \sec\ref{sub:qrsp-fractional}; the result is essentially the multi-Rees algebra of the Nygaard and $\phi^\bullet(I)$-adic filtrations on $\prismcpl_S[1/\phi^\bullet(I)]$, so that we encounter \emph{prismatic fractional ideals} like $I^2 \tnsr \phi(I)^{-3}$. We also give integral generalizations of some of our results to even $\E_\infty$-rings with $\T$-action. These computations ignore the distinction between $\lambda^n$ and $\lambda^{p^{v_p(n)}}$; in \sec\ref{sub:integral} we explain how to put this distinction back in, which does not affect module structures but \emph{does} twist restriction ($F$) maps, and becomes important when studying norms. Finally, in \sec\ref{sub:iterated-nygaard} we extend our computations to the spectra $\TC^{-n} \defeq (\TR^n)^{h\T/C_{p^n}}$.

\subsection{Prismatic fractional ideals}
\label{sub:qrsp-fractional}
In this section we explain how to see Theorem \ref{thm:qrsp}. We fix a $\delta$-ring $A$ with spherical lift $\S_A$, a relatively quasiregular semiperfectoid $A$-algebra $S$, and abbreviate $\TCmin\defeq\TCmin(S\rel\S_A)$, $\TP\defeq\TP(S\rel\S_A)$, $\vprism \defeq\vprism_{S/A}$, etc.

Our main tools will be the cell structure
\begin{equation}
  \label{eq:cell}
  \T/C_{n+} \to S^0 \to S^{\lambda^n}
\end{equation}
and the motivic filtrations
\begin{align*}
  \TCmin_{2i} &= \Nyg^{\ge i}\vprism\{i\} & \TP_{2i} &= \vprism\{i\}\\
  \THH^{hC_{p^n}}_{2i} &= \frac{\Nyg^{\ge i}\vprism\{i\}}{I_n\tnsr\Nyg^{\ge i+1}\vprism\{i\}} & \THH^{tC_{p^n}}_{2i} &= \frac{\vprism\{i\}}{I_n\vprism\{i\}}
\end{align*}
from \cite[Theorem 1.12(4)]{BMS2} and \cite[Corollary 3.3]{RiggenbachTruncated}. We will just explain the case of $\TCmin$, leaving the (easier) case of $\TP$ to the reader.

\begin{remark}
  \label{rmk:no-extns}
  Since $\THH$ of perfectoid rings is even, there are invertible classes $u_{\lambda_k}\in\TCmin_{2-\lambda_k}$ giving non-canonical isomorphisms
  \[ \TCmin_\alpha \cong \TCmin_{2d_0(\alpha)} \]
  for every $\alpha\in\RU(\T)$ \cite[Proposition 4.2]{RiggenbachDoublePoints}. Of course the whole point of this section is to avoid these classes and instead give canonical identifications. However, this observation is useful for solving extensions and identifying images.
\end{remark}

\begin{example}
  \label{ex:qrsp-pos}
  The cell structure \eqref{eq:cell} gives a short exact sequence
  \[\xymatrix{
      0 \ar[r] & \TCmin_{2i+\lambda_k} \ar[r]^-{a_{\lambda_k}} \ar@{=}[d] & \TCmin_{2i} \ar[r] \ar@{=}[d] & \THH^{hC_{p^k}}_{2i} \ar[r] \ar@{=}[d] & 0\\
      0 \ar[r] & ? \ar[r] & \Nyg^{\ge i}\vprism\{i\} \ar[r] & \dfrac{\Nyg^{\ge i}\vprism\{i\}}{I_k\tnsr\Nyg^{\ge i+1}\vprism\{i\}} \ar[r] & 0
  }\]
  which shows that $\TCmin_{2i+\lambda_k}=I_k\tnsr\Nyg^{\ge i+1}\vprism\{i\}$, as claimed by Theorem \ref{thm:qrsp}.
\end{example}

\begin{example}
  The cell structure \eqref{eq:cell} gives a short exact sequence
  \[\xymatrix{
      0 \ar[r] & \TCmin_{2i} \ar[r]^-{a_{\lambda_k}} \ar@{=}[d] & \TCmin_{2i-\lambda_k} \ar[r] \ar@{=}[d] & \THH^{hC_{p^k}}_{2(i-1)} \ar[r] \ar@{=}[d] & 0\\
      0 \ar[r] & \Nyg^{\ge i}\vprism\{i\} \ar[r] & ? \ar[r] & \dfrac{\Nyg^{\ge i-1}\vprism\{i-1\}}{I_k\tnsr\Nyg^{\ge i}\vprism\{i-1\}} \ar[r] & 0
  }\]
  Rewriting the final term as
  \[
    \dfrac{\Nyg^{\ge i-1}\vprism\{i-1\}}{I_k\tnsr\Nyg^{\ge i}\vprism\{i-1\}} 
    =
    \dfrac{I_k^{-1}\tnsr\Nyg^{\ge i-1}\vprism\{i\}}{\Nyg^{\ge i}\vprism\{i\}} 
  \]
  and using Remark \ref{rmk:no-extns}, we see that $\TCmin_{2i-\lambda_k} = I_k^{-1}\tnsr\Nyg^{\ge i-1}\vprism\{i\}$, as claimed by Theorem \ref{thm:qrsp}.

  Alternatively, we can use the multiplication map
  \[\xymatrix{
      \TCmin_{2i-\lambda_k} \tnsr \TCmin_{\lambda_k-2} \ar[r] \ar@{=}[d] & \TCmin_{2(i-1)} \ar@{=}[d]\\
      \TCmin_{2i-\lambda_k} \tnsr I_k\tnsr\vprism\{-1\} \ar[r] & \Nyg^{\ge i-1}\vprism\{i-1\}
  }\]
  where the left vertical identification is Example \ref{ex:qrsp-pos}. Using Remark \ref{rmk:no-extns}, we see that the rows are isomorphisms. We thus get
  \[ \TCmin_{2i-\lambda_k} = \Nyg^{\ge i-1}\vprism\{i-1\}\tnsr\left(I_k\tnsr\vprism\{-1\}\right)^{-1} = I_k^{-1}\Nyg^{\ge i-1}\vprism\{i\}. \]
\end{example}

\begin{example}
  We now explain how to mix representations. If $k\le n$, then $\res^\T_{C_{p^k}} \lambda_n$ is trivial, but $\res^\T_{C_{p^n}} \lambda_k$ is not unless $k=n$. We will thus induct by using $\lambda_n$ before we use $\lambda_k$. Now the cell structure \eqref{eq:cell} gives a short exact sequence
  \[\xymatrix{
      0 \ar[r] & \TCmin_{2i+\lambda_k-\lambda_n} \ar[r]^-{a_{\lambda_k}} \ar@{=}[d] & \TCmin_{2i-\lambda_n} \ar[r] \ar@{=}[d] & \THH^{hC_{p^k}}_{2(i-1)} \ar[r] \ar@{=}[d] & 0\\
      0 \ar[r] & ? \ar[r] & I_n^{-1}\Nyg^{\ge i-1}\vprism\{i\} \ar[r] & \dfrac{\Nyg^{\ge i-1}\vprism\{i-1\}}{I_k\tnsr\Nyg^{\ge i}\vprism\{i-1\}} \ar[r] & 0.
  }\]
  Rewriting the final term as
  \[
    \dfrac{\Nyg^{\ge i-1}\vprism\{i-1\}}{I_k\tnsr\Nyg^{\ge i}\vprism\{i-1\}}
    =
    \dfrac{I_k^{-1}\tnsr\Nyg^{\ge i-1}\vprism\{i\}}{\Nyg^{\ge i}\vprism\{i\}}
    =
    \dfrac{I_n^{-1}\tnsr\Nyg^{\ge i-1}\vprism\{i\}}{I_k\tnsr I_n^{-1}\tnsr\Nyg^{\ge i}\vprism\{i\}},
  \]
  we see that $\TCmin_{2i+\lambda_k-\lambda_n}=I_k\tnsr I_n^{-1}\tnsr\Nyg^{\ge i}\vprism\{i\}$, as claimed by Theorem \ref{thm:qrsp}.

  Alternatively, we can use the multiplication map
  \[\xymatrix{
      \TCmin_{2i+\lambda_k-\lambda_n} \tnsr \TCmin_{\lambda_n-2} \ar[r] \ar@{=}[d] & \TCmin_{2(i-1)+\lambda_k} \ar@{=}[d]\\
      \TCmin_{2i+\lambda_k-\lambda_n} \tnsr I_n\tnsr\vprism\{-1\} \ar[r] & I_k\tnsr\Nyg^{\ge i}\vprism\{i-1\}
  }\]
  where the left vertical identification is Example \ref{ex:qrsp-pos}. Using Remark \ref{rmk:no-extns}, we see that the rows are isomorphisms. We thus get
  \[ \TCmin_{2i+\lambda_k-\lambda_n} = I_k\tnsr\Nyg^{\ge i}\vprism\{i-1\}\tnsr\left(I_n\tnsr\vprism\{-1\}\right)^{-1} = I_k\tnsr I_n^{-1}\tnsr\Nyg^{\ge i}\vprism\{i\}. \]
\end{example}

Combining the techniques in the above examples, the general case of Theorem \ref{thm:qrsp} follows immediately.

\begin{remark}
  \label{rmk:even}
  The only inputs needed for Theorem \ref{thm:qrsp} were:
  \begin{itemize}
    \item $\THH(S\rel\S_A)$ is orientable in the Borel category \cite[Proposition 4.2]{RiggenbachDoublePoints};
    \item a description of the Mackey functors $\mpi_\rog\THH(S\rel\S_A)^h$ and $\mpi_\rog\THH(S\rel\S_A)^t$ \cite[Corollary 3.3]{RiggenbachTruncated}.
  \end{itemize}
  Riggenbach in fact shows the first point for all $\T$-equivariant even $\E_\infty$-rings, while he and Antieau provide a version of the second in this generality in \cite[Corollary 3.20]{CycSyn}. We now explain this integral generalization (Theorem \ref{thm:even}).

  Let $E$ be an even $\E_\infty$-ring with $\T$-action. Set $A = \pi_0 E^{t\T}$, $A\{1\} = \pi_2 E^{t\T}$, and $A\{n\} = A\{1\}^{\tnsr n}$; by the Tate spectral sequence, we have
  \[ \pi_* E^{t\T} = \bigoplus_{i\in\Z} A\{i\}[2i] \]

  Let $\Fil^\bullet A\{i\}$ be the abutment filtration of the Tate spectral sequence. By the homotopy fixed points spectral sequence, we have
  \[ \pi_* E^{h\T} = \bigoplus_{i\in \Z} \Fil^i A\{i\}[2i] \]
  The natural map
  \[
    \Fil^i A\{i\} \tnsr_A A\{-i\} \to \Fil^i A
  \]
  induces an isomorphism
  \[ \pi_{2i} E = \gr^i A\{i\} \iso \gr^i A \]
  and therefore the short exact sequences
  \[ 0 \to E^{h\T}_{2i+(j+1)\lambda}\xra{a_\lambda} E^{h\T}_{2i+j\lambda} \to \pi_{2(i+j)} E \to 0 \]
  give identifications
  \[ E^{h\T}_{2i+j\lambda} = \Fil^{i+j}A\{i\} \]
  for all $i,j\in\Z$, as before.

  To go further, let $\G$ be the formal group law of $E$, and let $t\in A\{-1\}$ be a choice of complex orientation. Write $[n]_\G(t)\in A\{-1\}$ for the $n$-series of $\G$, $\<n\>_\G(t) = \frac{[n]_\G(t)}{t} \in A$ for the divided $n$-series; we sometimes suppress the dependence on $t$ and just write $\<n\>_\G$ (since we only care about the ideal $\<n\>_\G(t)A$). Assume furthermore that $\<n\>_\G(t)$ is a non-zerodivisor for all $n$. By \cite[Corollary 3.20]{CycSyn}, we have
  \begin{align*}
    \pi_{2i} E^{hC_n} &= \frac{\Fil^i A\{i\}}{\<n\>_\G\Fil^{i+1}A\{i\}} &
    \pi_{2i}E^{tC_n} &= \frac{A\{i\}}{\<n\>_\G}
  \end{align*}
  It follows that
  \begin{align*}
    E^{h\T}_{2i+j\lambda^n} &= \Fil^{i+j} A\{i\} \tnsr_A (\<n\>_\G A)^{\tnsr j}\\
    E^{t\T}_{2i+j\lambda^n} &= A\{i\} \tnsr_A (\<n\>_\G A)^{\tnsr j}
  \end{align*}

  To go fully general, we need a way to describe the geometric fixed points. Recall $\Phi_{d,\G}(t)$ from Definition \ref{def:fgl-cyclotomic}. In our topological situation, we have $\Phi_{d,\G}(t) \in A$ for $d>1$, and we suppress the dependence on $t$ and simply write $\Phi_{d,\G}\in A$. Note that $\Phi_{n,\G}=0$ in $\pi_0 E^{\tau C_n}$ (proper Tate construction), and heuristically $A\{1\} = \bigotimes_{n>1} \Phi_{n,\G}A$. With these notations, Theorem \ref{thm:even} follows immediately.
\end{remark}

\begin{example}
  Let $E$ be an even $\E_\infty$-ring with $\T$-action, and assume that $\<n\>_\G \in A=\pi_0 E^{t\T}$ is a non-zerodivisor for all $n$. Then for all $i\in\Z$,
  \begin{align*}
    E^{h\T}_{(n)_{\lambda^m} + 2i}
      &= \Fil^{i+n} A\{i+1\} \tnsr \<n\>_{\G^{(m)}}? &
    E^{t\T}_{(n)_{\lambda^m} + 2i}
      &= A\{1+i\} \tnsr \<n\>_{\G^{(m)}}?\\
    E^{h\T}_{\lambda^m(n)_{\lambda^m} + 2i}
      &= \Fil^{i+n} A\{i\} \tnsr \<n\>_{\G^{(m)}}! &
    E^{t\T}_{\lambda^m(n)_{\lambda^m} + 2i}
      &= A\{i\} \tnsr \<n\>_{\G^{(m)}}!
  \end{align*}
  From here we can get more complicated examples like
  \[
    E^{h\T}_{\psi^m(\lambda((n)_\lambda - (k)_\lambda - (n-k)_\lambda)) + 2i}
      = \Fil^i A\{i\} \tnsr {\binom nk}_{\G^{(m)}}
  \]
\end{example}

\subsection{Integral representations}
\label{sub:integral}

Above, we have restricted to $\T$-representations of the form $\lambda^0$ and $\lambda^{p^i}$. For most purposes this is fine, since the divided Euler class $S^{\lambda^{p^n}} \to S^{\lambda^n}$ becomes an isomorphism after $p$-completion. However, for pinning down units we need to \emph{keep track of this isomorphism}. In particular, need this when identifying the norm map: for example, there is a norm map
\[ \pi_2 \THH \morph^N \THH^{hC_p}_{\ind^{C_p}_e 2} = \THH^{hC_p}_{(p)_\lambda} \cong \THH^{hC_p}_{\lambda^0 + (p-1)\lambda^1} \]
where $(p)_\lambda = \lambda^0 + \lambda^1 + \dotsb + \lambda^{p-1}$. The stated isomorphism multiples the $F$ map by $(p-1)!$, which must be taken into account when considering the axiom $FNx=x^p$.

We recall the following theorem of Angeltveit.

\begin{definition}[{\cite[Definition 4.5]{AngeltveitPicSpCn}}]
  Let $\ds a=(a_d)\in\prod_{1\ne d\mid n} \Z$ with each $a_d\ne0$. Let ${}^a \m W(\Z)$ denote the Mackey functor defined in the same way as the Burnside Mackey functor $\m W(\Z)$, but with restriction maps determined by
  \[ F^n_{n/k}(C_n/C_n) = \left(\prod_{d\mid k} a_d\right)C_{n/k}/C_{n/k} \]
\end{definition}

\begin{example}
  For $n=p$, the $C_p$-Mackey functor ${}^a \m W(\Z)$ is given by
  \[ \Lewis{\Z^2}{\Z}{\begin{bmatrix}a&p\end{bmatrix}}{\begin{bmatrix}0\\1\end{bmatrix}} \]
  for $a\in\{1,\dotsc,p-1\}$. This is a transversal Mackey functor, so it can equivalently be described by the data
  \[\xymatrix{
      \fade{\Z^2} \ar@[lightgray][r]^-{\fade R} \ar@[lightgray][d]_-{\fade F} \ar@{}[dr]|<<{\mbox{\fade{$\lrcorner$}}} & \Z \ar[d]^-a \\
      \Z \ar[r]_-\can & \F_p
  }\]
\end{example}

\begin{example}
  The $C_{p^2}$-Mackey functor ${}^{(a_p,a_{p^2})} \m W(\Z)$ is given by
  \[
    \Lewiss
      {\Z^3}
      {\Z^2}
      {\Z}
      {\begin{bmatrix}a_p&p&0\\0&0&p\end{bmatrix}}
      {\begin{bmatrix}a_{p^2} & p\end{bmatrix}}
      {\begin{bmatrix}0&0\\1&0\\0&1\end{bmatrix}}
      {\begin{bmatrix}0\\1\end{bmatrix}}
  \]
  for $a\in\{1,\dotsc,p-1\}$. This is a transversal Mackey functor, so it can equivalently be described by the data
  \[\xymatrix{
      \fade{\Z^2} \ar@[lightgray][r]^-{\fade R} \ar@[lightgray][d]_-{\fade F} \ar@{}[dr]|<<{\mbox{\fade{$\lrcorner$}}} & \Z \ar[d]^-a \\
      \Z \ar[r]_-\can & \F_p
  }\]
\end{example}

\begin{theorem}[{\cite[Theorem 5.1]{AngeltveitPicSpCn}}]
  Let $X$ be the $C_n$-spectrum $S^{\lambda^{n/d} - \lambda^{an/d}}$ with $d\mid n$ and $(a,d)=1$. Then $X$ is equivalent to the Mackey functor ${}^{a} \m W(\Z)$ with $a_e=\begin{cases}a&e=d\\ 1&\text{else.}\end{cases}$
\end{theorem}

\begin{definition}
  \label{dfn:prismatic-factorial}
  Let $(A, I)$ be a prism, and let $n>0$. Write $n = p^v j$ with $(p,j)=1$. We define the symbol $(n)_\prism A$ as follows: as an $A$-module, it is the ideal $I_v = I \phi(I) \dotsm \phi^{v-1}(I)$. Using the Tambara functor $\m{(A, I)}$ constructed in \cite{PrismsTambara1}, there is a natural Mackey functor structure on $\dfrac{I_v}{I_v \tnsr I_\bullet} = I_v \tnsr_A \m{(A, I)}$. The difference between $(n)_\prism A$ and $I_v$ is that the Mackey functor $\dfrac{(n)_\prism A}{(n)_\prism A \tnsr I_\bullet}$ has its $C_{p^{v+1}} \to C_{p^v}$ restriction map multiplied by $j$. To understand this formula, note that in the $q$-de Rham case we have $(n)_q = (p^v)_q (j)_{q^{p^v}}$, so $(j)_{q^{p^v}} \equiv j \in A/(q^{p^v}-1)$. Since this is natural in the prism $(A, I)$, we often just write $(n)_\prism$.

  We similarly define $(n)_{\prism^m} = (mn)_\prism \tnsr (m)_\prism^{\tnsr -1}$, $(\Phi_n)_\prism = \bigotimes_{d\mid n} {(n)_\prism ^{\tnsr\mu(n/d)}}$, $(n)_\prism! = \bigotimes_{k=1}^n {(k)_\prism}$ and $(n)_\prism? = (n-1)_\prism!$.
  Over the $q$-de Rham prism, these correspond to the ideals generated by $(n)_{q^m}$, $\Phi_n(q)$, $(n)_q!$, and $\Gamma_q(n)$ respectively.
\end{definition}

\begin{proposition}
  \renewcommand{\int}{\mathrm{int}}
  Let $\alpha$ be an integral virtual $\T$-representation with irreducible decomposition $\alpha = \sum_{i=0}^\infty k_i \lambda^i$ and dimension-sequence $d_n = \dim_\C \alpha^{C_n}$ (with $d_\infty=k_0$). If $A$ is a $\delta$-ring with spherical lift $\S_A$, and $S$ is a relatively quasiregular semiperfectoid $A$-algebra, then there are canonical identifications
  \begin{align*}
    \TCmin_\alpha(S\rel \S_A)
    &= \Nyg^{\ge d_1} \prismcpl^{(1)}_{S/A}\{k_0\} \tnsr \bigotimes_{n=1}^\infty {(n)_\prism^{\tnsr k_n}}
    &
    \TP_\alpha(S\rel\S_A)
    &= \prismcpl^{(1)}_{S/A}\{k_0\} \tnsr \bigotimes_{n=1}^\infty {(n)_\prism^{\tnsr k_n}}
    \\
    &=\Nyg^{\ge d_1} \prismcpl^{(1)}_{S/A} \tnsr \bigotimes_{n=1}^\infty {(\Phi_n)_\prism^{\tnsr d_n}}
    &
    &= \prismcpl^{(1)}_{S/A} \tnsr \bigotimes_{n=1}^\infty {(\Phi_n)_\prism^{\tnsr d_n}}
  \end{align*}
\end{proposition}

\subsection{Iterated Nygaard filtrations}
\label{sub:iterated-nygaard}

\allowdisplaybreaks

The computation of $\TCmin_\rog$ and $\TP_\rog$ was simple and elegant, but somewhat boring. We can get more interesting invariants by going ``more genuine'' and considering the spectra $\TC^{(-n)} \defeq (\TR^n)^{h\T/C_{p^n}}$. To avoid confusing ourselves with Frobenius twists, we will only consider absolute prismatic cohomology; but generalizing to the relative case is just a matter of notation. (Then again, many of the advances in this paper are ``just'' notation.)

This subsection is more in the spirit of \sec\ref{sec:dRW} and can be read after that. 

\begin{definition}[Iterated Nygaard filtrations]
  \label{rmk:iterated-nygaard}
  Define \emph{higher topological negative cyclic homology} by
  \[ \TC^{-r} = (\TR^r)^{h\T/C_{p^r}} \]
  This is equivalently the pullback 
  \[ \TCmin\times_\TP \dotsb \times_\TP \TCmin, \]
  where the structure maps are $\can\colon \TCmin\to \TP$ on the right and $\varphi\colon \TCmin\to\TP$ on the left. In particular we have $\TCmin=\TC^{-0}$, while morally $\TP=\TC^{-(-1)}$. Andriopoulos studied these spectra in his thesis \cite{AndriopoulosHigherNygaard}.
\end{definition}

There are several motivations for considering these spectra:
\begin{enumerate}
  \item Nygaard's original paper \cite{Nygaard} considered not just $\Nyg^{\ge i}\dRW=\varphi^{-1}(p^i\dRW)$, but also
    \[ \Nyg^{\ge i}_r \dRW \defeq \varphi^{-r}(p^{ir}\dRW). \]
    The modern prismatic version would be
    \[ \Nyg^{\ge i}_r \prism^{(r)}_{S/A} = \varphi^{-r}(I_r^i\prism_{S/A}) \]
    where $I_r = I\phi(I)\dotsm\phi^{r-1}(I)$; in particular, the elements $\xi_r=\xi\phi^{-1}(\xi)\dotsm\phi^{-(r-1)}(\xi)$ from \cite[Lemma 3.12]{BMS1} appear naturally in this context. When $R$ is a perfectoid ring, we have
    \[ \TC^{-r}_*(R;\Z_p) = \frac{\Ainf(R)[\sigma, t_r]}{\sigma t_r - \xi_{r+1}} \quad \text{ with }|\sigma|=2,\, |t_r|=-2 \]
    by \cite[Proposition 3.4]{AndriopoulosHigherNygaard}.

  \item \label{enum:tc-r} notably, Andriopoulos introduces the idea of viewing $\TR^r$ as the quotient $\TC^{-r}/t_r$, generalizing the identification $\TCmin/t = \THH$ \cite[Proposition 3.6]{AndriopoulosHigherNygaard}. Just as $\TCmin$ is easier to understand than $\THH$ (since it has less torsion), $\TC^{-r}$ and $\TR^{h\T}$ are easier to understand than $\TR^r$ and $\TR$, and this should clarify the passage between the cyclotomic spectra perspective and the topological Cartier module/polygonic spectra perspective \cite{TCart,KMNPolygonic}.

  \item it is hard to understand what the norm
    \[ N\colon \THH^{hC_{p^n}} \to \THH^{hC_{p^{n+1}}} \]
    is really doing arithmetically, since the finite-level homotopy groups have a lot of torsion. If we could take $(-)^{h\T}$ of the norm $N\colon \THH \to \TR^1$, we might be able to get a norm $N\colon \TCmin \to \TC^{-1}$, which should be easier to understand. However, it is unclear if this works.

  \item it is natural to ask if $\RO(\T)$-grading makes sense for $\TC$. The answer turns out to be no, but you'd stumble upon $\TC^{-r}$ while trying. (However, \cite[Corollary 1.2]{DMCurves} suggests that it might be possible to define \emph{polyrepresentation}-graded $\TC$; we do not explore this in this paper.)

  \item (integral versions of) these spectra arise when constructing Habiro cohomology via refined topological Hochschild homology: see \cite[1.23]{RefinedTC-} and \cite[Theorem 6.5]{HabiroCohomology}. Relatedly, they arise in Devalapurkar-Hahn-Raksit-Yuan's work on prismatization of ring spectra \cite{PrismatizationRingSpectra} for the purpose of decompleting the Nygaard filtration, see the discussion in \cite[Chapter 7]{DevalapurkarThesis}. 

  \item Lahoti and Manam study the stacks associated to $\TR^{h\T}$ \cite[\sec1.3]{LMCohomologyModuliRingStacks}, which they call the stack of \emph{polyfiltered Cartier-Witt divisors}. A variant in positive characteristic is studied by Zhang in \cite{StackyDRW}.
\end{enumerate}

\begin{notation}
  If $\alpha=(d_0,d_1,\dotsc)$, recall that $\alpha'=\alpha^{C_p}=(d_1,d_2,\dotsc)$ and $\alpha^{(r)}=\alpha^{C_{p^r}}=(d_r,d_{r+1},\dotsc)$. To ease notation, set
  \begin{align*}
    \varphi^{-n}(I^\alpha\prismcpl_S)
      &= \Nyg^{\ge d_0,\dotsc,d_{n-1}}_n I^{d_n} \phi(I)^{d_{n+1}}\dotsm\prismcpl_S\\
      &= \left\{f \in \prismcpl_S\{d_\infty\} \,\middle|\, \varphi^k(f) \in \left(\prod_{i=0}^\infty \phi^i(I)^{d_{n-k+i}-d_\infty}\right)\prismcpl_S\{d_\infty\} \text{ for all } 0\le k\le n\right\}\\
      &= \{f\in\prismcpl_S\{d_\infty\} \mid \varphi^k(f) \in I^{\alpha^{(n-k)}}\prismcpl_S \text{ for all }0\le k\le n\}
  \end{align*}
  We also write
  \[ \~\varphi^n(\Nyg^{\ge d_0,\dotsc,d_{n-1}}) \defeq \frac{\varphi^n}{I^{d_0}\dotsm \phi^{n-1}(I)^{d_{n-1}}}(\Nyg^{\ge d_0,\dotsc,d_{n-1}}) \]
  for the image of a divided Frobenius (the map itself depends on orienting $I$, but its image does not).
\end{notation}

\begin{example}
  \label{ex:higher-tc-min}
  First we explain how to compute $\pi_\rog$ of $\TC^{-1}=(\THH^{C_p})^{h\T/C_p}$. Write $\TC^{-,\can}=\THH^{h\T}$ for the copy of $\TCmin$ that maps to $\TP$ via the canonical map, and $\TC^{-,\varphi}=(\THH^{\Phi C_p})^{h\T/C_p}$ for the copy of $\TC^-$ that maps to $\TP$ via Frobenius. Then
  \begin{align*}
    \TC^{-,\can}_\alpha(S;\Z_p) &= \varphi^{-1}(I^\alpha\prismcpl_S) = \Nyg^{\ge d_0}I^{\alpha'}\prismcpl_S\\
    \TP_\alpha(S;\Z_p) &= I^{\alpha'}\prismcpl_S\\
    \TC^{-,\varphi}_\alpha(S;\Z_p) &= \varphi^{-1}(I^{\alpha'}\prismcpl_S) = \Nyg^{\ge d_1}I^{\alpha''}\prismcpl_S
  \end{align*}
  The spectral sequence computing $\TC^{-1}_{\alpha-*}(S;\Z_p)$, where $*\in\{0,1\}$, is
  \begin{align*}
    E_0 &=\dbox{$\vcenter{\xymatrix@R=0.2em{
          \TP_\alpha(S;\Z_p) & \TC^{-,\can}_\alpha(S;\Z_p) \ar[l]\\
                       & \TC^{-,\varphi}_\alpha(S;\Z_p)
                   }}$}\\
       &=
    \dbox{$\vcenter{\xymatrix@R=0.2em{
  I^{\alpha'}\prismcpl_S & \varphi^{-1}(I^\alpha\prismcpl_S) \ar[l]_-\can\\
                          & \varphi^{-1}(I^{\alpha'}\prismcpl_S)
                                                               }}$}\\
    E_1 &=
    \dbox{$\vcenter{\xymatrix@R=0.2em{
      \dfrac{I^{\alpha'}\prismcpl_S}{\varphi^{-1}(I^\alpha)} & \\
      & \varphi^{-1}(I^{\alpha'}\prismcpl_S) \ar[ul]_-\varphi
                                                               }}$}\\
    E_\infty &=
    \dbox{$\vcenter{\xymatrix@R=0.2em{
  \dfrac{I^{\alpha'}\prismcpl_S}{\varphi^{-1}(I^\alpha) + \varphi(\varphi^{-1}(I^{\alpha'}))} & \\
      & \varphi^{-2}(I^\alpha\prismcpl_S)
                                    }}$}\\
    &=
    \dbox{$\vcenter{\xymatrix@R=0.2em{
  \dfrac{\prismcpl_S}{\Nyg^{\ge d_0} + \~\varphi(\Nyg^{\ge d_1})} \tnsr I^{\alpha'}& \\
      & \Nyg^{\ge d_0,d_1}_2 \prismcpl_S \tnsr I^{\alpha''}
    }}$}
  \end{align*}
  We conclude that
  \begin{align*}
    \TC^{-1}_\alpha(S;\Z_p)
      &= \Nyg^{\ge d_0,d_1}_2 \prismcpl_S \tnsr I^{\alpha''}\\
    \TC^{-1}_{\alpha-1}(S;\Z_p)
      &= \dfrac{\prismcpl_S}{\Nyg^{\ge d_0} + \~\varphi(\Nyg^{\ge d_1})} \tnsr I^{\alpha'}
  \end{align*}
\end{example}

For later use, it will be helpful to rewrite the second group as the iterated pushout of
\[\xymatrix{
    & \Nyg^{\ge d_0}\prismcpl_S \tnsr I^{d_1} \ar[d]_-\can \ar[r] & 0\\
    \Nyg^{\ge d_1}\prismcpl_S \ar[r]^-\varphi \ar[d] & \prismcpl_S \tnsr I^{d_1}\\
    0
}\]
(which should further be tensored with $\phi(I^{\alpha''})$, but we omit that since it is not the interesting part).

\begin{example}
  Let $R = \Z_p[\zeta_{p^\infty}]^\land_p$, and let $S=(R[x^{1/p^\infty}]/x)^\land_p$. Let $A=\Ainf(R) = \Z_p[q^{1/p^\infty}]^\land_{(p,\~\xi)}$ where $\~\xi=(p)_q$. By \cite[\sec12.2]{Prismatic}, we have
\[
\prism_S = \bigoplusp_{\nu \in \N[1/p]} A \cdot \frac{y^\nu}{\lf \nu\rf_q!}
\]
where $y=x^{1/p}$ and the completion is $(p,\~\xi)$-adic. By \cite[Lemma 12.7(1)]{Prismatic}, the Nygaard filtration is
\[
  \Nyg^{\ge i} \prism_S = \bigoplusp_\nu \xi^{\pos{i - \lf\nu\rf}} A \cdot \frac{y^\nu}{\lf\nu\rf_q!}
\]
with $\xi=\phi^{-1}(\~\xi) = (p)_{q^{1/p}}$ and $\pos x \defeq \max(x, 0)$. Using Legendre's formula to identify $\~\varphi(\Nyg^{\ge d_1}\prism_S)$, we find that
\begin{align*}
  \Nyg^{\ge d_0, d_1} \prism_S &=
  \bigoplusp_\nu \phi^{-1}(\xi)^{\pos{d_0 - \lf\nu\rf}} \xi^{\pos{d_1 - \lf\nu\rf}} A \fade{{}\cdot \frac{y^\nu}{\lf\nu\rf_q!}}\\
  \dfrac
    {\prism_S}
    {\Nyg^{\ge d_0} \prism_S + \~\varphi(\Nyg^{\ge d_1} \prism_S)}
  &=
  \bigoplusp_{d_1 \le \lf\nu\rf \le d_0} A/(\xi^{\pos{d_0 - \lf\nu\rf}}, \phi(\xi)^{\pos{\lf\nu\rf - d_1}} ) \fade{{}\cdot \frac{y^\nu}{\lf\nu\rf_q!}}
\end{align*}
\end{example}

\begin{example}
  Next we compute $\TC^{-2}_\rog(S; \Z_p)$; everything below should be tensored with $I^{\alpha^{(3)}}$, but for brevity we omit that until the end.
\begin{align*}
  &\phantom{{}\leadsto}\dbox{$\vcenter{\xymatrix@C=3em@R=2em{
  \prismcpl_S\otimes I^{d_1}\phi(I)^{d_2}
    & \Nyg^{\ge d_0}\prismcpl_S\otimes I^{d_1}\phi(I)^{d_2}
       \ar[l]_{\can}\\
  \prismcpl_S\otimes I^{d_2}
    & \Nyg^{\ge d_1}\prismcpl_S\otimes I^{d_2}
       \ar[l]_{\can}\\
  & \Nyg^{\ge d_2}\prismcpl_S
}}$}\\ 
&\leadsto
\dbox{$\vcenter{\xymatrix@C=3em@R=2em{
  \dfrac{\prismcpl_S}
        {\Nyg^{\ge d_0}\prismcpl_S} \tnsr I^{d_1}\phi(I)^{d_2} &\\
  \dfrac{\prismcpl_S}
        {\Nyg^{\ge d_1}\prismcpl_S} \tnsr I^{d_2} &\\
  & \Nyg^{\ge d_2}\prismcpl_S \ar[ul]^\varphi
}}$}\\[1em]
&\leadsto \dbox{$\vcenter{
\xymatrix@C=3em@R=2em{
  \dfrac{\prismcpl_S}
        {\Nyg^{\ge d_0}\prismcpl_S} \tnsr I^{d_1}\phi(I)^{d_2} &\\
  \dfrac{\prismcpl_S}
        {\Nyg^{\ge d_1}\prismcpl_S
          +\~\varphi(\Nyg^{\ge d_2}\prismcpl_S)} \tnsr I^{d_2} &\\
  & \Nyg^{\ge d_1,d_2}\prismcpl_S \ar@/_1em/[uul]_-{\varphi^2}
}}$}\\[1em]
&\leadsto
\dbox{$\vcenter{\xymatrix@C=3em@R=2em{
  \dfrac{\prismcpl_S}
        {\Nyg^{\ge d_0}\prismcpl_S
          +\~\varphi^2(\Nyg^{\ge d_1,d_2}\prismcpl_S)} \tnsr I^{d_1}\phi(I)^{d_2}
  \fade{{}\tnsr\phi^2(I^{\alpha^{(3)}})} &\\
  \dfrac{\prismcpl_S}
        {\Nyg^{\ge d_1}\prismcpl_S
        +\~\varphi(\Nyg^{\ge d_2}\prismcpl_S)} \tnsr I^{d_2}
  \fade{{}\otimes\phi(I^{\alpha^{(3)}})} &\\
  & \Nyg^{\ge d_0,d_1,d_2}\prismcpl_S
  \fade{{}\tnsr I^{\alpha^{(3)}}}
    }}$}
\end{align*}
The extension is the iterated pushout of
\[
\xymatrix@R=2em{
  && \Nyg^{\ge d_0}\prismcpl_S\otimes I^{d_1}\phi(I)^{d_2}
      \ar[r] \ar[d]_{\can} & 0\\
  & \Nyg^{\ge d_1}\prismcpl_S\otimes I^{d_2}
      \ar[r]^-\varphi \ar[d]_{\can} & \prismcpl_S\otimes I^{d_1}\phi(I)^{d_2}\\
  \Nyg^{\ge d_2}\prismcpl_S
      \ar[r]^-\varphi \ar[d] & \prismcpl_S\otimes I^{d_2}\\
  0 &
}
\]
\end{example}

\begin{theorem}
Let $S$ be a quasiregular semiperfectoid ring, and let $\alpha=(d_0,\ldots,d_n;d_\infty)$. Then there is a canonical identification
\[
  \TC^{-n}_{\alpha}(S;\Z_p)
  =\Nyg^{\ge d_0,\ldots,d_n}\prismcpl_S
  \fade{{}\tnsr I^{\alpha^{(n+1)}}}.
\]
The group $\TC^{-n}_{\alpha-1}(S;\Z_p)$ has a filtration with associated graded
\[
  \gr^i\TC^{-n}_{\alpha-1}(S;\Z_p)=
  \frac
    {\prismcpl_S}
    {
      \Nyg^{\ge d_{n-i}}\prismcpl_S
      +
      \~\varphi^i(\Nyg^{\ge d_{(n+1-i)}, \dotsc, d_n}\prismcpl_S)
    }
  \tnsr
  I^{\alpha^{(n+1-i)}}
\]
for $1\le i\le n$ (indexed such that $\gr^0=0$). It is given explicitly by the iterated pushout
\[
\xymatrix{
    &&& \Nyg^{\ge d_0}\prismcpl_S
    \tnsr I^{\alpha^{(1)}}
      \ar[r] \ar[d]_-\can & 0\\
                          && \Nyg^{\ge d_{(n-k-1)}}\prismcpl_S
                    \tnsr I^{\alpha^{(n-k)}}
      \ar[r]^-\varphi \ar[d]_-\can & \cdots\\
                                   & \Nyg^{\ge d_{n-k}} \prismcpl_S \tnsr I^{\alpha^{(n-k+1)}} \ar[r]^-\varphi \ar[d]_-\can & \prismcpl_S \tnsr I^{\alpha^{(n-k)}}\\
  \Nyg^{\ge d_n}\hat{\prism}_S
  \tnsr I^{\alpha^{(n+1)}}
      \ar[r]^-\varphi \ar[d] & \dots\\
  0 &
}
\]
\end{theorem}

\subsubsection{Reindexing and $\TR^{h\T}$}

To tackle $\TR^{h\T}$, it will be convenient to reindex our notation. $\TR^{h\T}$ is naturally graded on $\lim_{\alpha\mapsto\alpha'}\RO(\T)$, which we identify with bi-infinite sequences $\m{\alpha}=(\ldots,d_{-2},d_{-1},d_0,d_1,d_2,\ldots)$ such that $d_i=d_{i+1}$ for $i\gg0$. For typographical reasons, we write $\m{d}_i$ instead of $d_{-i}$. Next, we shift our indexing so that
\[
  \TCmin_{\m\alpha}(S;\Z_p)
    =\Nyg^{\ge\m d_1}\prismcpl_S\otimes I^\alpha,
  \qquad
  \TP_{\m\alpha}(S;\Z_p)
    =\prismcpl_S\otimes I^\alpha.
\]

Furthermore, we write $\prismcpl_S^{(-n)}$ for $\prismcpl_S$ considered as a $\prismcpl_S$-algebra via the map $\varphi^n$. We then write $\prismcpl_S^{(-n)}\otimes\phi^{i-n}(I)$ instead of $\prismcpl_S\otimes \phi^i(I)$, and we set
\[ I^{\m\alpha(i)} = \bigotimes_{j=-i}^\infty \phi^j(I)^{d_j}. \]
To do a computation in $\TC^{-n}_\rog$, $d_i$ in the old notation corresponds to $\m d_{(n+1-i)}$ in the new notation. Thus our previous computation of $\TC^{-1}_\rog(S;\Z_p)$ becomes
\[
\begin{aligned}
  \TC^{-1}_{\m{\alpha}}(S;\Z_p)
    &=\Nyg^{\ge\m{d}_2,\m{d}_1}\prismcpl_S^{(0)}
    \tnsr I^{\m\alpha(0)}\\
  \TC^{-1}_{\m{\alpha}-1}(S;\Z_p)
    &=\frac{\prismcpl_S^{(-1)}}
      {\Nyg^{\ge\m{d}_2}\prismcpl_S^{(-1)} +\~\varphi\left(\Nyg^{\ge\m{d}_1}\prismcpl_S^{(0)}\right)}
       \tnsr I^{\m\alpha(1)}.
\end{aligned}
\]

\begin{theorem}
Let $S$ be a quasiregular semiperfectoid ring. Then there is a canonical identification
\[
  \TC^{-n}_{\m{\alpha}}(S;\Z_p)
  =\Nyg^{\ge\m{d}_{n+1},\ldots,\m{d}_1}
  \prismcpl_S^{(0)} \tnsr I^{\m\alpha(0)}.
\]
The group $\TC^{-n}_{\m{\alpha}-1}(S;\Z_p)$ has a filtration with associated graded
\[
  \gr^i\TC^{-n}_{\m{\alpha}-1}(S;\Z_p)=
  \frac
    {\prismcpl_S^{(-i)}}
    {
      \Nyg^{\ge\m{d}_{(i+1)}}\prismcpl_S^{(-i)}
      +
      \~\varphi^i\left(\Nyg^{\ge\m{d}_i,\ldots,\m{d}_1}\prismcpl_S^{(0)}\right)
    }
  \tnsr
  I^{\m\alpha(i)}
\]
for $1\le i\le n$ (indexed such that $\gr^0=0$). It is given explicitly by the iterated pushout
\[
\xymatrix@C=4em@R=2em{
  && \Nyg^{\ge\m{d}_{n+1}}\prismcpl_S^{(-n)}
      \tnsr I^{\m\alpha(n)}
      \ar[r] \ar[d] & 0\\
  & \Nyg^{\ge\m{d}_2}\prismcpl_S^{(-1)}
      \tnsr I^{\m\alpha(1)}
      \ar[r] \ar[d] & \cdots\\
  \Nyg^{\ge\m{d}_1}\prismcpl_S^{(0)}
  \tnsr I^{\m\alpha(0)}
      \ar[r] \ar[d] & \prismcpl_S^{(-1)}
      \tnsr I^{\m\alpha(1)}\\
  0 &
}
\]
\end{theorem}

In general we do not know how to take the limit to get $\TR^{h\T}$; but see \cite[Lemma 3.23]{BMS1} for controlling $\lim\nolimits^1$ terms.

\begin{remark}
Let $\Sol_p=\lim(\cdots\xrightarrow{p}\T\xrightarrow{p}\T)$ be the $p$-adic solenoid, and let $K_i$ be the kernel of the map to the $i^\th$ copy of $\T$. If $U$ is a (possibly infinite-dimensional) representation of $\Sol_p$, say that $U$ is \emph{quasifinite-dimensional} if each $U^{K_i}$ is finite-dimensional (terminology by analogy with \cite{KMNPolygonic}). Then the grading group $\lim_{\alpha\mapsto\alpha'}\RO(\T)$ can be identified with the Grothendieck group of quasifinite-dimensional representations of $\Sol_p$.
\end{remark}

%% file: dRW.tex
%
\section{\tops{$\RO(\Z_p)^\land$}{RO(Z\_p)\^}-graded \tops{$\TR$}{TR} of smooth algebras}
\label{sec:dRW}

\newcommand{\untwist}{\smash{{}^{\fade{\scriptscriptstyle\phi^{-1}}}}}
\newcommand{\thh}{\THH(S;\Z_p)}

In this section we study the polyrepresentation-graded $\TR$ of a $p$-completely smooth algebra $S$ over a perfectoid base ring $R$, deducing Norm maps for the de Rham-Witt complex as our main application. We recall the explicit description of the de Rham-Witt complex in this setting in \sec\ref{sub:dRW-perfd}. In \sec\ref{sub:ctgs}, we extend the homotopy orbits to TR spectral sequence to compute $\TR^\bullet_\rog(S;\Z_p)$ in the contiguous range. We then record some special cases in \sec\ref{sub:special}. We identify the Norm on the de Rham-Witt complex (borrowed from the Norm on $\TR$) in \sec\ref{sub:dRW-norm}.

Let us explain the strategy. As a $\T$-spectrum, the representation sphere $S^{\lambda^k}$ has a 1-step cell structure with associated graded
\[ S^0,\quad\T/C_{k+}\tnsr S^1. \]
As a $C_n$-spectrum, it instead has a 2-step cell structure with associated graded
\[ S^0,\quad C_n/C_{k+}\tnsr S^1,\quad C_n/C_{k+}\tnsr S^2. \]

The ``short'' cell structure was very effective at computing $\TCmin_\rog$ and $\TP_\rog$ in \sec\ref{sec:qrsp} since they are even. If one uses the ``long'' cell structure to compute $\TR^\bullet_\rog$, one finds that the 1-skeleton contributes $V^k$ and the 2-skeleton contributes $dV^k$ to $\Fil^k\dRW[\bullet]^*=V^k + dV^k$. While instructive, this computation is extremely messy, so we defer it to a companion paper. Instead, we will use the BMS filtration to import the results of \sec\ref{sec:qrsp} and compute $\TR^\bullet_\rog$ via the homotopy orbits to $\TR$ spectral sequence from \cite[\sec3]{ROS1TR} and \cite[\sec3]{SulROGTHH}.

\subsection{Structure of the de Rham-Witt complex over a perfectoid base}
\label{sub:dRW-perfd}

Let $R$ be a perfectoid ring, and let $S$ be a $p$-completely smooth $R$-algebra. We are interested in the relative de Rham-Witt complex $\dRW[\bullet]^*_{S/R}$ of Langer-Zink \cite{LangerZinkDRW}. We recall from \cite{BMS1} the structure of this complex and its connection to the $q$-de Rham complex.

\begin{example}[Cartier isomorphism]
  \label{ex:cartier}
  Let $R$ be a perfectoid ring admitting all $p$-power roots of unity. Fix a compatible choice of roots $\{\zeta_{p^n}\}$, and define the following elements of $A=\Ainf(R)=\rmW(R^\flat)$:
  \begin{align*}
    q &= [\epsilon] & \~\xi_n &= (p^n)_q & \xi_n &= (p^n)_{q^{1/p^{n+1}}}\\
    \mu &= 1-q &&= \frac{\phi^n(\mu)}\mu &&= \frac{\phi^{-1}(\mu)}{\phi^{-(n-1)}(\mu)}
  \end{align*}
  In this context, prismatic cohomology of $R$-algebras can be computed as a $q$-de Rham complex by \cite[Lemma 9.6]{BMS1} and \cite[Theorem 17.2]{Prismatic}.

  We work this out for the affine line $S=R[x]^\land_p$. The $q$-de Rham complex is given by
  \[ \AOmega_S = \widehat{\bigoplus_n} \left[A{\<x^n\>} \xra{\nabla_q} A{\<x^n\dlog x\>}\right] \]
  where $\widehat\bigoplus$ denotes a $p$-completed direct sum, and the differential sends $\nabla_q(x^n) = (n)_q x^n \dlog x$.

  Let us compute the cohomology of $\~\Omega_S = \AOmega_S/\~\xi = \AOmega_S/(p)_q$. Now the differential is $0$ for exponents divisible by $p$, and an isomorphism otherwise. We conclude that
  \begin{align*}
    \H^*(\~\Omega_S) &= \widehat{\bigoplus_n} A/\~\xi{\<x^{pn}\>}[0] \oplus A/\xi{\<x^{pn}\dlog x\>}[-1]
  \end{align*}
  These are isomorphic to
  \[\xymatrix@R=2em@C=-4em{
      *!R{A/\xi{\<x^n\>}\ } \ar[r]^-\phi_-\sim & *!L{A/\~\xi{\<x^{pn}\>}}\\
      *!R{A/\xi{\<x^n\dlog x\>}\tnsr \xi^{-1}\ } \ar[r]^-\phi_-\sim & *!L{A/\~\xi{\<x^{pn}\dlog x\>}}
  }\]
  The twist by $\xi^{-1}$ is necessary in the second line because $\phi(\dlog x)=\~\xi\dlog x$. In other words,
  \[ \H^*(\~\Omega_S) = \Omega^*_{S/R}\{-*\}=\Omega^*_{S/R}\tnsr \xi^{-*}A, \]
  in accordance with \cite[Theorem 8.7]{BMS1}.

Next we compute the cohomology of $\widetilde{\dRW[1]}_S = \AOmega_S/\~\xi_2 = \AOmega_S/(p^2)_q$. We find that
  \begin{align*}
    \H^*(\widetilde{\dRW[1]})
      &= \phantom{{}\oplus} \bigoplusp_n A/\~\xi_2\<x^{p^2 n}\>[0] \oplus A/\~\xi_2\<x^{p^2 n}\dlog x\>[-1]\\
      &\phantom{{}=} \oplus \bigoplusp_{p\nmid j} A/\~\xi\<\phi(\~\xi) x^{pj}\>[0] \oplus A/\~\xi\<x^{pj}\dlog x\>[-1].
  \end{align*}
  These are isomorphic to
  \[\xymatrix@R=2em@C=-5em{
      *!R{A/\xi_2\<x^n\>\ } \ar[r]^-{\phi^2}_-\sim & *!L{A/\~\xi_2\<x^{p^2n}\>}\\
      *!R{A/\phi^{-1}(\xi)\<\xi x^{j/p}\>\ } \ar[r]^-{\phi^2}_-\sim & *!L{A/\~\xi\<\phi(\~\xi) x^{pj}\>}\\
      *!R{A/\xi_2\<x^n\dlog x\>\tnsr \xi_2^{-1}\ } \ar[r]^-{\phi^2}_-\sim & *!L{A/\~\xi_2\<x^{p^2n}\dlog x\>}\\
      *!R{A/\phi^{-1}(\xi)\<x^{j/p}\dlog x\>\tnsr \xi_2^{-1}\ } \ar[r]^-{\phi^2}_-\sim & *!L{A/\~\xi\<x^{pj}\dlog x\>}
  }\]
  In other words, $\H^*(\widetilde{\dRW[1]}_S) = \dRW[1]^*_{S/R}\{-*\}=\dRW[1]^*_{S/R} \tnsr \xi_2^{-*}A$, in accordance with \cite[Theorem 11.1]{BMS1}.
\end{example}

We need a similar calculation which is not covered by \cite[Theorem 11.1]{BMS1}.

\begin{proposition}
  \label{prop:h*-nyg}
  If $S$ is a $p$-completely smooth $R$-algebra, then
  \[
    \H^j\left(\frac{\Nyg^{\ge i}\AOmega_S}{I_n \Nyg^{\ge i+1}\AOmega_S}\right)
    =
    \left\{
      \begin{array}{r@{}ll}
        \dRW[n]^j_{S/R} &{}\tnsr \xi_{n+1}^{-j} \phi^{-n}(\xi)^{i}A & i\ge j\\
        \dRW[n-1]^j_{S/R} &{}\tnsr \xi_n^{-j}A & i < j
      \end{array}
    \right.
  \]
\end{proposition}
\begin{proof}
  By \cite[Remark 9.11]{BMS2}, the Nygaard filtration is given explicitly by
  \[ \Nyg^{\ge i}\AOmega^j_S = \xi^{\pos{i-j}}\AOmega^j_S \]
  where $\pos x \defeq \max(x, 0)$. If $i<j$, then this is just $\H^j(\AOmega_S/I_n)=\dRW[n-1]^j_{S/R} \tnsr \xi_n^{-j}A$ by the existing Cartier isomorphism \cite[Theorem 11.1]{BMS1}. For $j \le i$, the Frobenius induces an equivalence
  \[ \tau^{\le j} \Nyg^{\ge i} \AOmega_S \morph^\varphi_\sim \tau^{\le j} I^i A\Omega_S \]
  by (applying $\tau^{\le j}$ to both sides of) \cite[Theorem 2.16]{BouisCartier}. Thus we get
  \begin{align*}
    \H^j\left(\frac{\Nyg^{\ge i}\AOmega_S}{I_n\Nyg^{\ge i+1}\AOmega_S}\right)
    &= \H^j\left(\frac{\tau^{\le j}\Nyg^{\ge i}\AOmega_S}{\tau^{\le j+1}I_n\Nyg^{\ge i+1}\AOmega_S}\right)\\
    &\morph^\varphi_\sim \H^j\left(\frac{\tau^{\le j}I^i\AOmega_S}{\tau^{\le j+1}I^i I_{n+1}\AOmega_S}\right)\\
    &= \H^j\left(\frac{\AOmega_S}{I_{n+1}\AOmega_S}\right) \tnsr I^i\\
    &\xra[\sim]{\phi^{-{n+1}}} \dRW[n]^j_{S/R} \tnsr \xi_{n+1}^{-j}A \tnsr \phi^{-n}(\xi)^i A. \qedhere
  \end{align*}
\end{proof}

\begin{proofsketch}
  We sketch a modern proof of Hesselholt's theorem, as formulated in Notation \ref{not:canonical}:
  \[ \TR^n_{2i,j}(S;\Z_p) = \dRW[n]^j_{S/R} \tnsr \xi_{n+1}^i A. \]

  The base case $n=0$ is \cite[Corollary 6.9]{BMS2}. Inductively, we consider the isotropy separation square
  \[\xymatrix{
      \TR^n(S;\Z_p) \ar[r]^-R \ar[d] \pullback & \TR^{n-1}(S;\Z_p) \ar[d]\\
      \thh^{hC_{p^n}} \ar[r] & \thh^{tC_{p^n}}
  }\]
  By \cite[Lemma 3.1]{RiggenbachTruncated}, the motivic filtrations are
  \begin{align*}
    \gr^i_\mot \thh^{hC_{p^n}}
      &= \frac{\Nyg^{\ge i}\hat\prism_S}{\~\xi_n \tnsr \Nyg^{\ge i+1}\hat\prism_S}\tnsr \~\xi_n^i A[2i]\\
      &= \frac{\Nyg^{\ge i}A\Omega_S}{\~\xi_n \tnsr \Nyg^{\ge i+1}A\Omega_S}\tnsr \~\xi_n^i A[2i]\\
    \gr^i_\mot \thh^{tC_{p^n}}
      &= \frac{\hat\prism_S}{\~\xi_n}\tnsr \~\xi_n^i A[2i]\\
      &= \frac{A\Omega_S}{\~\xi_n}\tnsr \~\xi_n^i A[2i]\\
      &= \~{\dRW[n-1]_S}\tnsr \~\xi_n^i A[2i]\\
  \end{align*}
  Taking cohomology and applying the Cartier isomorphism (\cite[Theorem 8.7]{BMS1} and Proposition \ref{prop:h*-nyg}), we find that $\thh^{hC_{p^n}}_*$ is computed by the following $E_1$ spectral sequence:
  \[\xymatrix@C=1em@R=1em{
      \fade\dots & \fade{\dRW[n-1]^5 \tnsr \xi_n^{-3}} & \fade{\dRW[n-1]^4 \tnsr \xi_n^{-2}} & \fade{\dRW[n-1]^3 \tnsr \xi_n^{-1}} & \dRW[n]^2 & \dRW[n]^1 \tnsr \xi_{n+1} \\
      \fade\dots & \fade{\dRW[n-1]^3 \tnsr \xi_n^{-2}} & \fade{\dRW[n-1]^2 \tnsr \xi_n^{-1}} & \dRW[n]^1 & \dRW[n]^0 \tnsr \xi_{n+1}\\
      \fade{\dRW[n-1]^2 \tnsr \xi_n^{-2}} &\fade{\dRW[n-1]^1 \tnsr \xi_n^{-1}} &
      \dRW[n]^0 \ar@[lightgray][ul]\\
    \fade{\dRW[n-1]^0 \tnsr \xi_n^{-1}}\\
      \dots & \pi_{-1} & \pi_0 & \pi_1 & \pi_2 & \dots
  }\]
  Observe how the negative twists from the Cartier isomorphism cancel out the Breuil-Kisin twists from the motivic filtration, shifting de Rham-Witt terms from cohomological to homological degree. Hesselholt's proof \cite{LarsPTyp} based on the cyclic bar construction establishes an additive decomposition which shows that there are no further differentials or extensions (although it should be possible to give a direct proof within the current formalism).

  The solid region is where $\TR^{n-1}_* \to \THH^{tC_{p^n}}_*$ is an isomorphism, which implies the same for $\TR^n_* \to \THH^{hC_{p^n}}_*$; and the faded region is an isomorphism for $\THH^{hC_{p^n}} \to \THH^{tC_{p^n}}$, which implies the same for $\TR^n_* \to \TR^{n-1}_*$. Hesselholt's theorem follows. \qed
\end{proofsketch}

\subsection{Contiguous range}
\label{sub:ctgs}

In this section we compute $\TR_\rog(S; \Z_p)$ in the contiguous range. Our method is the homotopy orbits to TR spectral sequence from \cite{ROS1TR,SulROGTHH}. In fact, the computation for $\TR_{\rog,j}(S;\Z_p)$ works the same way as the computation for $\TR_\rog(R;\Z_p)$.

\begin{method}
  We explain the strategy for computing $\TR^1_\rog(S;\Z_p)$. We will use the isotropy separation square
  \begin{equation}
    \label{eq:isotropy}
    \vcenter{\xymatrix{
      \TR^1(S;\Z_p) \ar[r] \ar[d] \pullback & \thh^{\Phi C_p} \ar[d]\\
      \thh^{hC_p} \ar[r] & \thh^{tC_p}
    }}
  \end{equation}
  Recall that with this method, the computation breaks up as a direct sum over de Rham-Witt degrees $j$, where $\pi_{\alpha,j}$ denotes the $\dRW^j$ contribution to $\pi_{\alpha+j}$. With this notation, the spectral sequence associated to the isotropy separation square takes the form
  \begin{align*}
    E_0 &=\dbox{$\vcenter{\xymatrix@R=0.2em{
          \thh^{tC_p}_{\alpha,j} & \thh^{hC_p}_{\alpha,j} \ar[l]\\
                              & \thh^{\Phi C_p}_{\alpha,j}
                   }}$}\\
    E_1 &=
    \dbox{$\vcenter{\xymatrix@R=0.2em{
    (\thh_{hC_p})_{\alpha-1,j} & (\thh_{hC_p})_{\alpha,j} \\
                            & \thh^{\Phi C_p}_{\alpha,j} \ar[ul]^-\varphi
                                                               }}$}\\
    E_\infty &=
    \dbox{$\vcenter{\xymatrix@R=0.2em{
    \TR^1_{\alpha-1,j}(S;\Z_p) & \gr^1 \TR^1_{\alpha,j}(S;\Z_p)\\
                               & \gr^0\TR^1_{\alpha,j}(S;\Z_p)
                                    }}$}
  \end{align*}

  We need to identify these terms canonically. Suppose that $\alpha=(d_0;d_\infty)$ is a complex $C_p$-representation. Since the Euler class $a_\lambda$ is invertible in the right column of \eqref{eq:isotropy}, the corresponding groups canonically depend only on $d_\infty$: we have canonical identifications
  \begin{align*}
    \thh^{\Phi C_p}_{\alpha,j} &= \thh^{\Phi C_p}_{2d_\infty,j}\\
                               &= \begin{cases}
                                 \Omega^j_{S/R} \tnsr \xi^{d_\infty} & d_\infty \ge 0\\
                                 0 & d_\infty < 0,
                               \end{cases}\\
    \thh^{tC_p}_{\alpha,j} &= \thh^{t C_p}_{2d_\infty,j} \\
                           &= \widetilde{\Omega}^j_{S/R} \tnsr \xi^{d_\infty}.
  \end{align*}
  
  Since the Thom class $u_\lambda$ is invertible in the bottom row of \eqref{eq:isotropy} (and in particular on $\thh_{hC_p}$), the corresponding groups \emph{non-canonically} depend only on $d_0$: we have \emph{non-canonical} identifications
  \begin{align*}
    \thh^{hC_p}_{\alpha,j} &\mathrel{\red\isom} \thh^{hC_p}_{2d_0,j}\\
                           &= \left\{\!\!\!
                             \begin{array}{r@{}ll}
                               \dRW[1]^j_{S/R} &{}\tnsr (\xi\phi^{-1}(\xi))^{d_0} & d_0\ge0\\
                               \Omega^j_{S/R} &{}\tnsr \xi^{d_0} & d_0<0
                           \end{array}\right.
  \end{align*}
  To make this canonical, we need to match the twist on $\thh^{tC_p}$. The map
  \[ \thh^{hC_p}_{2d_0,j} \morph \thh^{tC_p}_{2d_0,j} \]
  is $\phi^{-1}(\xi)^{d_0}$ times the canonical projection when $d_0\ge 0$, and an isomorphism when $d_0<0$. Thus, the canonical identification is
  \[
    \THH^{hC_p}_{\alpha,j} =
    \left\{\!\!\!
      \begin{array}{r@{}ll}
        \dRW[1]^j_{S/R} &{}\tnsr \xi^{d_0} \phi^{-1}(\xi)^{d_\infty} & d_0\ge 0\\
        \Omega^j_{S/R} &{}\tnsr \xi^{d_\infty} & d_0 < 0
    \end{array}
    \right.
  \]

  Finally, we switch to poly-representation notation: if $\alpha = (d_0;d_\infty)$ is a $\T$-representation, the corresponding $\Z_p$-polyrepresentation is $\m\alpha=(\m d_0, \m d_1)$ with $\m d_0=d_\infty$ and $\m d_1 = d_0$. Thus our canonical identifications are
  \[
    \dbox{$\vcenter{\xymatrix@R=0.2em{
          \thh^{tC_p}_{\alpha,j} & \thh^{hC_p}_{\alpha,j} \\
                              & \thh^{\Phi C_p}_{\alpha,j}
                   }}$}
    =
    \dbox{$\vcenter{\xymatrix@R=0.2em@C=-4em{
          \Omega^j_{S/R} \tnsr \xi^{\m d_0} & *!L{\begin{cases}\dRW[1]^j_{S/R} \tnsr \xi^{\m d_0} \phi^{-1}(\xi)^{\m d_1} & \m d_1 \ge 0\\ 
              \Omega^j_{S/R} \tnsr \xi^{\m d_0} & \m d_1 < 0
                                      \end{cases}} \\
                                      & *!L{\begin{cases}\Omega^j_{S/R} \tnsr \xi^{\m d_0} & \m d_0 \ge 0\\ 
                                          0 & \m d_0 < 0
                                      \end{cases}}
                                                               }}$}
      \]
\end{method}

\begin{example}
  First we compute $\TR^1_\rog(S\rel\SAinf{R^+})$ when $R^+$ is a $p$-torsionfree perfectoid ring and $S$ is a $p$-completely smooth $R^+$-algebra. Let $\m\alpha = (\m d_0,\m d_1)$ be a virtual $\Z/p$-polyrepresentation. There are three cases, depending on the signs of $\m d_0$ and $\m d_1$ (the result is zero if both are negative).

  \begin{center}
    \begin{tabular}{c|c}
      (1) & (2)\\\hline
      0 & (3)
    \end{tabular}
  \end{center}

  \begin{enumerate}
    \item $\m d_0\ge0$, $\m d_1<0$:
      \[
        \left(\vcenter{\xymatrix@R=0em{
              \Omega^j \tnsr \xi^{\m d_0} & \Omega^j \tnsr \xi^{\m d_0} \ar[l]_-\sim\\
                                          & \Omega^j \tnsr \xi^{\m d_0}
        }}\right)
        \leadsto
        \left(\vcenter{\xymatrix@R=0em{
              \phantom0 & \phantom0\\
               & \Omega^j \tnsr \xi^{\m d_0}
        }}\right)
      \]
      In this region, $\TR^1_\rog$ coincides with $\THH^{\Phi C_p}_\rog$.

    \item $\m d_0\ge 0$, $\m d_1\ge 0$: since $\rmW_0(R^+)=A/\xi$ is $\phi^{-1}(\xi)$-torsionfree,
      \begin{align*}
        &\phantom{{}\leadsto}\left(\vcenter{\xymatrix@R=0em{
        \Omega \tnsr \xi^{\m d_0} & {\dRW[1]^j \tnsr \xi^{\m d_0} \phi^{-1}(\xi)^{\m d_1}} \ar[l]\\
                                        & {\Omega^j \tnsr \xi^{\m d_0}}
        }}\right)\\
        &\leadsto
        \left(\vcenter{\xymatrix@R=0em{
          \dfrac{\Omega^j}{\phi^{-1}(\xi)^{\m d_1}} \tnsr \xi^{\m d_0} & \Fil^1\dRW[1]^j \tnsr \xi^{\m d_0} \phi^{-1}(\xi)^{\m d_1}\\
                                                                  & \Omega^j \tnsr \xi^{\m d_0} \ar[ul]
        }}\right)\\
        &\leadsto
        \left(\vcenter{\xymatrix@R=0em@C=10em{
        \phantom0 & *!R{\Fil^1\dRW[1]^j \tnsr \xi^{\m d_0} \phi^{-1}(\xi)^{\m d_1}}\\
              \phantom0 & *!R{\Omega^j \tnsr \xi^{\m d_0} \phi^{-1}(\xi)^{\m d_1}}
        }}\right)
      \end{align*}
      The extension is $\dRW[1]^j_{S/R^+} \tnsr \xi^{\m d_0} \phi^{-1}(\xi)^{\m d_1}$. In this region, $\TR^1_\rog$ coincides with $\THH^{hC_p}_\rog$.

    \item $\m d_0<0$, $\m d_1\ge0$: again using that $R^+$ is $\phi^{-1}(\xi)$-torsionfree,
      \[
        \left(\vcenter{\xymatrix@R=0em{
          \Omega \tnsr \xi^{\m d_0} & {\dRW[1]^j \tnsr \xi^{\m d_0} \phi^{-1}(\xi)^{\m d_1}} \ar[l]\\
        }}\right)
        \leadsto
        \left(\vcenter{\xymatrix@R=0em{
              \dfrac{\Omega^j}{\phi^{-1}(\xi)^{\m d_1}} \tnsr \xi^{\m d_0} & \Fil^1\dRW[1]^j \tnsr \xi^{\m d_0} \phi^{-1}(\xi)^{\m d_1}
        }}\right)
      \]
      In this region, $\TR^1_\rog$ coincides with $(\THH_{hC_p})_\rog$.
  \end{enumerate}

  In summary:
  \begin{center}
    \renewcommand{\arraystretch}{2} 
    \includegraphics[width=0.85\textwidth]{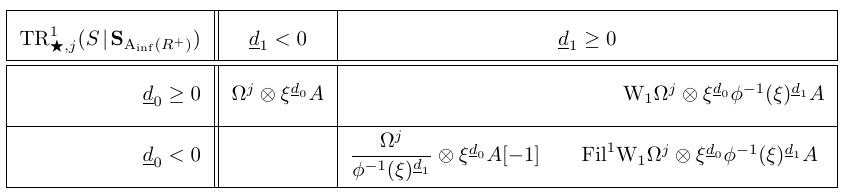}
  \end{center}
\end{example}

For an arbitrary perfectoid ring $R$, the argument in the first region goes through. In the second region, the intermediate pages may look different, but the final answer is the same: regardless of $R$, $\THH^{\Phi C_p}_\rog \to \THH^{tC_p}_\rog$ is an isomorphism in this region, so a priori so must be $\TR^1_\rog \to \THH^{hC_p}_\rog$. The third region (as well as the $E^1$ page of the second region) depends more delicately on the $p$-torsion of $R$.

\begin{example}
  Next we compute $\TR^1_{\rog,j}(S\rel\SW k)$ when $k$ is a perfect $\F_p$-algebra and $S$ is a smooth $k$-algebra. By the preceding remark, it suffices to consider the second region. Note that
  \begin{align*}
    \ker\left(\dRW[1]^j \tnsr p^{\m d_0 + \m d_1} \to \Omega^j \tnsr p^{\m d_0} \right)
    &=
    \begin{cases}
      \Fil^1 \dRW[1]^j \tnsr p^{\m d_0 + \m d_1} & \m d_1 = 0\\
      \phantom{\Fil^1} \dRW[1]^j \tnsr p^{\m d_0 + \m d_1} & \m d_1 > 0
    \end{cases}\\
    &= \Fil^{1 - \m d_1} \dRW[1]^j \tnsr p^{\m d_0 + \m d_1}
  \end{align*}
  recalling that by convention $\Fil^r \dRW[\bullet]^* = \dRW[\bullet]^*$ for $r\le0$. This may look strange, but it will prevent us from losing our minds with casework. We similarly have the peculiar identity
  \[ \ker\left( \Omega^j \to \frac{\Omega^j}{p^{\m d_1}} \right) = \dRW[-\m d_1]^j \tnsr p^{\m d_1}W \]
  where the right-hand side is zero for $\m d_1>0$.

  Now the spectral sequence goes
  \begin{align*}
    &\phantom{{}\leadsto}\left(\vcenter{\xymatrix@R=0em{
    \Omega \tnsr p^{\m d_0} & {\dRW[1]^j \tnsr p^{\m d_0 + \m d_1}} \ar[l]\\
                                    & {\Omega^j \tnsr p^{\m d_0}}
    }}\right)\\
    &\leadsto
    \left(\vcenter{\xymatrix@R=0em{
  \dfrac{\Omega^j}{p^{\m d_1}} \tnsr p^{\m d_0} & \Fil^{1-\m d_1}\dRW[1]^j \tnsr p^{\m d_0 + \m d_1}\\
                                                              & \Omega^j \tnsr p^{\m d_0} \ar[ul]
    }}\right)\\
    &\leadsto
    \left(\vcenter{\xymatrix@R=0em@C=10em{
 \phantom0 & *!R{\Fil^{1 - \m d_1}\dRW[1]^j \tnsr p^{\m d_0 + \m d_1}}\\
 \phantom0 & *!R{\dRW[-\m d_1]^j \tnsr p^{\m d_0} }
    }}\right)
  \end{align*}
  and the extension is $\dRW[1]^j_S \tnsr p^{\m d_0 + \m d_1}W$.

  In summary:
  \begin{center}
    \renewcommand{\arraystretch}{2} 
    \includegraphics[width=0.76\textwidth]{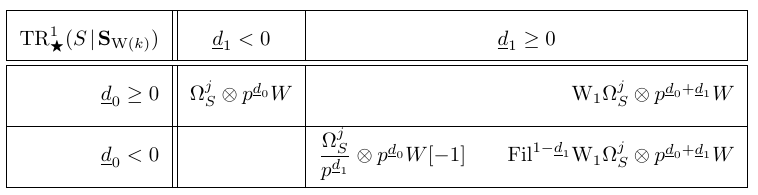}
\end{center}
\end{example}

\begin{example}
  We can also compute these as Mackey functors. The basic Mackey functors in play are
  \[\xymatrix{
      0 \ar[r] & \Fil^1\dRW[1]^j \mdown F \ar[r] & \dRW[1]^j \mdown F \ar[r] & \Omega^j \mdown{} \ar[r] & 0 \\
      0 \ar[r] & \Omega^j \ar@{=}[r] \mup V & \Omega^j \mup V \ar[r] & 0 \mup{} \ar[r] & 0
  }\]
  (we do not need the fact that these fit into a short exact sequence). The actions on the twists (remembering that $F$ is a \emph{divided} Frobenius map) are given by
  \[
    F(1\tnsr f\xi_2) = 1\tnsr\phi(f)\xi, \quad\qquad  F(1\tnsr f\xi) = 1 \tnsr \phi(f), \quad\qquad  F(1\tnsr f\phi^{-1}(\xi)) = 1\tnsr \phi(f)\xi.
  \]
\end{example}

\begin{definition}
  \label{def:poly-ideal}
  Given a poly-representation $\m\alpha=(\m d_0,\m d_1,\dotsc)$, let $I^{\m\alpha(r)}$ denote the fractional ideal
  \[ I^{\m\alpha(r)} \defeq \bigotimes_{i=0}^r \phi^{-i}(\xi)^{\m d_i}A. \]
  When $\m d_0,\dotsc,\m d_r\ge0$, this agrees with the iterated Nygaard filtration $\Nyg^{\ge \m d_r,\dotsc,\m d_0}A$ from \sec\ref{sub:iterated-nygaard}.
\end{definition}

\begin{example}
  We now study $\TR^2_{\rog,j}(S\rel \SAinf{R^+})$ in the range $\m d_0,\,\m d_1,\,\m d_2\ge0$. By Tsalidis' theorem \cite[Theorem 2.4]{Tsalidis}, we know a priori that the result will be isomorphic to $\dRW[2]^j_{S/R^+}$. 
  The spectral sequence goes
  \begin{align*}
    E_0 &=\dbox{$\vcenter{\xymatrix@R=0em{
        \dRW[1]^j \tnsr I^{\m\alpha(1)} & \dRW[2]^j \tnsr I^{\m\alpha(2)} \ar[l]\\
        \dRW[0]^j \tnsr I^{\m\alpha(0)} & \dRW[1]^j \tnsr I^{\m\alpha(1)} \ar[l]\\
                         & \dRW[0]^j \tnsr I^{\m\alpha(0)}
                   }}$}\\
    E_1 &=
    \dbox{$\vcenter{\xymatrix@R=0.2em@C=-5em{
    *!R{\dRW[1]^j/\phi^{-2}(\xi)^{\m d_2} \tnsr I^{\m\alpha(1)}} & *!L{\Fil^2\dRW[2]^j \tnsr I^{\m\alpha(2)}}\\
    *!R{\dRW[0]^j/\phi^{-1}(\xi)^{\m d_1} \tnsr I^{\m\alpha(0)}} & *!L{\Fil^1\dRW[1]^j \tnsr I^{\m\alpha(1)}} \ar[ul]\\
                                                                     & *!L{\fade{\Fil^0{}}\dRW[0]^j \tnsr I^{\m\alpha(0)}} \ar[ul]
                                                               }}$}\\
    E_2 &=
    \dbox{$\vcenter{\xymatrix@R=0.2em@C=-5em{
  *!R{\dRW[0]^j/\phi^{-2}(\xi)^{\m d_2} \tnsr I^{\m\alpha(1)}} & *!L{\Fil^2\dRW[2]^j \tnsr I^{\m\alpha(2)}}\\
                                          & *!L{\Fil^1\dRW[1]^j \tnsr I^{\m\alpha(2)}}\\
                                          & *!L{\fade{\Fil^0{}}\dRW[0]^j \tnsr I^{\m\alpha(1)}} \ar[uul]
                                    }}$}\\
    E_\infty&=
    \dbox{$\vcenter{\xymatrix@R=0.2em@C=8.2em{
            & *!L{\Fil^2\dRW[2]^j \tnsr I^{\m\alpha(2)}}\\
            & *!L{\Fil^1\dRW[1]^j \tnsr I^{\m\alpha(2)}}\\
            & *!L{\fade{\Fil^0{}}\dRW[0]^j \tnsr I^{\m\alpha(2)}}
      }}$}
  \end{align*}
  The extension is $\dRW[2]^j_{S/R^+} \tnsr I^{\m\alpha(2)}A$.
\end{example}

\begin{example}
  We repeat the previous calculation in the characteristic $p$ case. We use the identification
  \[ \Fil^r\dRW[n]^* \cap p^i\dRW[n]^* = \Fil^{r-i} \dRW[n-i]^*\tnsr p^i\rmW. \]
  The spectral sequence goes
  \begin{align*}
    E_0 &=\dbox{$\vcenter{\xymatrix@R=0em@C=4em{
          \dRW[1]^j \tnsr I^{\m\alpha(1)} & \dRW[2]^j \tnsr I^{\m\alpha(2)} \ar[l]\\
          \dRW[0]^j \tnsr I^{\m\alpha(0)} & \dRW[1]^j \tnsr I^{\m\alpha(1)} \ar[l]\\
                                               & \dRW[0]^j \tnsr I^{\m\alpha(0)}
                   }}$}\\
    E_1 &=
    \dbox{$\vcenter{\xymatrix@R=0.2em@C=-4.1em{
    *!R{\dRW[1]^j/p^{\m d_2} \tnsr I^{\m\alpha(1)}} & *!L{\Fil^{2-\m d_2}\dRW[2]^j \tnsr I^{\m\alpha(2)}}\\
    *!R{\dRW[0]^j/p^{\m d_1} \tnsr I^{\m\alpha(0)}} & *!L{\Fil^{1-\m d_1}\dRW[1]^j \tnsr I^{\m\alpha(1)}} \ar[ul]\\
                                                         & *!L{\dRW[0]^j \tnsr I^{\m\alpha(0)}} \ar[ul]
                                                               }}$}\\
    E_2 &=
    \dbox{$\vcenter{\xymatrix@R=0.2em@C=-6.2em{
  *!R{\dRW[-\m d_1]^j/p^{\m d_2} \tnsr I^{\m\alpha(1)}} & *!L{\Fil^{2-\m d_2}\dRW[2]^j \tnsr I^{\m\alpha(2)}}\\
                                          & *!L{\Fil^{1-(\m d_1+\m d_2)}\dRW[1-\m d_2]^j \tnsr I^{\m\alpha(2)}}\\
                                          & *!L{\dRW[-\m d_1]^j \tnsr I^{\m\alpha(1)}} \ar[uul]
                                    }}$}\\
    E_\infty&=
    \dbox{$\vcenter{\xymatrix@R=0.2em@C=4em{
            & *!L{\Fil^{2-\m d_2}\dRW[2]^j \tnsr I^{\m\alpha(2)}}\\
            & *!L{\Fil^{1-(\m d_1+\m d_2)}\dRW[1-\m d_2]^j \tnsr I^{\m\alpha(2)}}\\
            & *!L{\dRW[-(\m d_1+\m d_2)]^j \tnsr I^{\m\alpha(2)}}
      }}$}
  \end{align*}
  The extension is $\dRW[2]_S^j \tnsr I^{\m\alpha(2)} = \dRW[2]_S^j \tnsr p^{\m d_0 + \m d_1 + \m d_2}W$. The reader may find it instructive to run the spectral sequence for various values of $\m d_1,\,\m d_2$ and see which terms vanish.
\end{example}

\begin{figure}
  \label{fig:tr-cp2-perfd}
  \begin{center}
    \includegraphics[width=\textwidth]{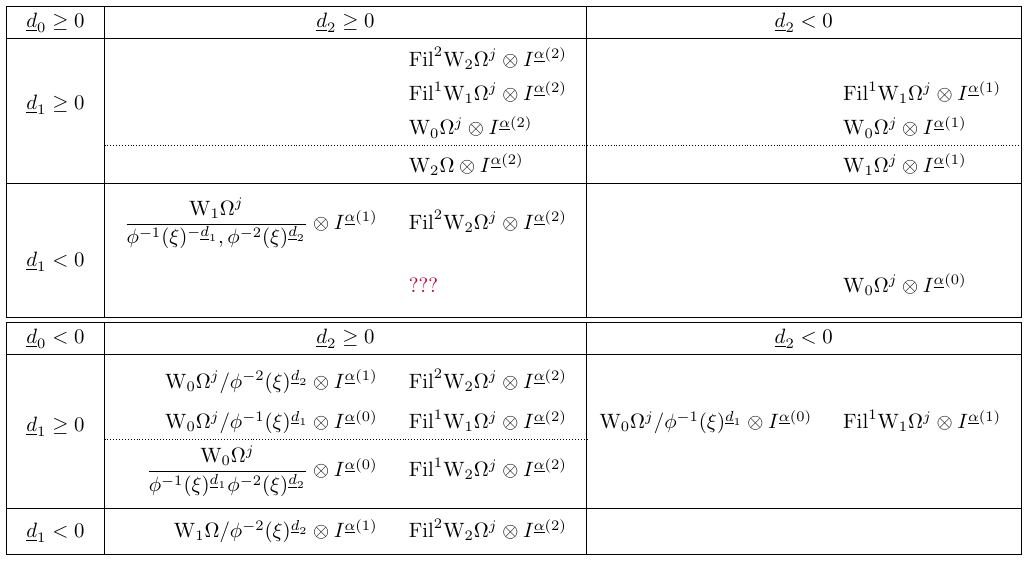}
    \caption{$E_\infty$ pages of the HOTRSS for $\TR^2_\rog(S\rel\SAinf{R^+})$. Extensions are marked below the dotted lines.}
  \end{center}
\end{figure}

\begin{figure}
  \label{fig:tr-cp2-crys}
  \begin{center}
    \includegraphics[width=\textwidth]{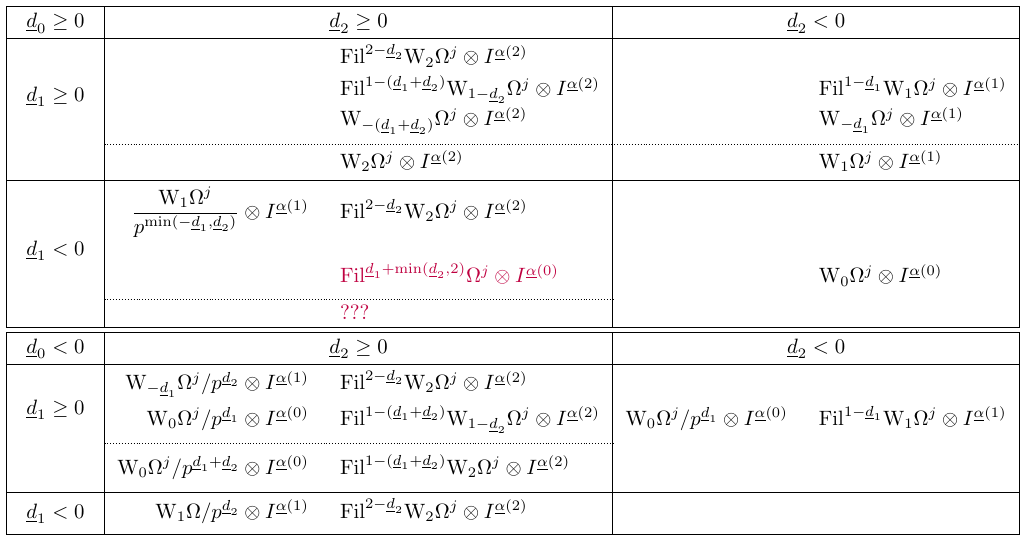}
    \caption{$E_\infty$ pages of the HOTRSS for $\TR^2_\rog(S\rel\SW k)$. Extensions are marked below the dotted lines.}
  \end{center}
\end{figure}

\begin{example}
  Now let $\m\alpha=(\m d_0,\m d_1,\m d_2)$ be any $\Z/p^2$-polyrepresentation. Using the same ideas as the previous examples, it is straightforward to evaluate the HOTRSS for $\m\alpha$, provided that $\m\alpha$ is contiguous. The results are listed in Figures \ref{fig:tr-cp2-crys} and \ref{fig:tr-cp2-perfd}. The exotic discontiguous parts are highlighted in purple.
\end{example}

\begin{theorem}[Contiguous range of $\TR_\rog$]
  \label{thm:tr-ctg}
  Let $\m\alpha = (\m d_0,\dotsc, \m d_n)$ be a $\Z/p^n$-polyrepresentation which is non-negatively concentrated in $[-\ell, -k]$.

  If $S$ is $p$-completely smooth over a $p$-torsionfree perfectoid ring $R^+$, then there are canonical identifications
  \begin{align*}
    \TR^n_{\m\alpha,j}(S\mid\SAinf{R^+})
      &= \Fil^k \dRW[\ell]^j_{S/R^+} \tnsr \bigotimes_{i=0}^\ell (\phi^{-i}(\xi)A)^{\tnsr\m d_i}\\
      \TR^n_{\m\alpha-1,j}(S\mid\SAinf{R^+})
      &= \frac{\dRW[k-1]^j_{S/R^+}}{\phi^{-k}(\xi)^{\m d_k}\dotsm\phi^{-\ell}(\xi)^{\m d_\ell}}\tnsr\bigotimes_{i=0}^{k-1} (\phi^{-i}(\xi)A)^{\tnsr\m d_i}
  \end{align*}

  If $S$ is smooth over a perfect $\F_p$-algebra $k$, then there are canonical identifications
  \begin{align*}
    \TR^n_{\m\alpha,j}(S\mid\SW k)
      &= \Fil^{k-(\m d_k+\dotsb+\m d_\ell)} \dRW[\ell]^j_S\tnsr \left(p^{\m d_0 + \dotsb + \m d_\ell}W\right)\\
    \TR^n_{\m\alpha-1,j}(S\mid\SW k)
      &= \frac{\dRW[k-1]^j_S}{p^{\m d_k+\dotsb + \m d_\ell}}\otimes \left(p^{\m d_0 + \dotsb + \m d_{(k-1)}}W\right)
  \end{align*}
\end{theorem}

\begin{proof}
  We begin with the easier case that $S$ is $p$-completely smooth over $R^+$.

  First we show that for $r\ge1$, the $E^r$ page of the HOTRSS is given by
  \begin{align*}
    E^r_{s,\m\alpha,j} &= \Fil^s\dRW[s]^j_{S/R^+} \tnsr I^{\m \alpha(\min(s+r-1,\ell))}\\
    E^r_{s,\m\alpha-1,j} &= \frac{\dRW[\max(s-r,k-1)]^j_{S/R^+}}{\phi^{-s}(\xi)^{\m d_s}} \tnsr I^{\m\alpha(s-1)}
  \end{align*}

  The $E^0$ page is given by
  \[
    E^0_{s,\m\alpha-*,j} =
    \left(\xymatrix@1{
        \ds\dRW[s-1]^j_{S/R^+} \tnsr I^{\m\alpha(s-1)}
        &
        \ar[l]
        \ds\dRW[s]^j_{S/R^+} \tnsr I^{\m\alpha(s)}
      }\right)
      \quad k\le s\le \ell
  \]
  from which we see that the $E^1$ page is as claimed. Inductively, on page $E^r$ there are differentials
  \[\xymatrix@C=-15em@R=1em{
      *!R{\fade{[s+r\le\ell]\left(\ds\frac{\dRW[s]^j_{S/R^+}}{\fade[100]{\phi^{-(s+r)}(\xi)^{\m d_{(s+r)}}}} \tnsr I^{\m\alpha(s+r-1)}
      \right)}}\\
      *!R{\left(\ds\frac{\dRW[s-r]^j_{S/R^+}}{\phi^{-s}(\xi)^{\m d_s}} \tnsr
          I^{\m\alpha(s-1)}
  \right)}
     & *!L{\ \ds\Fil^s \dRW[s]^j_{S/R^+} \tnsr
       I^{\m\alpha(s+r-1)}
   } \ar[ul]\\
     & *!L{\fade{\left(\ds\fade[100]{\Fil^{s-r}} \dRW[s-r]^j_{S/R^+} \tnsr
           I^{\m\alpha(s-1)}
\right)[k\le s-r]}}
     \ar[ul]
  }\]
  where we use the Iverson bracket to denote that the terms may be zero. Since $\m\alpha$ is contiguous, the terms stabilize as soon as the target/source of the relevant differential vanishes. From this we see that the $E^{r+1}$ page is as claimed.

  Now we have
  \begin{alignat*}2
    E^\infty_{s,\m\alpha,j} &= \Fil^s\dRW[s]^j_{S/R^+} \tnsr
    I^{\m\alpha(\ell)}
                            &\qquad& k\le s\le\ell\\
    E^\infty_{s,\m\alpha-1,j} &= \frac{\dRW[k-1]^j_{S/R^+}}{\phi^{-s}(\xi)^{\m d_s}} \tnsr
    I^{\m\alpha(s-1)}
                              && k \le s\le\ell
  \end{alignat*}
  from which the extensions are clear.

  Now we tackle the harder case that $S$ is smooth over $k$. Let $\hat r=\min(r,\ell-s+1)$ and let $\check r=\min(r,s-k+1)$. First we show that for $r\ge1$, the $E^r$ page of the HOTRSS is given by
  \begin{align*}
    E^r_{s,\m\alpha,j} &= \Fil^{s-\sum_{i=s}^{s+\hat r-1}\m d_i}\dRW[[s-\sum_{i=s+1}^{s+\hat r-1}\m d_i{]}]^j_S \tnsr 
   I^{\m\alpha(s+\hat r-1)}
    \\
    E^r_{s,\m\alpha-1,j} &= \frac{\dRW[[s-\check r-\sum_{i=s-\check r+1}^{s-1} \m d_i{]}]^j_S}{p^{\m d_s}} \tnsr
    I^{\m\alpha(s-1)}
  \end{align*}

  The $E^0$ page is given ($k\le s\le \ell$) by
  \begin{align*}
    E^0_{s,\m\alpha-*,j} &=
    \left[\xymatrix@1{
        \ds\dRW[s-1]^j_S \tnsr
        I^{\m\alpha(s-1)}
        &
        \ar[l]
        \ds\dRW[s]^j_S \tnsr
        I^{\m\alpha(s)}
      }\right]
      \\
      \intertext{from which we see that the $E^1$ page is}
      E^1_{s,\m\alpha-*,j} &=
        \left[\xymatrix@1{
                \ds\frac{\dRW[s-1]^j_S}{p^{\m d_s}} \tnsr
                I^{\m\alpha(s-1)}
                &
                \ds\Fil^{s-\m d_s}\dRW[s]^j_S \tnsr
                I^{\m\alpha(s)}
  }\right]
  \end{align*}
  as claimed. Inductively, on page $E^r$ there are differentials
  \small\[\hspace{-2em}\xymatrix@C=-22em@R=1em{
      *!R{\fade{[s+r\le\ell]\left(\ds\frac{\dRW[[s-\sum_{i=s+1}^{s+r-1}\m d_i{]}]^j_S}{\fade[100]{p^{\m d_{(s+r)}}}} \tnsr
            I^{\m\alpha(s+r-1)}
\right)}}\\
      *!R{\left(\ds\frac{\dRW[[s-r-\sum_{i=s-r+1}^{s-1}\m d_i{]}]^j_S}{p^{\m d_s}} \tnsr
          I^{\m\alpha(s-1)} 
  \right)}
     & *!L{\ \ds\Fil^{s - \sum_{i=s}^{s+r-1}\m d_i} \dRW[[s-\sum_{i=s+1}^{s+r-1}\m d_i{]}]^j_S \tnsr
       I^{\m\alpha(s+r-1)} 
   } \ar[ul]\\
     & *!L{\fade{\left(\ds\fade[100]{\Fil^{s-r-\sum_{i=s-r}^{s-1}\m d_i}} \dRW[[s-r-\sum_{i=s-r+1}^{s-1}\m d_i{]}]^j_S \tnsr
           I^{\m\alpha(s-1)} 
     \right)[k\le s-r]}}
     \ar[ul]
}\]\normalsize
where we use the Iverson bracket to denote that the terms may be zero. Since $\m\alpha$ is contiguous, the terms stabilize as soon as the target/source of the relevant differential vanishes. From this we see that the $E^{r+1}$ page is as claimed.

  Now we have
  \begin{alignat*}2
    E^\infty_{s,\m\alpha,j} &= \Fil^{s-\sum_{i=s}^\ell \m d_i}\dRW[[s-\sum_{i=s+1}^\ell \m d_i{]}]^j_S \tnsr
    I^{\m\alpha(\ell)}
                            &\qquad& k\le s\le\ell\\
    E^\infty_{s,\m\alpha-1,j} &= \frac{\dRW[[k-1-\sum_{i=k}^{s-1}\m d_i{]}]^j_S}{p^{\m d_s}} \tnsr
    I^{\m\alpha(s-1)}
                              && k \le s\le\ell
  \end{alignat*}
  from which the extensions are clear.
\end{proof}

\subsection{Special cases}
\label{sub:special}

Here we record some special cases of the theorem that are of interest for applications. In the positive and negative ranges, it is often more convenient to work with coefficients of irreducible decompositions rather than with dimension sequences. Note that there is some redundancy in the following statements since the integer grading $2i$ could be absorbed into $\m d_0$; we have phrased them this way since in applications one usually fixes $\m\alpha$ and considers $\TR_{2i\pm\m\alpha}$ for varying $i$.

\begin{corollary}[Negative range]
  Let $\m\alpha = \sum_{r=0}^n \m c_r \m\lambda_r$ be an actual poly-representation. Given $i\ge0$, let $\ell$ be such that
  \[ i - \m d_{\ell+1}(\m\alpha) < 0 \le i-\m d_\ell(\m\alpha). \]

  Then if $S$ is $p$-completely smooth over an arbitrary perfectoid ring $R$, we have
  \begin{align*}
    \TR^n_{2i-\m\alpha,j}(S;\Z_p)
      &= \dRW[\ell]^j_{S/R} \tnsr \bigotimes_{r=0}^\ell \phi^{-r}(\xi)^{i-\m d_r}A\\
      &= \dRW[\ell]^j_{S/R} \tnsr \xi_{\ell+1}^iA \tnsr \bigotimes_{r=0}^\ell \phi^{-r}(\xi_{\ell+1-r})^{-\m c_r}A\\
    \TR^n_{(2i-\m\alpha)-1,j}(S;\Z_p) &= 0
  \end{align*}
\end{corollary}
\begin{proof}
  When $R$ is either a perfect $\F_p$-algebra or $p$-torsionfree, we apply Theorem \ref{thm:tr-ctg}. For this, note that the polyrepresentation $2i-\m\alpha$ is non-negatively concentrated in $[-\ell,0]$, and $\m d_r(\alpha) = \sum_{s=0}^r \m c_r$. Since the result is the same in both cases, it must be true for arbitrary perfectoid rings using the fracture square of \cite[Proposition IV.3.2]{PrismNotes}.
\end{proof}

\begin{corollary}[Positive range]
  Let $\m\alpha = \sum_{r=0}^n \m c_r \m\lambda_r$ be an actual polyrepresentation. Given $i\in\Z$, let $k=k(i)$ be such that
  \[ \m d_{(k-1)}(\m\alpha) - i < 0 \le \m d_k(\m\alpha) - i. \]

  If $S$ is $p$-completely smooth over a $p$-torsionfree perfectoid ring $R^+$, then there are canonical identifications
  \begin{align*}
    \TR^n_{\m\alpha-2i,j}(S\rel\SAinf{R^+})
      &= \Fil^k \dRW[n]^j_{S/R^+} \tnsr \bigotimes_{r=0}^n \phi^{-r}(\xi)^{\m d_r-i}A\\
      &= \Fil^k \dRW[n]^j_{S/R^+} \tnsr (\xi_{n+1} A)^{-i} \tnsr \bigotimes_{r=0}^n \phi^{-r}(\xi_{n+1-r})^{\m c_r} A\\
      \TR^n_{(\m\alpha-2i)-1,j}(S\rel\SAinf{R^+})
      &= \frac{\dRW[k-1]^j_{S/R^+}}{\phi^{-k}(\xi)^{\m d_k-i}\dotsm\phi^{-n}(\xi)^{\m d_n-i}}\tnsr\bigotimes_{r=0}^{k-1} \phi^{-r}(\xi)^{\m d_r-i}A\\
      &= \frac{\dRW[k-1]^j_{S/R^+}}{\ds\frac1{\phi^{-k}(\xi_{n+1-k})^i}\prod_{r=k}^n \phi^{-r}(\xi_{n+1-r})^{\m c_r}} \tnsr (\xi_k A)^{-i} \tnsr \bigotimes_{r=0}^{k-1} \phi^{-r}(\xi_{k-r})^{\m c_r} A
  \end{align*}
  Using the first form, this also holds for $n=\infty$ if $R^+=\O_C$.

  If $S$ is smooth over a perfect $\F_p$-algebra $k$, then there are canonical identifications
  \begin{align*}
    \TR^n_{\m\alpha-2i,j}(S\rel\SW k)
      &= \Fil^{k - \left[\sum_{r=k}^n (\m d_r - i) \right]} \dRW[n]^j_S\tnsr p^{\left[\sum_{r=0}^n (\m d_r - i) \right]}W\\
      &= \Fil^{k - \left[\sum_{r=k}^n (n+1-r)\m c_r\right] + (n-k+1)i} \dRW[n]^j_S\tnsr p^{\left[\sum_{r=0}^n (n+1-r)\m c_r\right] - (n+1)i}W\\
    \TR^n_{\m\alpha-2i-1,j}(S\rel\SW k)
      &= \frac{\dRW[k-1]^j_S}{p^{\left[\sum_{r=k}^n (\m d_r-i)\right]}}\tnsr p^{\left[\sum_{r=0}^{k-1} (\m d_r-i)\right]}W\\
      &= \frac{\dRW[k-1]^j_S}{p^{\left[\sum_{r=k}^n (n+1-r)\m c_r\right]-(n+1-k)i}}\tnsr p^{\left[\sum_{r=0}^{k-1}(k-r)\m c_r\right]-ki}W
  \end{align*}
  This also holds for $n=\infty$ using the first form.
\end{corollary}

\begin{corollary}
  Let $\m\alpha$ be an actual poly-representation with $\m d_0(\m\alpha)\ge i$. If $S$ is $p$-completely smooth over an arbitrary perfectoid ring $R$, then there are canonical identifications
  \begin{align*}
    \TR^n_{\m\alpha-2i,j}(S\rel\SAinf{R})
      &= \dRW[n]^j_{S/R} \tnsr \bigotimes_{r=0}^n \phi^{-r}(\xi)^{\m d_r-i}A\\
      &= \dRW[n]^j_{S/R} \tnsr (\xi_{n+1} A)^{-i } \tnsr \bigotimes_{r=0}^n\phi^{-r}(\xi_{n-r+1})^{\m c_r}A\\
      \TR^n_{\m\alpha-2i-1,j}(S\rel\SAinf{R})
      &= 0
  \end{align*}
  Using the first form, this also holds for $n=\infty$ if $n<\infty$ or if $R\in\{k,\O_C\}$.
\end{corollary}

Here is the only other infinite-level computation where we know how to control the limit.

\begin{corollary}[Eventually constant contiguous range]
  Let $\m\alpha$ be a $\Z_p$-polyrepresentation non-negatively concentrated in $(-\infty,-k]$, and suppose that $\m d_r=\m d_\ell$ for all $r\ge\ell$.

  If $S$ is $p$-completely smooth over $\O_C$, then
  \begin{align*}
    \TR_{\m\alpha,j}(S\rel\SAinf{\O_C})
      &= \Fil^k \dRW^j_{S/\O_C} \tnsr \left(\bigotimes_{i=0}^{\ell-1} \phi^{-i}(\xi)^{\m d_i}A\right) \tnsr \left(\phi^{-\ell}(\mu)^{\m d_\ell}A\right)\\
      \TR_{\m\alpha-1,j}(S\rel\SAinf{\O_C})
      &= \frac{\dRW[k-1]^j_{S/\O_C}}{\phi^{-\ell}(\mu)^{\m d_\ell}\prod_{i=0}^{\ell-1}\phi^{-i}(\xi)^{\m d_i}} \tnsr\bigotimes_{i=0}^{k-1} \phi^{-i}(\xi)^{\m d_i}A
  \end{align*}

  If $S$ is smooth over $k$,
  \begin{align*}
    \TR_{\m\alpha,j}(S\mid\SW k)
      &= \begin{cases}
        \Fil^{k-(\m d_k+\dotsb+\m d_{\ell-1})} \dRW^j_S\tnsr \left(p^{\m d_0 + \dotsb + \m d_{\ell-1}}W\right) & \m d_\ell = 0\\
        0 & \m d_\ell > 0
      \end{cases}\\
    \TR_{\m\alpha-1,j}(S\mid\SW k)
      &= \left\{\!\!\!
        \begin{array}{r@{}ll}
          \dfrac{\dRW[k-1]^j_S}{p^{\m d_k+\dotsb + \m d_{\ell-1}}} &{}\tnsr \left(p^{\m d_0 + \dotsb + \m d_{(k-1)}}W\right) & \m d_\ell = 0\\
          \dRW[k-1]^j_S &{}\tnsr \left(p^{\m d_0 + \dotsb + \m d_{(k-1)}}W\right) & \m d_\ell >0
        \end{array}
        \right.
  \end{align*}
\end{corollary}

%% file: norms.tex

\section{Norms}
\label{sec:norms}

In this section we study the Norm map in $\RO(G)$-graded degrees. We briefly study $p$-primary Adams operations in \sec\ref{sub:adams}. We determine Norms for the de Rham-Witt complex in \sec\ref{sub:dRW-norm}.

\subsection{\tops{$p$}{p}-primary Adams operations}
\label{sub:adams}

The identity $FNf = f^p$ nearly determines the norm in ghost coordinates; to fully determine it, we need to understand the composite $R^n N$. This will turn out to be the $p^\th$ Adams operation $\psi^p$, so in this subsection we recall this operation and prove some facts about it.

\begin{definition}
  For an integer $\ell$, the natural map
  \[ \psi^\ell\colon \THH(R) = R^{\tnsr\T} \to R^{\tnsr\T/C_\ell} = \THH(R)_{\tnsr C_\ell} \]
  is called the \emph{$\ell^\th$ Adams operation}.
\end{definition}

Beware that $\THH(R)$ is equivalent to $\THH(R)_{\tnsr C_\ell}$ in $\Sp$, but not in $\Sp^{B\T}$. Therefore, $\psi^\ell$ may be regarded as an endomorphism of $\THH$, but it does not necessarily pass to $\TC^-$ and $\TP$. However, $\THH(R)$ and $\THH(R)_{\tnsr C_\ell}$ do become equivalent in $\Sp[1/\ell]^{B\T}$, so that in this case $\psi^\ell$ may be regarded as an automorphism of $\TC^-$ and $\TP$. This holds in particular if we $p$-complete at a prime $p$ not dividing $\ell$ (in which case $\ell$ may be generalized to any element of $\Z_p^\times$).

  For this reason, the analysis in the literature has mostly focused on the case $p\nmid\ell$. For example, in this case $\psi^\ell$ acts as multiplication by $\ell^i$ on $\gr^i$ of $\THH(-;\Z_p)$, $\TCmin(-;\Z_p)$, and $\TP(-;\Z_p)$ by \cite[Proposition 9.14]{BMS2} (see also \cite[Theorem 10.5]{ABGHLM}). However, we will need to study $\psi^p$ in the $p$-complete setting.

\begin{remark}[Lurie, c.f.\ the paragraph preceding {\cite[\sec IV.3]{NikolausScholze}}]
  The composite $RN$ is a map of $\T$-equivariant $\E_\infty$-rings, so it factors through the $C_p$ orbits of $\THH$ in the category of $\E_\infty$-rings, which we denote by $\THH_{\tnsr hC_p}$. The induced map $\THH_{\tnsr hC_p} \to \THH^{\Phi C_p}$ is an equivalence. If we identify the domain with $\THH$ and further map the codomain to $\THH^{tC_p}$, we obtain the Nikolaus-Scholze Frobenius.
  \[\xymatrix{
      \THH \ar[r]^-N \ar[d]_-{\psi^p} & \TR^1 \ar[r]^-R & \THH^{\Phi C_p} \ar[d]^-{\varphi_{\mathrm{NS}}}\\
      \THH_{\tnsr hC_p} \ar@{..>}[urr]_-\simeq \ar[rr]_-{\varphi_{\mathrm{NS}}} && \THH^{tC_p}
  }\]
  It follows that $RN\colon \THH \to \THH^{\Phi C_p} \simeq \THH_{\tnsr hC_p}$ is given by the $p$th Adams operation $\psi^p$.
\end{remark}


\begin{lemma}
  \label{lem:adams-perfd}
  Let $R$ be a perfectoid ring. Then $\psi^p$ acts as multiplication by $p^i$ on $\pi_{2i}\THH(R;\Z_p)$.
\end{lemma}
\begin{proof}
  The proof of \cite[Proposition 6.4.12]{RaksitDDR} goes through to show that $\psi^p$ acts as $p^i$ on $\gr^i_\HKR \HH_{\fil}(B/A)$ for any map $A\to B$ of derived commutative rings (the $\HC^-$ and $\HP$ parts of his proof do not go through). Applying this to $\HH(R;\Z_p)$, we learn that $\psi^p$ acts as $p^i$ on $\HH_{2i}(R;\Z_p)$. But the natural map $\THH_2(R;\Z_p) \to \HH_2(R;\Z_p)$ is an isomorphism by \cite[Proposition IV.4.2]{NikolausScholze}, so $\psi^p$ acts as $p^i$ on $\THH_{2i}(R;\Z_p) = \Sym^i \pi_2 \THH(R;\Z_p) = R\tnsr \xi^i A$.
\end{proof}

Let $S$ be $p$-completely smooth over any perfectoid ring $R$. In this case, recall from \cite[Corollary 6.9]{BMS2} (but essentially due to Hesselholt \cite{LarsPTyp}) that
\[ \pi_* \THH(S;\Z_p) = \Omega^*_{S/R} \tnsr_R \pi_*\THH(R;\Z_p). \]

\begin{lemma}
  \label{lem:adams-smooth}
  Let $S$ be $p$-completely smooth over any perfectoid ring $R$. Then $\psi^p$ acts as multiplication by $p^{i+j}$ on $\Omega^j_{S/R} \tnsr \xi^i A \subset \THH_{2i+j}(S;\Z_p)$.
\end{lemma}
\begin{proof}
  Given Lemma \ref{lem:adams-perfd}, we just need to show that $\psi^p$ acts as $p^j$ on $\Omega^j_{S/R}$, which follows by applying \cite[Proposition 6.4.12]{RaksitDDR} to $\HH(S/R;\Z_p)$.
\end{proof}

\subsection{de Rham-Witt norms}
\label{sub:dRW-norm}

In this subsection, we describe the Norm maps that the de Rham-Witt complex inherits from $\TR$.

First, we record how induction works.

\begin{proposition}
  In $p$-typical notation, induction of polyrepresentations is given by
  \[ \ind_n^{n+k}(\m d_0,\dotsc, \m d_n) = (\m d_0, p\m d_0, \dotsc, p^k \m d_0, p^k \m d_1, \dotsc, p^k \m d_n). \]
  Integrally, we have
  \[ \ind^{\Z/mn}_{\Z/n}(\lambda^j) = \sum_{i=0}^{m-1} \lambda^{in+j} \]
\end{proposition}
\begin{proof}
  These follow from $\res^G_H \ind^G_H \alpha = |G/H|\alpha$ and $(\ind^G_H\alpha)^G = \alpha^H$.
\end{proof}

It follows that there are Norm maps
\[ \dRW[n]^j_{S/R} \mono \TR^n_j(S;\Z_p) \morph^N \TR^{n+1}_{\ind j}(S;\Z_p) \hookleftarrow \dRW[n+1]^j_{S/R} \fade{{}\tnsr({}\cdots)} + \Fil^1\dRW[n+1]^{pj}_{S/R} \fade{{}\tnsr({}\cdots)}\]
where $({}\cdots)$ is some complicated invertible module which we ignore for the moment. There are other terms in $\TR^{n+1}_{\ind j}$, but only these two are relevant to the identities $FN\omega = \omega^p$ and $N \equiv \psi^p \bmod \Fil^1$.

As explained in the Introduction, the operation
\begin{align*}
  N(\omega)
    &\defeq (p-V(1))\frac{\psi^p(\omega)}p + V\frac{\omega^p}p\\
    &= \psi^p(\omega) - V\underbrace{\frac{F\psi^p\omega - \omega^p}p}_{{}\defeq \delta(\omega)}
\end{align*}
satisfies these identities. Beware that $V$ here is the $\RO(G)$-graded Verschiebung
\[ V\colon \TR^n_{pj}(S;\Z_p) \to \TR^{n+1}_{\ind(j)}(S;\Z_p). \]
Before identifying the twists $({}\cdots)$, let us explain why these identities suffice to pin down the Norm.

\begin{lemma}
  \label{lem:norm-ghost}
  Suppose $N\colon \TR^n_j(S;\Z_p) \to \TR^{n+1}_{\ind(j)}(S;\Z_p)$ is a natural transformation of functors
  \[ \{\text{\normalfont rings $p$-completely smooth over some perfectoid ring}\} \to \{\text{\normalfont sets}\}.\]
  satisfying $FN\omega = \omega^p$ and $N \equiv \psi^p \bmod \Fil^1$. Then $N$ coincides with the Hill-Hopkins-Ravenel norm. 
\end{lemma}
\begin{proof}
  Let $S$ be $p$-completely smooth over the perfectoid ring $R$. When $R$ is $p$-torsionfree, so is $S$, and thus the maps $F$ and $R^{n+1}$ are jointly injective, so the identities uniquely determine $N$.

  For a general perfectoid ring $R$, there exists a surjection $R^+ \epi R$ from a $p$-torsionfree perfectoid ring $R^+$ by \cite[Proposition 2.12]{CSPurity}. The $p$-completed cotangent complex $\L_{\Spf R/\Spf R^+} = (\L_{R/R^+})^\land_p$ vanishes by \cite[Lemma 3.14]{BMS1}, so by deformation theory $S$ descends to a $p$-completely smooth $R^+$-algebra $S^+$. By \cite[Proposition 10.14]{BMS1}, the natural map
  \[ \dRW[n]^*_{S^+/R^+} \mathbin{\hat\tnsr}_{\rmW_n(R^+)} \rmW_n(R) \to \dRW[n]^*_{S/R} \]
  is an isomorphism. Since $\rmW_n(-)$ preserves surjections, $\rmW_n(R^+) \to \rmW_n(R)$ is surjective, hence so is $\dRW[n]^*_{S^+/R^+} \to \dRW[n]^*_{S/R}$. This shows that the Norm for $S^+$ uniquely determines the Norm for $S$.
\end{proof}

In the following theorem, it is not obvious that the various forms of the codomain are equivalent; this will be explained in the proof. Recall the prismatic factorial notation from Definition \ref{dfn:prismatic-factorial}, and also that we write $n? \defeq (n-1)!$ for the Gamma function.

\begin{theorem}[Norms for the de Rham-Witt complex]
  \label{thm:norms-dRW}
  Let $R$ be any perfectoid ring, and $S$ a $p$-completely smooth $R$-algebra. When $p>2$ and $j\in\N$, there are Norm maps
  \begin{align*}
    N^k = N^{n+k}_n\colon
    \dRW[n]^{j}_{S/R}
    &\morph
    \sum_{r=0}^k
    \Fil^r \dRW[n+k]^{jp^r}_{S/R} \tnsr \left(\frac{\phi^{-k}\left(\xi_{n+1}^{h_1} (h_1)_\prism!\right)}{\xi_{n+k+1}^{(p^r-1)/2}}A\right)^{\tnsr j}\\
    &\hspace{.6em}=\fade{\sum_{r=0}^k
    \Fil^r \dRW[n+k]^{jp^r}_{S/R} \tnsr{}} \frac{\phi^{-k}\left(\xi_{n+1}^{h_j} (h_j)_\prism!\right)}{\xi_{n+k+1}^{(p^rj-1)/2}}A \qquad\text{\normalfont for $j$ odd}
    \intertext{where $h_m\defeq\frac{p^k m-1}2$, and $\sum$ is $\bigoplus$ for $j>0$. For even degrees, this can equivalently be written}
    \fade{N^k = N^{n+k}_n\colon}
    \dRW[n]^{2j}_{S/R}
    &\morph
    \sum_{r=0}^k
    \Fil^r \dRW[n+k]^{2jp^r}_{S/R} \tnsr \left(\frac{\phi^{-k}\left(\xi_{n+1}^{p^k-1} (p^k)_\prism?\right)}{\xi_{n+k+1}^{p^r-1}}A\right)^{\tnsr j}\\
    &\hspace{.6em}=\fade{
    \sum_{r=0}^k
  \Fil^r \dRW[n+k]^{2jp^r}_{S/R} \tnsr{}} \frac{\phi^{-k}\left(\xi_{n+1}^{p^kj-1} (p^k j)_\prism?\right)}{\xi_{n+k+1}^{p^rj-1}}A\\
  \intertext{and this form also works for $p=2$.}
  \intertext{When $R=k$ is a perfect $\F_p$-algebra and $j>0$, these reduce to an operation}
    N^k \colon
    \dRW[n]^{j}_S
    &\morph
    \Fil^k \dRW[n+k]^{jp^k}_S \tnsr \left(\frac{h_1!}{p^{kh_1}}W\right)^{\tnsr j}\\
    &\hspace{.6em}=\fade{
    \Fil^k \dRW[n+k]^{jp^k}_S \tnsr{}} \frac{h_j!}{p^{kh_j}}W \qquad \text{\normalfont for $j$ odd}\\
    \fade{N^k \colon}
    \dRW[n]^{2j}_S
    &\morph
    \Fil^k \dRW[n+k]^{2jp^k}_S \tnsr \left(\frac{p^k?}{p^{k(p^k-1)}}W\right)^{\tnsr j}\\
    &\hspace{.6em}=\fade{
    \Fil^k \dRW[n+k]^{2jp^k}_S \tnsr{}} \frac{(p^kj)?}{p^{k(p^kj-1)}}W\\
    \intertext{which also exists for $n=\infty$.}
    \intertext{When $R=\O_C$ is spherically complete, they assemble to an infinite Norm map}
    N^k \colon
    \dRW^{j}_{S/\O_C}
    &\morph
    \sum_{r=0}^k
    \Fil^r \dRW^{jp^r}_{S/\O_C} \tnsr \left(\frac{\phi^{-k}\left(\mu^{h_1} (h_1)_q!\right)}{\mu^{(p^r-1)/2}}\Ainf(\O_C)\right)^{\tnsr j}\\
    &\hspace{.6em}=\fade{\sum_{r=0}^k
    \Fil^r \dRW^{jp^r}_{S/\O_C} \tnsr{}} \frac{\phi^{-k}\left(\mu^{h_j} (h_j)_q!\right)}{\mu^{(p^rj-1)/2}}\Ainf(\O_C) \qquad\text{\normalfont for $j$ odd}\\
    \fade{N^k \colon}
    \dRW^{2j}_{S/\O_C}
    &\morph
    \sum_{r=0}^k
    \Fil^r \dRW^{2jp^r}_{S/\O_C} \tnsr \left(\frac{\phi^{-k}\left(\mu^{p^k-1} (p^k)_q?\right)}{\mu^{p^r-1}}\Ainf(\O_C)\right)^{\tnsr j}\\
    &\hspace{.6em}=\fade{
    \sum_{r=0}^k
  \Fil^r \dRW^{2jp^r}_{S/\O_C} \tnsr{}} \frac{\phi^{-k}\left(\mu^{p^kj-1} (p^k j)_q?\right)}{\mu^{p^rj-1}}\Ainf(\O_C)
  \end{align*}

  These maps satisfy the following properties:
  \begin{itemize}
    \item when $j=0$, $N(f)=f - V\delta(f)$ is the Norm map of Witt vectors \cite{AngeltveitNorm};

    \item when $j>0$, $N$ is given by 
      \begin{align*}
        N(\omega)
          &= (p-V(1))\frac{\psi^p(\omega)}p + V\frac{\omega^p}p
      \end{align*}
      or more generally
      \begin{align*}
        N^k(\omega)
          &= \sum_{r=0}^k V^r(p^{k-r} - V{\lf p^{k-r-1}\rf})\frac{\psi^{p^{k-r}}(\omega)^{p^{r}}}{p^k}
      \end{align*}

    \item when $R=k$ is a perfect $\F_p$-algebra, the preceding formula simplifies to
      \[ N^k(\omega) = V^k \frac{\omega^{p^k}}{p^k} \]
      for $|\omega|>0$, and in particular $N(\omega)=0$ whenever $\omega$ is of odd degree;

    \item $N$ satisfies the following identities for all $\omega,\,\tau$:
      \begin{align*}
        FN\omega &= \omega^p\\
        NV\omega &= p^{p-2}V^2\omega^p\\
        NR\omega &= RN\omega\\
        N(\omega\tau) &= N(\omega)N(\tau)\\
        N(\omega+\tau) &= N(\omega)+N(\tau)+Vs_1(\omega, \tau)\\
                       &\phantom{={}}\text{ where } (a+b)^p = a^p + b^p + ps_1(a, b)
      \end{align*}
  \end{itemize}
\end{theorem}

\begin{proof}
  We have already explained why the formula for $N$ must hold, and the final list of identities is true for any Tambara functor. It remains to identify the twists in the codomain of $N$. First suppose that $p$ is odd and $j=1$. We have
  \[ \ind^{n+k}_n(1_\R) = 1_\R + \lambda^{p^n} + \dotsb + \lambda^{h_1p^n} \]
  so that $\ind^{n+k}_n(1) - 1$ corresponds in prismatic notation to the twist by
  \begin{align*}
    (p^{n+1})_{\prism^{1/p}} \tnsr \dotsb \tnsr \left(h_1 p^{n+1}\right)_{\prism^{1/p}}
    &=
    \xi^{h_1} I_n^{h_1} \tnsr (h_1)_{\prism^{p^n}}!
  \end{align*}
  To switch to polygonic notation, we apply $\phi^{-(n+k)}$, turning this into
  \[ \phi^{-k}\left(\xi_{n+1}^{h_1} \tnsr (h_1)_\prism!\right) \]
  This is the twist for the $\dRW[n+k]^1_{S/R}$ term in the codomain of $N^k\colon \dRW[n]^1_{S/R} \to ({}\cdots)$. Each additional $\Fil^r \dRW[n+k]^{p^r}_{S/R}$ term moves $\frac{p^r-1}2\lambda^0$ from the topological part to the de Rham-Witt part, which has the effect of dividing by $\xi_{n+k+1}^{(p^r-1)/2}$.

 That is the analysis for $j=1$; for $j>1$, we can just raise everything to the $j^\th$ power, which yields the first form. But for $j$ odd, we could also write
 \[ \ind^{n+k}_n(j_\R) = 1_\R + \lambda^{p^n} + \dotsb + \lambda^{h_j p^n} \]
 which yields a twist of
 \[\xymatrix@C=4em{
     \xi^{h_j} I_n^{h_j} \tnsr (h_j)_{\prism^{p^n}}!
    \ar@{~>}[r]^-{\phi^{-(n+k)}} &
    \phi^{-k}\left(\xi_{n+1}^{h_j} \tnsr (h_j)_{\prism}!\right)
 }\]
 Then we need to subtract $\frac{p^rj-1}2 \lambda^0$ to land in $\Fil^r\dRW[n+k]^{jp^r}_{S/R}$, which yields the second form.

 Similarly, for even degrees $2j$ we can write
 \begin{align*}
   \ind^{n+k}_n(j\lambda^0)
    &= j\ind^{n+k}_n(\lambda^0)\\
    &= j\left(\lambda^0 + \lambda^{p^n} + \dotsb + \lambda^{(p^k-1)p^n}\right)
    \intertext{or}
    \ind^{n+k}_n(j\lambda^0)
    &= \lambda^0 + \lambda^{p^n} + \dotsb + \lambda^{(p^kj-1)p^n}
\end{align*}
which give the third and fourth forms respectively.

It is not obvious that the various forms of the twist are equivalent. In fact, they are not if the symbols are interpreted as $C_{p^\infty}$-Mackey functors. But they are equivalent when restricted to $C_{p^{n+k}}$, which is all that matters. For a simple example, $((\frac{5-1}2)!)^2=4$ is equal to $(5-1)!=24$ in $\Z/5$ but not in $\Z/25$; beware that $(2)_q$ is not equal to $(3)_q (4)_q$ in $\Z[q]/(5)_q$. Evidently many other formulations of the twist are possible.

  The claims about the Norm simplifying in characteristic $p$ are nontrivial, and proved in Corollary \ref{cor:norms-Fp} below. But we can identify the codomain of the $n=\infty$ Norms over $k$ and over $\O_C$. Over $\O_C$, we just replace every instance of $\xi_\bullet$ by $\mu$. In characteristic $p$, we have
  \[
    \frac{\phi^{-(k+1)}\left(\xi_n^{h_1} (h_1)_\prism!\right)}{\xi_{n+k+1}^{(p^r-1)/2}}
    =
    \frac{p^{nh_1} h_1!}{p^{(n+k+1)(p^r-1)/2}}
  \]
  and
  \[ n\frac{p^k-1}2 - (n+k+1)\frac{p^r-1}2 = n\frac{p^k - p^r}2 - (k+1)\frac{p^r-1}2 \]
  which goes to $\infty$ for $r<k$. Thus only the $r=k$ summand survives the limit in characteristic $p$.
\end{proof}

Next we want to justify the claim that $N(\omega)=V\frac{\omega^p}p$ in characteristic $p$. Intuitively this is true because $p-V(1)=0$ in $\rmW(k)$; however, we must be careful since some of this factor gets absorbed into the twist. Instead, as in Lemma \ref{lem:norm-ghost}, we must use $p$-torsionfree perfectoid rings to control the Norm in characteristic $p$.

For that matter, even in mixed characteristic, there is an implicit claim that $p-V(1)$ is divisible by the ideal twist. The resolution is a subtle property that has come up repeatedly when considering the Norm.

\begin{remark}[Refraction]
  \label{rmk:snell}
  There are two familiar instances of the optical phenomenon of \emph{refraction}\footnote{Previously the author referred to the first as diffraction and the second as refraction, but was informed that this is incorrect.}, with different corresponding mathematical incarnations. First, a prism breaks up white light into constituent wavelengths; according to Bhatt, the famous depiction of this on the cover of \cite{DarkSideMoon} was the inspiration for the name \emph{prismatic} cohomology, specifically due to the axiom $p \in I + \phi(I)$. In terms of Mackey data, the Tambara functor $\m{(A, I)}$ constructed in \cite{PrismsTambara1} has the feature that $\ker F^{n+k}_n$ and $\im V^{n+k}_n$ are \emph{annihilators} of each other; for a generic Green functor, all one can say is that $(\ker F^{n+k}_n)(\im V^{n+k}_n) = 0$ (by Frobenius reciprocity). We refer to this incarnation of refraction as \emph{Floyd refraction}, and say that the Green functor $\m{(A, I)}$ is \emph{Floyd refractive}.

  The second instance is that objects placed in water seem larger than they actually are; this distortion is described by Snell's law. Mathematically, this corresponds to being able to strengthen the prism condition to $p \in I^p + \phi(I)$; a generic Green functor always has $p - V(1) \in \ker F$, but $\res^{C_{p^\infty}}_{C_p} \m{(A, I)}$ obeys the scaled condition $p-V(1) \in (\ker F)^p$. We refer to this phenomenon as \emph{Snell refraction}.

  We highlighted this phenomenon in \cite{PrismsTambara1}, but there we did not know how to generalize Snell refraction beyond $C_p$. Mao discovered the correct generalization when using the norm to construct the (refined) prismatic logarithm \cite[Lemma 2.1]{MaoLogarithm}. His formulation is that for a transversal prism $(A, I)$,
  \[ \frac{\phi^n(I)}{I \cdot I^{p-1}_n} \to \frac{A}{I \cdot I_r^{p-1} } \]
  is a monomorphism with image $(p)$. Examining the proof, we find that in fact
  \[ p = \phi^n(\~\xi)\fade u + \~\xi \~\xi_n^{p-1}\fade v \]
  for units $u$ and $v$. Applying $\phi^{-(n+1)}$, we get
  \begin{equation}
    \label{eq:snell-polygonic}
    p = \xi \fade u + \phi^{-1}(\xi_n)^{p-1} \phi^{-n}(\xi) \fade v.
  \end{equation}
  This is the form we will use.
\end{remark}

\begin{corollary}
  \label{cor:norms-Fp}
  In characteristic $p$, for $|\omega|>0$ the Norm is given simply by $N^k(\omega)=V^k\frac{\omega^{p^k}}{p^k}$. Furthermore, $N(\omega)=0$ for $|\omega|$ odd.
\end{corollary}
\begin{proof}
  Let $R$ be a $p$-torsionfree perfectoid ring, and let $j>0$. First suppose that $p$ is odd, and let $h=\frac{p-1}2$ as above. Then the Norm goes
  \[
    N\colon \dRW[n]^j_{S/R} \to
    \dRW[n+1]^j_{S/R} \tnsr \phi^{-1}\left(\xi_{n+1}^h (h)_\prism!\right)^{\tnsr j}
    \fade{{}\oplus
    \Fil^1\dRW[n+1]^{pj}_{S/R} \tnsr ({}\cdots)}
  \]
  We can ignore the factor $\phi^{-1}((h)_\prism!)^j$ which is a unit. Now the expression ``$(p-V(1))\frac{\psi^p(\omega)}p$'' means an expression that goes to $0$ under $F$ and to $\psi^p(\omega)$ under $R^{n+1}$. Using Snell refraction \eqref{eq:snell-polygonic},
  \[
    \phi^{-1}(\xi_{n+1})^{hj} \phi^{-(n+1)}(\xi)^j \omega \tnsr \phi^{-1}(\xi_{n+1})^{hj} \approx p^j \omega \bmod \xi
  \]
  is (up to units) such an expression, since the $\phi^{-1}(\xi_{n+1})$ term is killed by $F$. Specializing to characteristic $p$, the coefficient of $\omega$ in the first tensor factor becomes
  \[ p^{(n+1)hj+j} \ge p^{(n+1)h+1} \ge p^{n+2} = 0 \in \rmW_{n+1}(k). \]

  When $p=2$ and $j$ is even, $hj$ still makes sense, and the above argument still goes through, except the last line becomes $p^{(n+1)hj+j} \ge p^{(n+1)+2} = 0 \in \rmW_{n+1}(k)$. We do not make any claims about odd-degree Norms when $p=2$. 

  To handle iterates of the Norm, one can either generalize the above argument, or just induct using the identity $NV\omega = p^{p-2} V^2 \omega^p$. To see that $N(\omega)=0$ for odd $|\omega|$, it suffices to show that $\omega \mapsto \frac{\omega^{p^k}}{p^k}$ is zero for odd $|\omega|$. This is easy to see by reduction to the $p$-torsionfree case.
\end{proof}